\documentclass[11pt]{amsart}

\usepackage[left=1.4in,top=1.3in,right=1.4in,nofoot]{geometry}
\usepackage{amssymb,amsxtra,eucal}
\usepackage{hyperref}
\usepackage[all,cmtip]{xy}
\usepackage{xcolor}

\def\varddots{\mathinner{\raise7pt\vbox{\kern3pt\hbox{.}}\mkern1mu\smash{\raise4pt\hbox{.}}
\mkern1mu\smash{\raise1pt\hbox{.}}}}

\newcommand{\p}{\mathpzc p}

\newcommand{\A}{\mathbb A}
\newcommand{\cA}{\mathcal A}
\renewcommand{\b}{\mathfrak b}
\newcommand{\B}{\mathbb B}

\newcommand{\C}{\mathcal C}
\newcommand{\E}{\mathcal E}
\newcommand{\cF}{\mathcal F}
\newcommand{\F}{\mathcal F}

\newcommand{\g}{\mathfrak g}

\renewcommand{\H}{\mathbb H}

\renewcommand{\L}{\mathbb L}
\newcommand{\cL}{\mathcal L}
\newcommand{\lb}{\overline L}

\newcommand{\m}{\mathfrak m}
\newcommand{\M}{\mathbb M}
\newcommand{\N}{\mathbb N}

\renewcommand{\O}{{\mathcal O}}

\renewcommand{\p}{\mathfrak p}

\newcommand{\Q}{\mathbb Q}

\newcommand{\QQ}{\mathbb Q}

\newcommand{\qp}{\QQ_p}

\newcommand{\qpb}{\overline\QQ_p}

\newcommand{\R}{\mathbb R}
\newcommand{\RR}{\mathcal R}
\newcommand{\cS}{\mathcal S}

\newcommand{\T}{\mathbb T}
\renewcommand{\t}{\mathfrak t}

\newcommand{\tv}{{\widetilde v}}
\newcommand{\tvprime}{{\widetilde v'}}
\newcommand{\X}{\mathbb X}
\newcommand{\Y}{\mathbb Y}
\newcommand{\vp}{\varphi}
\newcommand{\W}{\mathcal W}

\newcommand{\Z}{\mathbb Z}

\newcommand{\ve}{\varepsilon}

\renewcommand{\o}[1]{\overline{#1}}
\renewcommand{\u}[1]{\underline{#1}}
\newcommand{\wt}[1]{\widetilde{#1}}
\newcommand{\wh}[1]{\widehat{#1}}

\newcommand{\into}{\hookrightarrow}
\newcommand{\onto}{\twoheadrightarrow}
\newcommand{\congto}{\xrightarrow{\,\sim\,}}

\newcommand{\s}{^\times}

\newcommand{\ab}{^{\mathrm{ab}}}
\newcommand{\alg}{^{\mathrm{alg}}}

\newcommand{\dual}{^\vee}
\newcommand{\an}{^{\mathrm{an}}}
\newcommand{\kan}{^\text{\textup{$K$-an}}}
\newcommand{\qpan}{^\text{\textup{$\qp$-an}}}
\newcommand{\rig}{^{\mathrm{rig}}}
\newcommand{\cts}{^{\C^0}}
\newcommand{\ord}{^{\mathrm{ord}}}
\newcommand{\cris}{_{\mathrm{cris}}}
\newcommand{\lp}{^{\mathrm{lp}}}
\newcommand{\fs}{^{\mathrm{fs}}}

\makeatletter
\def\lan{\@ifnextchar[{\@with}{\@without}}
\def\@with[#1]{_{\text{\textup{$#1$-an}}}}
\def\@without{_{\mathrm{an}}}

\def\lrig{\@ifnextchar[{\@withr}{\@withoutr}}
\def\@withr[#1]{_{\text{\textup{$#1$-rig}}}}
\def\@withoutr{_{\mathrm{an}}}
\makeatother

\DeclareMathOperator{\Hom}{Hom}
\DeclareMathOperator{\Res}{Res}

\DeclareMathOperator{\GL}{GL}

\DeclareMathOperator{\WD}{WD}

\DeclareMathOperator{\diag}{diag}

\DeclareMathOperator{\Gal}{Gal}

\DeclareMathOperator{\coker}{coker}

\DeclareMathOperator{\Ind}{Ind}

\DeclareMathOperator{\soc}{soc}

\DeclareMathOperator{\im}{im}
\DeclareMathOperator{\Frob}{Frob}

\DeclareMathOperator{\supp}{supp}

\DeclareMathOperator{\id}{id}

\DeclareMathOperator{\LA}{LA}
\DeclareMathOperator{\PS}{PS}
\DeclareMathOperator{\Fil}{Fil}
\DeclareMathOperator{\Rep}{Rep}
\DeclareMathOperator{\Ref}{Ref}
\DeclareMathOperator{\colim}{colim}
\newcommand{\replac}{\Rep_{\text{\upshape{la.c}}}}
\newcommand{\repad}{\Rep_{\text{\upshape{ad}}}}
\newcommand{\plim}{\varprojlim}  
\newcommand{\ilim}{\varinjlim}

\newcommand{\whotimes}{\mathop{\widehat\otimes}\limits}

\newcommand{\sanmr}{\wh S\lan[\qp][\m_r]}

\newcommand{\pInd}{\Ind_P^G}
\newcommand{\bInd}{\Ind_B^G}

\newcommand{\Ws}{\W^{\soc}}
\newcommand{\Cs}{\C^{\soc}}

\renewcommand{\(}{\textup{(}}
\renewcommand{\)}{\textup{)}}

\theoremstyle{plain} 
\newtheorem{lm}[equation]{Lemma}

\newtheorem{prop}[equation]{Proposition}
\newtheorem{thm}[equation]{Theorem}
\newtheorem{coroll}[equation]{Corollary}
\newtheorem{conj}[equation]{Conjecture}
\theoremstyle{definition}
\newtheorem{df}[equation]{Definition}
\newtheorem{rk}[equation]{Remark}

\newtheorem{ex}[equation]{Example}
\newtheorem{hyp}[equation]{Hypothesis}
\newtheorem{cond}[equation]{Condition}

\numberwithin{equation}{section}
\numberwithin{figure}{section}

\begin{document}

\title{Towards the finite slope part for $\GL_n$}
\author{Christophe Breuil}
\address{L.M.O.\\ C.N.R.S.\\ Universit\'e Paris-Sud\\
Universit\'e Paris-Saclay\\ 91405 Orsay\\ France}
\email{christophe.breuil@math.u-psud.fr}
\author{Florian Herzig}
\address{Department of Mathematics\\ University of Toronto\\ 
Toronto, ON M5S 2E4\\ Canada}
\email{herzig@math.toronto.edu}
\thanks{The second author was partially supported by an NSERC grant and a Simons Fellowship (Simons Foundation grant \#504284).}

\maketitle

\begin{abstract}
Let $L$ be a finite extension of $\qp$ and $n\geq 2$. We associate to a crystabelline $n$-dimensional representation of $\Gal(\lb/L)$ satisfying mild genericity assumptions a finite length locally $\qp$-analytic representation of $\GL_n(L)$. In the crystalline case and in a global context, using the recent results on the locally analytic socle from \cite{BHS3} we prove that this representation indeed occurs in spaces of $p$-adic automorphic forms. We then use this latter result in the ordinary case to show that certain 
``ordinary'' $p$-adic Banach space representations constructed in our previous 
work appear in spaces of $p$-adic automorphic forms. This gives strong new 
evidence to our previous conjecture in the $p$-adic case.
\end{abstract}

\setcounter{tocdepth}{2}
\tableofcontents

\section{Introduction}\label{sec:introduction}

Let $\u G$ be a unitary group over a totally real number field $F^+$ which becomes $\GL_n$ over a totally imaginary quadratic extension $F$ of $F^+$ and such that $\u G(F^+\otimes _{\Q}\R)$ is compact. Then the $p$-adic Banach spaces of continuous functions $\C^0(\u G(F^+) \backslash \u G(\A_{F^+}^\infty),E)$ for $E$ a (varying) finite extension of $\qp$ can be seen as a $p$-adic analogue of the complex Hilbert space $L^2(\u G(F^+) \backslash \u G(\A_{F^+}))$. Assume that all places $v\vert p$ of $F^+$ split in $F$ and choose $\tv\vert v$ in $F$ for each $v\vert p$. Choose also a prime-to-$p$ compact open subgroup $U^p=\prod_{v\nmid p}U_v$ of $\u G(\A_{F^+}^{\infty,p})$, then the $U^p$-invariant vectors $ \wh S(U^p, E):=\{f:\u G(F^+) \backslash \u G(\A_{F^+}^\infty)/U^p \rightarrow E\}$ forms an admissible continuous representation of $\u G(F^+\otimes_{\Q}\qp)\cong \prod_{v\vert p}\GL_n(F_\tv)$ over $E$ (with $\u G(F^+\otimes_{\Q}\qp)$ acting by right translation of functions). Moreover to any absolutely irreducible automorphic Galois representation $r : \Gal(\o F/F) \to \GL_n(E)$ of level $U^p$ one can associate a non-zero invariant closed subspace $\wh S(U^p, E)[\m_r]$ of $\wh S(U^p, E)$.

The continuous representations $\wh S(U^p, E)[\m_r]$, as well as their locally $\qp$-analytic vectors $\wh S(U^p,E)\lan[\qp][\m_r]$ have attracted some attention over the past years (see e.g. \cite{Em1}, \cite{fondamental}, \cite{MR3529394}, \cite{BHS3} and all the references therein). With the notable exception of $\GL_n(F_\tv)=\GL_2(\qp)$ for all $v\vert p$ (see for instance \cite{CDP} and \cite{CEGGPS2}), and despite several partial results, these representations remain mysterious, e.g.\ one doesn't even know if they have finite length. In this article, we focus on the locally $\qp$-analytic representations $\wh S(U^p,E)\lan[\qp][\m_r]$ when $r$ is {\it crystalline} at all $p$-adic places and we use the recent results and techniques in \cite{BHS3} and \cite{breuil-ext1} to determine an explicit subrepresentation of $\wh S(U^p,E)\lan[\qp][\m_r]$, which is, to the knowledge of the authors, the largest known subrepresentation so far. Going beyond this subrepresentation will almost certainly require (seriously) new ideas.\\

Start with an arbitrary finite extension $L$ of $\qp$, and a crystalline representation $\rho : \Gal(\lb/L) \to \GL_n(E)$ (here $\Hom(L,E)$ has cardinality $[L:\qp]$) with distinct $\sigma$-Hodge--Tate weights for all embeddings $\sigma:L\hookrightarrow E$ and such that the eigenvalues $\varphi_1,\dots,\varphi_n$ of $\varphi^{[L_0:\qp]}$ on $D\cris(\rho)$ satisfy $\varphi_i \varphi_j^{-1} \not\in \{1, p^{\pm [L_0:\qp]}\}$ for all $i \ne j$, where $L_0$ is the maximal unramified extension contained in $L$. Then one can associate to $\rho$ a certain list of distinct irreducible constituents (the ``locally analytic socle'') $\Cs(\rho) = \{C(w\alg,\F)\}$ depending on two parameters: $w\alg$, which is a permutation of the $\sigma$-Hodge--Tate weights for each $\sigma:L\hookrightarrow E$, and $\F$ which is a refinement, i.e.\ an ordering of the eigenvalues $\varphi_i$. These two parameters $(w\alg,\F)$ satisfy a certain relation that involves the Hodge filtration on $D\cris(\rho)\otimes_{L_0}L$ (see \S\ref{sec:local-setup-crystab}). In fact $C(w\alg,\F)$ is the socle of a certain locally $\qp$-analytic principal series $\PS(w\alg,\F)$ of $\GL_n(L)$, see (\ref{PS}), and $C(1,\F)$ is the usual locally algebraic representation associated to $\rho$ by the classical local Langlands correspondence (it is in $\Cs(\rho)$ and doesn't actually depend on $\F$).

Then one defines a finite length admissible locally $\qp$-analytic representation $\Pi(\rho)\fs$ of $\GL_n(L)$ over $E$ (``fs'' for ``finite slope'') as follows (see Definition \ref{our-amalgam}):
\begin{enumerate}
\item For any $(w\alg, \F)$ such that $C(w\alg,\F)\in \Cs(\rho)$, we let $M(w\alg,\F)$ be the largest subrepresentation of $\PS(w\alg,\F)$ such that none of the irreducible constituents of $M(w\alg,\F)/C(w\alg,\F)$ is in $\Cs(\rho)$.
\item For any $C \in \Cs(\rho)$ we define the following amalgam over the common socle~$C$
\begin{equation*}
M(\rho)_C:= \bigoplus_{C} \{ M(w\alg,\F) : C(w\alg,\F) \cong C\}
\end{equation*}
and denote by $\Pi(\rho)\fs_C$ the unique quotient of $M(\rho)_C$ with socle isomorphic to~$C$.
\item We finally set $\Pi(\rho)\fs := \bigoplus_{C \in \Cs(\rho)} \Pi(\rho)\fs_C$.
\end{enumerate}

The representation $\Pi(\rho)\fs$ has socle $\bigoplus_{C \in \Cs(\rho)} C$ and in general is {\it not} multiplicity free. We give two explicit examples for $\GL_3(\qp)$ in \S\ref{sec:examples}. Note that $\Pi(\rho)\fs$ {\it does not} allow one to recover the Hodge filtration on $D\cris(\rho)\otimes_{L_0}L$, though it depends on it.

Now, for $v\vert p$ let $r_\tv := r|_{\Gal(\o {F_\tv}/F_\tv)}$ (where $r$ is automorphic of level $U^p$ as above) which we assume crystalline satisfying the above conditions for each $v\vert p$. One of us conjectured in \cite[Conj.\ 6.1]{socle2} that $\whotimes_{v|p}C(w\alg_\tv, \F_\tv)(\ve^{n-1})$ appears as a subrepresentation of $\wh S(U^p,E)\lan[\qp][\m_r]$ if and only if $C(w\alg_\tv, \F_\tv)\in \Cs(r_\tv)$ for each $v\vert p$ (see Conjecture \ref{conj:loc-an-socle}, here $(\ve^{n-1})$ is the twist by the $(n-1)$th power of the $p$-adic cyclotomic character $\ve$ on each factor $C(w\alg_\tv, \F_\tv)$). This conjecture was proven in \cite[Thm.~1.3]{BHS3} under the usual Taylor--Wiles assumptions when $U^p$ is sufficiently small and $r$ is residually absolutely irreducible (see Theorem \ref{BHS3}).

\begin{thm}[see Corollary \ref{cor:adjunct-max-subrep-ps-coho}]\label{intro}
Suppose that \cite[Conj.\ 6.1]{socle2} holds and for each $C = \wh \otimes_{v|p} C_\tv$ with $C_\tv \in \Cs(r_\tv)$ let
$$n_C := \dim_E \Hom_{\u G(F^+\otimes_{\Q}\qp)}\big(\whotimes_{v|p} C_\tv(\ve^{n-1}),\wh S(U^p,E)\lan[\qp][\m_r]\big) \in \Z_{> 0}.$$
Then there exists an injective $\u G(F^+\otimes_{\Q}\qp)$-equivariant linear map
\begin{equation}\label{fsembedding}
\bigoplus_{C = \whotimes C_\tv} \big(\whotimes_{v|p} \Pi(r_\tv)\fs_{C_\tv}(\ve^{n-1})\big)^{\oplus n_C} \hookrightarrow \wh S(U^p,E)\lan[\qp][\m_r].
\end{equation}
In particular there exists a $\u G(F^+\otimes_{\Q}\qp)$-equivariant injection of admissible locally $\qp$-analytic representations $\whotimes_{v|p} \Pi(r_\tv)\fs(\ve^{n-1}) \hookrightarrow \wh S(U^p,E)\lan[\qp][\m_r]$.
\end{thm}

In fact, Theorem \ref{intro} extends {\it verbatim} (and by the same proof) to the case when $r_{\tv}$ is crystabelline (and not just crystalline) for each $v\vert p$ (satisfying conditions analogous to the ones above). But in the crystalline case at least it becomes unconditional thanks to \cite[Thm.~1.3]{BHS3} (under the hypothesis of {\it loc.cit.}). Note that the embedding (\ref{fsembedding}) should be quite far from being an isomorphism in general, for instance because one can reasonably expect that the locally $\qp$-analytic representation $\wh S(U^p,E)\lan[\qp][\m_r]$ does determine all the Galois representations $r_\tv$ for $v\vert p$, and we know that this isn't the case for the left-hand side of (\ref{fsembedding}).

The proof of Theorem \ref{intro} proceeds as follows: (i) one first deduces from a general adjunction formula (Proposition \ref{prop:adjunct-max-subrep-ps-abstract} in the text) that  any equivariant homomorphism $\whotimes_{v|p} C_\tv(\ve^{n-1}) \hookrightarrow \wh S(U^p,E)\lan[\qp][\m_r]$ extends uniquely to an equivariant homomorphism $\whotimes_{v|p} M(r_\tv)_{C_\tv}(\ve^{n-1}) \rightarrow \wh S(U^p,E)\lan[\qp][\m_r]$, (ii) one proves (using \cite[Conj.\ 6.1]{socle2}) that any such homomorphism necessarily factors through the quotient $\whotimes_{v|p} \Pi(r_\tv)\fs_{C_\tv}(\ve^{n-1})$, and is injective since it is injective in restriction to the socle $\whotimes_{v|p} C_\tv(\ve^{n-1})$. The proof of Proposition \ref{prop:adjunct-max-subrep-ps-abstract} itself relies on the same strategy that was already used in the proof of \cite[Thm.~1.3]{breuil-ext1} (based on an extension of functional analysis results of Emerton).\\

Theorem \ref{intro} and Proposition \ref{prop:adjunct-max-subrep-ps-abstract} (on which it crucially relies) ``unify'' various results and give evidence to several conjectures in the literature, which is the main reason why we wrote this article. First Proposition \ref{prop:adjunct-max-subrep-ps-abstract} generalises an adjunction formula of Bergdall and Chojecki (\cite[Thm.~B]{BC-adjunction}) in the case of a Borel subgroup and a locally algebraic character. When $n=2$ the representation $\Pi(\rho)\fs$ is exactly the representation $\Pi(D_{\rm cris}(\rho))$ in \cite[\S4(9)]{socle1}. In this case Theorem \ref{intro} was already proven by Ding in the setting of the completed $H^1$ of unitary Shimura curves (see \cite[Thm.~6.3.7]{Ding}) by a different argument. Theorem \ref{intro} also gives evidence to \cite[Conj.~6.1.1]{breuil-ext1}, which implies in particular that each constituent of the form $C(w\alg_\tv, \F_\tv)(\ve^{n-1})$ (assuming there is only one place $v\vert p$ in $F^+$ for simplicity) which is {\it not} in $\Cs(r_\tv)$ and where $w\alg_\tv$ is a {\it simple} reflection appears in the socle of $\wh S(U^p,E)\lan[\qp][\m_r]/C(1, \F_\tv)(\ve^{n-1})$ (see the end of \S\ref{sec:examples}). Finally, Theorem \ref{intro} allows us to give strong evidence to a previous conjecture of the two authors in the ordinary case (\cite[Conj.~4.2.2]{fondamental}). Consider the crystalline representation $\rho$ above and assume moreover that $L=\qp$ and that $\rho$ is {\it upper triangular}. In that case we have a canonical refinement $\F_\rho$ and in \cite[\S3.3]{fondamental} we associated to $\rho$ a finite length continuous admissible representation $\Pi(\rho)\ord$ of $\GL_n(\qp)$ over $E$ of the form $\Pi(\rho)\ord=\bigoplus_{w}\Pi(\rho)_{C_\rho,w^{-1}}$, where $w$ runs over those $w=w\alg$ such that $C(w, w(\F_\rho))\in \Cs(\rho)$ and where each $\Pi(\rho)_{C_\rho,{w}^{-1}}$ is indecomposable and is a successive extension of certain unitary continuous principal series of $\GL_n(\qp)$ over $E$.

\begin{thm}[see Theorem \ref{thm:main}]\label{intrord}
In the setting of Theorem \ref{intro}, assume moreover that $p$ is totally split in $F^+$ and that each $r_\tv$ for $v\vert p$ is {\it upper triangular}. For each $w=(w_\tv)_{v\vert p}$ such that $C(w_\tv, w_\tv(\F_{r_\tv})) \in \Cs(r_\tv)$ let
$$n_w := \dim_E \Hom_{\u G(F^+\otimes_{\Q}\qp)}\big(\whotimes_{v|p} C(w_\tv, w_\tv(\F_{r_\tv}))(\ve^{n-1}),\wh S(U^p,E)\lan[\qp][\m_r]\big) \in \Z_{> 0}.$$
Then there exists an injective $\u G(F^+\otimes_{\Q}\qp)$-equivariant linear map
\begin{equation*}
\bigoplus_{w = (w_\tv)_v} \big(\whotimes_{v|p} \Pi(r_\tv)_{C_{r_\tv}, w_\tv^{-1}}(\ve^{n-1})\big)^{\oplus n_w} \into \wh S(U^p,E)[\m_r].
\end{equation*}
In particular there exists a $\u G(F^+\otimes_{\Q}\qp)$-equivariant injection of admissible continuous representations $\whotimes_{v|p} \Pi(r_\tv)\ord(\ve^{n-1}) \into \wh S(U^p,E)[\m_r]$.
\end{thm}

As above, recall that Theorem \ref{intrord} becomes unconditional under the assumptions of \cite[Thm.~1.3]{BHS3}. One way to prove Theorem \ref{intrord} goes as follows: (i) one first proves that the locally analytic vectors $(\Pi(\rho)\ord)\lan[\qp]$ of $\Pi(\rho)\ord$ is a (closed) subrepresentation of $\Pi(\rho)\fs$ (Proposition \ref{prop:Qrho-ordinary}), (ii) one proves that the universal unitary completion of $(\Pi(\rho)\ord)\lan[\qp]$ gives back $\Pi(\rho)\ord$ (Proposition \ref{prop:univ-completion-of-analytic-vecs}), (ii) one then combines these two results with Theorem \ref{intro} to deduce Theorem \ref{intrord}. (We actually give an alternative proof, see \S\ref{sec:global-app-ord}.)

Along the way, we carefully prove several unsurprising but useful technical results (e.g. Lemma \ref{lm:OS-product-groups} or Lemma \ref{lm:universal-tensor}), some of which having already been tacitly used in previous references (e.g.\ in \cite{BHS2}, \cite{BHS3}). We also provide a complete proof to the crucial Proposition \ref{prop:breuil-633} in the text which was already stated (but without a complete proof) in \cite[Prop.~6.3.3]{breuil-ext1}.

\subsection{Notation}\label{sec:notation}

We let $E/\qp$ be a finite extension and $\O$ its ring of integers. The field $E$ will be the coefficient field for {\it all} representations and locally convex vector spaces, unless otherwise stated. In particular all (completed) tensor products of locally convex vector spaces will be over $E$. All locally analytic manifolds will be assumed to be paracompact.

If $V$ is a locally convex vector space then $V'_b$ is its continuous dual with the strong topology (\cite[\S9]{nfa}). A unit ball in a Banach space is any open bounded lattice (or equivalently the unit ball of some norm defining
the Banach topology). If $V$ is any BH-space (see \cite[Def.~1.1.1]{locallyanalytic-memoir}), then $\o V$ denotes the latent Banach space structure on $V$ (\cite[\S1.1]{locallyanalytic-memoir}).

If $V_i$ are locally convex vector spaces, we write $V_1 \otimes_\pi V_2$ for the tensor product equipped
with the projective topology (\cite[\S17]{nfa}). If the $V_i$ are of compact type or Fr\'echet, then this agrees with the inductive topology and we just write $V_1 \otimes V_2$ (\cite[Prop.\ 1.1.31]{locallyanalytic-memoir}).

If $L$ is a finite extension of $\qp$ we will tacitly identify (characters of) $L^\times$ with (characters of) $W_L\ab$ by local class field theory, normalised so that uniformisers correspond to geometric Frobenius elements, and where $W_L$ is the Weil group of $L$ and $W_L\ab$ its maximal abelian quotient. We let $|\cdot |_L$ be the normalised absolute value $x\in L^\times \mapsto p^{-[L:\qp]{\rm val}(x)}$, where val$(p)=1$. For instance, for $L=\qp$ the cyclotomic character $\ve$ is identified with the character $x \mapsto x|x|_{\qp}$ of $\qp\s$.

If $I$ is a finite set, we denote by $\# I$ its cardinality.

All other notation will be introduced in the course of the text.

\subsection{Acknowledgements}\label{sec:acknowledgements}

The second author thanks the Universities of Paris-Sud and Paris 6, where some of this work was carried out, for pleasant working conditions. We also thank the referee for helpful comments.

\section{Preliminaries}\label{sec:preliminaries}

We establish some results in non-archimedean functional analysis that we couldn't find in the literature.\\

In this section, $K$ will denote a subfield of $E$ containing $\qp$ (in particular, $K$ is a finite extension of $\qp$). For a locally $K$-analytic group $G$ let $\replac G$ denote the category of locally $K$-analytic representations of $G$ on locally convex vector spaces of compact type and $\repad G$ the full subcategory of admissible locally $K$-analytic representations (see \cite{MR1887640}, \cite{MR1990669}, \cite{locallyanalytic-memoir}). Recall that a continuous linear map $f : V \to W$ between locally convex vector spaces is \emph{strict} if the continuous bijection $V/\ker(f) \to \im(f)$ is a topological isomorphism. If $V$ is a locally convex vector space, we denote by $\wh{V}$ its Hausdorff completion (\cite[Prop.~7.5]{nfa}).

\begin{lm}\label{lm:tensor-products-plim}\
\begin{enumerate}
\item If $0 \to V' \xrightarrow f V \xrightarrow g V'' $ is a strict short left exact sequence of locally convex vector spaces, then the sequence $0 \to \wh{V'} \xrightarrow {\hat f} \wh V \xrightarrow {\hat g} \wh {V''}$ is strict exact.
\item If $W$ is a locally convex vector space the functor $W \otimes_\pi (-)$ \(resp.\ $W \whotimes_\pi (-)$\) is exact \(resp.\ left exact\), meaning that it sends strict short exact sequences \(resp.\ strict short left exact sequences\) of locally convex vector spaces to strict short exact sequences \(resp.\ strict short left exact sequences\).
\item Suppose we are given locally convex vector spaces $V_i$ \($i \in I$\) and $W_j$ \($j \in J$\). Then we have a natural isomorphism
\begin{equation*}
\big(\prod_{i \in I} V_i\big) \whotimes_\pi \big(\prod_{j \in J} W_j\big) \cong \prod_{(i,j) \in I\times J} V_i \whotimes_\pi W_j.
\end{equation*}
\end{enumerate}
\end{lm}
\begin{proof}
  (i) 
We first consider the case where $g$ is surjective.
On the level of vector spaces we have a left exact sequence with continuous maps because projective limits are left exact. The map $\hat f$ is strict by the correspondence between open lattices in a locally convex vector space and in its completion. The map $\hat g$ is strict because the universal property of completion shows that $\wh V/\im (\hat f) \to \wh{V''}$ is a completion map, and completion maps are strict (\cite[Prop.\ 7.5]{nfa}).

In the general case we factor $g$ as composition $V \onto \im(g) \into V''$ of a strict surjection and
a strict injection. It remains to note that the completion of a strict injection is a strict
injection (by the same argument as before), and that the composition of a strict map and a strict
injection is still strict.

(ii) Consider the sequence from (i) with $g$ surjective. Then $1 \otimes_\pi f$ is strict by \cite[Cor.\ 17.5]{nfa} and $1 \otimes_\pi g$ is strict by the definitions and by \cite[Cor.\ 17.5]{nfa} again. The statement for $\whotimes_\pi$ now follows from part (i).

(iii) It suffices to show that $(\prod_{i \in I} V_i) \whotimes_\pi W \cong \prod_{i \in I} (V_i \whotimes_\pi W)$ for any locally convex vector space $W$. Since both sides are complete and Hausdorff, it suffices to show that both spaces are canonically isomorphic after passing to $\Hom_{\O}\cts(-,M)$ for any $\O$-torsion module $M$ (with discrete topology), where $\Hom_{\O}\cts$ means the continuous $\O$-linear maps. We easily check
$$\Hom_{\O}\cts\big(\big(\prod_{i \in I} V_i\big) \whotimes_\pi W,M\big) \cong \Hom_{\O}\cts\big(\big(\prod_{i \in I} V_i\big) \otimes_\pi W,M\big).$$
Given any element $f$ of this space, it has to vanish on $(\prod_S \Lambda_i \times \prod_{I-S} V_i) \times \Lambda$ for some finite subset $S \subseteq I$ and some open lattices $\Lambda_i$ in $V_i$ ($i \in S$) and $\Lambda$ in $W$. By linearity, $f$ is zero on $ (\prod_{I-S} V_i) \otimes_\pi W$. By (ii) the function $f$ factors to give an element of $\Hom_{\O}\cts((\prod_{S} V_i) \otimes_\pi W,M)$. Hence
\begin{equation*}
  \Hom_{\O}\cts\big(\big(\prod_{i \in I} V_i\big) \otimes_\pi W,M\big) \cong \ilim_{S} \Hom_{\O}\cts\!\big(\big(\prod_{i \in S} V_i\big) \otimes_\pi W,M\big)
  \cong \bigoplus_{i \in I} \Hom_{\O}\cts(V_i \otimes_\pi W,M),
 \end{equation*}
where the direct limit is over finite subsets $S$ of $I$. Replacing for a moment $V_i$ by $V_i \whotimes_\pi W$ and $W$ by $E$ we also get
 \begin{equation*}
  \Hom_{\O}\cts\big(\prod_{i \in I} \big(V_i \whotimes_\pi W\big),M\big) \cong \bigoplus_{i \in I} \Hom_{\O}\cts(V_i \whotimes_\pi W,M).
 \end{equation*}
By combining all the above isomorphisms (and using again that $M$ is complete) we are done.
\end{proof}

We recall that a Hausdorff locally convex vector space is said to be \emph{hereditarily complete} if all its Hausdorff quotients are complete. This is true for Fr\'echet and compact type spaces, cf.\ \cite[Def.\ 1.1.39]{locallyanalytic-memoir}.

\begin{coroll}\label{cor:exactness}
The completion functor is exact on any short exact sequence of locally convex vector spaces whose middle term has the property that its completion is hereditarily complete. The completed tensor product is even exact in the full subcategories of Fr\'echet and compact type spaces.
\end{coroll}
\begin{proof}
The first part follows from Lemma~\ref{lm:tensor-products-plim}(i) and its proof. The second part then follows from Lemma~\ref{lm:tensor-products-plim}(ii), as both subcategories are stable under completed tensor products \cite[Prop.\ 1.1.32]{locallyanalytic-memoir}. (In the case of Fr\'echet spaces, see also \cite[Cor.\ 4.14]{schraen-gl3}.)
\end{proof}

\begin{lm}\label{lm:compl-tensor-finite-colimits}
Suppose that $\C$ is a finite category and $F$ a functor from $\C$ to $\repad G$. If $V$ is a compact type space, then we have a topological isomorphism
\[ (\colim F) \whotimes V \cong \colim (F \whotimes V). \]
\end{lm}
\begin{proof}
In the category $\repad G$ any homomorphism is strict. The functor $\repad G \to \replac G$ sending $W$ to $W \whotimes V$ commutes with finite direct sums and cokernels by Corollary~\ref{cor:exactness}, hence also with finite colimits.
\end{proof}

We will apply Lemma \ref{lm:compl-tensor-finite-colimits} in the case of colimits indexed by a partially ordered set.

\begin{coroll}\label{cor:kernel-complete-tens}
Suppose that $V$, $W$, $U$ are Hausdorff locally convex vector spaces such that $V$ and $V'_b$ are bornological, $W$ is hereditarily complete, and $U$ is complete. Then for any continuous linear map $f : W \to U$ we have $\ker(1 \whotimes_\pi f) \cong V \whotimes_\pi (\ker f)$, where $1$ denotes the identity of $V$.
\end{coroll}
\begin{proof}
By assumption, the map $f$ factors as $W \onto W/\ker f \to U$, where $X = W/\ker f$ is complete, the first map is a strict surjection and the second map $i$ is a continuous injection. By \cite[Prop.\ 1.1.26]{locallyanalytic-memoir}, the map $1 \whotimes_\pi i$ is injective. The claim then follows from Lemma~\ref{lm:tensor-products-plim}(ii) applied to $0\to \ker f\to W\to W/\ker f\to 0$. (Note that the proof of \cite[Prop.\ 1.1.26]{locallyanalytic-memoir} uses also that $V'_b$ is bornological, in applying \cite[Prop.\ 7.16]{nfa}.)
\end{proof}

If $X$ is a locally $K$-analytic manifold and $V$ a Hausdorff locally convex vector space, we denote by $\C\kan(X,V)$ the locally convex vector space of locally $K$-analytic functions from $X$ to $V$ (\cite[Satz~2.1.10]{feaux-diss}, \cite[\S2]{MR1887640}).

\begin{lm}\label{lm:product-spaces}
Suppose that $X_1$, $X_2$ are locally $K$-analytic manifolds and that $V_1$, $V_2$ are locally convex vector spaces of compact type. Then the natural map
\begin{align}
 \C\kan(X_1,V_1) \whotimes_\pi \C\kan(X_2,V_2) &\congto \C\kan(X_1 \times X_2,V_1 \whotimes V_2)\label{eq:13} \\
 f_1 \whotimes f_2 &\mapsto \big((x_1,x_2) \mapsto f_1(x_1) \whotimes f_2(x_2)\big)\notag
 \end{align}
is an isomorphism.
\end{lm}
\begin{proof}
Write $X_i = \coprod_{j \in J_i} X_{ij}$ as a disjoint union of compact open subsets. Then $\C\kan(X_i,V_i) \cong \prod_{j \in J_i} \C\kan(X_{ij},V_i)$ (see \cite[\S2]{MR1887640}) and similarly for $\C\kan(X_1 \times X_2,V_1 \whotimes V_2)$. By Lemma~\ref{lm:tensor-products-plim}(iii) we are thus reduced to the case where $X_1$, $X_2$ are compact. In this case the lemma follows from \cite[A.1, A.2]{ST-duality} and \cite[Prop.\ 2.1.28]{locallyanalytic-memoir} by checking on the dense set of functions $f_i(x_i) = \vp_i(x_i)v_i$ with $\vp_i \in \C\kan(X_i,E)$ and $v_i \in V_i$.
\end{proof}

Let $H,G$ be locally $K$-analytic groups such that $H$ is a cocompact closed subgroup of $G$. Let $V \in \replac H$, then we have an exact functor $(\Ind_H^G -)\kan : \replac H \to \replac G$. Explicitly (see \cite[\S4.1]{feaux-diss})
$$(\Ind_H^G V)\kan = \{f \in \C\kan(G,V) : f(hg) = hf(g)\; \forall h \in H, g \in G \}$$
with left action of $G$ by right translation of functions. Moreover, choosing a locally analytic splitting $s : H\backslash G \to G$ of the projection $G \to H\backslash G$ (\cite[Satz 4.1.1]{feaux-diss}), we obtain an isomorphism of locally $K$-analytic manifolds $H \times H\backslash G \congto G$, $(h,x) \mapsto h s(x)$ and hence an isomorphism of locally convex vector spaces (with \cite[Prop.\ 2.1.28]{locallyanalytic-memoir})
\begin{equation}\label{eq:11}
 (\Ind_H^G V)\kan \cong \C\kan(H\backslash G,V) \cong \C\kan(H\backslash G) \whotimes V.
\end{equation}
Then exactness follows from the last assertion in Corollary \ref{cor:exactness}.

\begin{lm}\label{lm:ind-product-groups}
Suppose that for $i = 1,2$ we are given locally $K$-analytic groups $G_i$ with closed and cocompact subgroups $H_i$, as well as locally analytic representation $V_i \in \replac H_i$. Then $V_1 \whotimes V_2 \in \replac (H_1 \times H_2)$ and we have a natural isomorphism
\begin{equation}\label{eq:12}
(\Ind_{H_1 \times H_2}^{G_1 \times G_2} V_1 \whotimes V_2)\kan \cong (\Ind_{H_1}^{G_1} V_1)\kan \whotimes (\Ind_{H_2}^{G_2} V_2)\kan
\end{equation}
in $\replac (G_1 \times G_2)$.
\end{lm}
\begin{proof}
By \cite[Prop.\ 3.6.18]{locallyanalytic-memoir} applied to the inflations of $V_i$ in $\replac (H_1 \times H_2)$ we deduce that $V_1 \whotimes V_2 \in \replac (H_1 \times H_2)$. Fix now sections $H_i\backslash G_i \to G_i$, which induce a section $(H_1\times H_2)\backslash (G_1\times G_2)\to G_1\times G_2$. From \eqref{eq:13} applied with $X_i:=H_i\backslash G_i$ and from \eqref{eq:11} we can identify the two sides of~\eqref{eq:12} and see that $f_1 \whotimes f_2$ for $f_i \in (\Ind_{H_i}^{G_i} V_i)\kan$ is identified with $(g_1,g_2) \mapsto f_1(g_1) \whotimes f_2(g_2)$ in $(\Ind_{H_1 \times H_2}^{G_1 \times G_2} V_1 \whotimes V_2)\kan$. In particular, the map \eqref{eq:12} is $G_1 \times G_2$-equivariant.
\end{proof}

We will now prove a compatibility of the construction of Orlik--Strauch \cite{OS}, as extended in \cite[\S2]{socle1} and
\cite[Rk.\ 5.1.2]{BHS3}, with respect to product groups.  We recall briefly that for $G$ a split connected reductive group
over $K$ and a standard parabolic subgroup $P \subseteq G$, Orlik--Strauch constructed admissible locally $K$-analytic
representations $\F_P^G(M,\pi)$, where $M$ is an object of $\O_{\alg}^{\p}$ (an algebraic analogue of the parabolic BGG
category $\O$) and $\pi$ an admissible smooth representation of the Levi quotient of $P$. They showed that $\F_P^G(\cdot,\cdot)$
is functorial and exact in each argument (contravariant in the first argument), and they established strong irreducibility properties.
The first author extended their constructions to the case where $G$ is the restriction of scalars from $K$ to $\qp$ of a split connected reductive
group over $K$,
or even a product of several such groups. In this case, which we consider in the following, note that $\F_P^G(M,\pi)$ is locally $\qp$-analytic.
We refer to \cite[\S2]{socle1} and \cite[Rk.\ 5.1.2]{BHS3} for details.

\begin{lm}\label{lm:OS-product-groups}
Suppose that for $i = 1,2$ we are given a locally $\qp$-analytic group $G_i$ which is of the form $\prod_{j=1}^{n_i} G_{ij}(K_{ij})$, where $K_{ij}/\qp$ is finite and $G_{ij}$ is a split connected reductive group over $K_i$ \(we assume $p>3$ as in \cite{OS} if at least one $G_{ij}$ has factors of type different from $A$\). Suppose for each $i$ that $P_i \subseteq G_i$ is a parabolic subgroup, that $M_i \in \O_{\alg}^{\p_i}$ \(where $\p_i$ is the Lie algebra of the locally-$\qp$-analytic group $P_i$ and we use the notation of \cite{OS}\), and that $\pi_{P_i}$ is an \(admissible\) smooth representation of finite length of the Levi quotient of $P_i$. Then \(see \cite{OS} or \cite[\S2]{socle1} for the notation\)
\begin{equation}\label{eq:17}
 \F_{P_1}^{G_1}(M_1,\pi_{P_1}) \whotimes \F_{P_2}^{G_2}(M_2,\pi_{P_2})\cong \F_{P_1 \times P_2}^{G_1 \times G_2}(M_1 \otimes M_2, \pi_{P_1} \otimes \pi_{P_2}).
\end{equation}
\end{lm}
\begin{proof}
We take $K = \qp$ for the purpose of this proof. We let $\g_i$ be the Lie algebra of $G_i$ and $U(\g_i)$, $U(\p_i)$ the enveloping algebras. We note that $M := M_1 \otimes M_2 \in \O_{\alg}^{\p_1 \times \p_2}$ and that $\pi := \pi_{P_1} \otimes \pi_{P_2}$ is an (admissible) smooth representation of finite length of the Levi quotient of $P_1 \times P_2$. Choose finite-dimensional subspaces $W_i \subseteq M_i$ that are $\p_i$-stable and generate $M_i$ as $U(\g_i)$-module. Then $W := W_1 \otimes W_2$ has the analogous properties for the product group. For any $\partial_i \in U(\g_i) \otimes_{U(\p_i)} W_i$ we get the following commutative diagram from Lemma~\ref{lm:product-spaces} and Lemma~\ref{lm:ind-product-groups}:
\begin{equation}\small\label{eq:18}
 \begin{gathered}
 \xymatrix{
 (\Ind_{P_1}^{G_1} W_1' \otimes \pi_{P_1})\qpan \whotimes (\Ind_{P_2}^{G_2} W_2'\otimes \pi_{P_2})\qpan \ar[r]^-\sim \ar@{^{(}->}[d] & (\Ind_{P_1 \times P_2}^{G_1 \times G_2} W' \otimes \pi)\qpan \ar@{^{(}->}[d] \\
 \C\qpan(G_1,W_1' \otimes \pi_{P_1}) \whotimes_\pi \C\qpan(G_2,W_2' \otimes
 \pi_{P_2})\ar[d]_{\partial_1 \whotimes\partial_2}\ar[r]^-{\sim} &
 \C\qpan(G_1 \times G_2,W' \otimes \pi)\ar[d]^{\partial_1 \otimes\partial_2} \\
 \C\qpan(G_1,\pi_{P_1}) \whotimes_\pi \C\qpan(G_2,\pi_{P_2})\ar[r]^-{\sim} & \C\qpan(G_1 \times
 G_2,\pi), }
 \end{gathered}
\end{equation}
where the vertical maps in the bottom square are as defined in \cite[\S4.4]{OS}. Let $\phi_i : U(\g_i) \otimes_{U(\p_i)} W_i \onto M_i$ denote the natural surjection, and note that the natural surjection $\phi : U(\g_1 \times \g_2) \otimes_{U(\p_1 \times \p_2)} W \onto M$ is identified with $\phi_1 \otimes \phi_2$. By considering the diagram~\eqref{eq:18} we deduce that the left-hand side of~\eqref{eq:17} is identified (inside the top left of~\eqref{eq:18}) with the simultaneous kernel of all $\partial_1 \whotimes\partial_2$ for $\partial_1 \whotimes\partial_2$ in $\ker \phi_1 \whotimes W_2$ and $W_1 \whotimes \ker \phi_2$, equivalently with the simultaneous kernel of all $\partial_1 \whotimes 1\in \ker \phi_1 \whotimes 1$ and $1 \whotimes \partial_2\in 1 \whotimes \ker \phi_2$ (where 1 denotes alternatively the identity map of $(\Ind_{P_2}^{G_2} W_2'\otimes \pi_{P_2})\qpan$ and $(\Ind_{P_1}^{G_1} W_1' \otimes \pi_{P_1})\qpan$). Since the $U(\g_i)$ are noetherian, we easily deduce the claim from Corollary~\ref{cor:kernel-complete-tens}.
\end{proof}

If $G$ is a locally $K$-analytic group and $\sigma$ a continuous representation of $G$ on a Banach space, we denote by $\sigma\lan[(G,K)]$ the subspace of $\sigma$ of locally $K$-analytic vectors for the action of $G$, which carries a natural locally convex topology finer than the subspace topology (\cite[Def.~3.5.3]{locallyanalytic-memoir}). If the group $G$ is clear from the context, we will write $\sigma\lan[K]$ instead of $\sigma\lan[(G,K)]$. If $X, Y$ are topological spaces we denote by $\C^0(X,Y)$ the space of continuous maps from $X$ to $Y$.

\begin{lm}\label{lm:an-vectors-in-parab-ind}
Suppose that $H,G$ are locally $K$-analytic groups such that $H$ is a closed subgroup of $G$. Assume that there exists a compact open subgroup $G_0$ of $G$ such that $G = H G_0$. If $\sigma$ is an admissible continuous representation of $H$ on a Banach space, then
\begin{equation*}
(\Ind_H^G \sigma\lan[(H,K)])\kan\congto \big((\Ind_H^G \sigma)\cts\big)\lan[(G,K)],
\end{equation*}
where $(\Ind_H^G \sigma)\cts$ is the Banach space $\{f \in \C^0(G,\sigma) : f(hg) = hf(g)\; \forall h \in H, g \in G \}$ 
with left action of $G$ by right translation of functions.
\end{lm}

A priori, $(\Ind_H^G \sigma)\cts$ is a closed subspace of $\C^0(G,\sigma)$, where the latter has the compact-open topology.
As the restriction map identifies $(\Ind_H^G \sigma)\cts$ with $(\Ind_{H_0}^{G_0} \sigma)\cts$, or by the identification
with $\C^0(H\backslash G,\sigma)$, its topology is Banach.

\begin{proof}
Note that the assumption is satisfied if $G$ is the group of $K$-points of a connected reductive group over $K$
and $H$ is a parabolic subgroup. Also note that $\sigma\lan[(H,K)]$ is of compact type, as $\sigma$ is admissible.

Recall from \cite[\S 2.1, \S 3.5]{locallyanalytic-memoir} that we have a continuous injection $\sigma\lan[K] \hookrightarrow \sigma$ and hence a continuous injection $\C\kan(G,\sigma\lan[K]) \hookrightarrow \C^0(G,\sigma)$. This induces a continuous injection $(\Ind_H^G \sigma\lan[K])\kan \hookrightarrow (\Ind_H^G \sigma)\cts$ of closed subspaces, which is clearly $G$-equivariant. By passing to locally analytic vectors we get a continuous injection $i : (\Ind_H^G \sigma\lan[K])\kan \!\hookrightarrow ((\Ind_H^G \sigma)\cts)\lan[K]$ by Prop.\ 2.1.30, Prop.\ 3.5.6 and Thm.\ 3.6.12 in \cite{locallyanalytic-memoir}.

We now show that $i$ is surjective. Let $H_0 := H \cap G_0$. Note that restriction maps identify $(\Ind_H^G \sigma)\cts$ with $(\Ind_{H_0}^{G_0} \sigma)\cts$ and $(\Ind_H^G\sigma\lan[K])\kan$ with $(\Ind_{H_0}^{G_0} \sigma\lan[K])\kan$ (equivariantly for the action of $G_0$). Therefore we may assume, without loss of generality, that $G$ is compact.

Suppose that $f \in ((\Ind_H^G \sigma)\cts)\lan[K]$. Choose an analytic open subgroup $L \subseteq G$ which is the $K$-points of an affinoid rigid analytic group variety $\L$ defined over $K$ (we use the notation of \cite[\S2.1]{locallyanalytic-memoir}) such that $f$ is $\L$-analytic, i.e.\ the orbit map $L \to (\Ind_H^G \sigma)\cts$ of $f$ is rigid analytic in the sense of \cite[Def.~2.1.9(ii)]{locallyanalytic-memoir}. In particular, for each $g \in G$, the map $o_g : L \to \sigma$, $l \mapsto f(gl)$ is rigid analytic. Now choose an analytic open subgroup $M \subseteq H$ which is the $K$-points of an affinoid rigid analytic group variety $\M$ over $K$ such that (i) $M \subseteq H \cap \bigcap_{g \in G} g L g^{-1}$ and (ii) for any $g \in G$, the map $M \to L$, $m \mapsto g^{-1} m g =: m^g$ is rigid analytic (in particular this applies to the inclusion $M \subseteq L$ by considering $g=1$). This is possible since $G$ is compact.

Define $\sigma\lrig[\M] :=\C\rig(\M,\sigma)^{\Delta_{1,2}(M)}$ as in \cite[\S 3.3]{locallyanalytic-memoir}. Fix coset representatives $g_1,\dots,g_n$ of $G/L$ and define maps $\wt f_i : g_i L \to \sigma\lrig[\M]$ by $\wt f_i(g_i l)(m) := f(mg_i l) = mf(g_i l)$ for $l \in L$, $m \in M$. Then $\wt f_i$ is rigid analytic, since $\C\rig(g_i \L,\C\rig(\M,\sigma)) \cong \C\rig(g_i \L \times \M, \sigma)$ and the function $(g_il,m) \mapsto f(m g_i l) = f(g_i m^{g_i} l) = o_{g_i}(m^{g_i} l)$ is rigid analytic. It follows that the $\wt f_i$ define a locally analytic function $\wt f:G \to \sigma\lrig[\M]$ whose composition with $j : \sigma\lrig[\M] \hookrightarrow \sigma$ is $f$. Note that $\wt f$ is $H$-equivariant, since this is true after applying the injection $j$. Then the composition of $\wt f$ with $\sigma\lrig[\M] \hookrightarrow \sigma\lan[K]$ is the desired preimage of $f$.

We deduce that $i : (\Ind_H^G \sigma\lan[K])\kan \rightarrow ((\Ind_H^G \sigma)\cts)\lan[K]$ is a continuous bijection. It is a topological isomorphism of LB-spaces by Prop.\ 3.5.6 and Thm.\ 1.1.17 in \cite{locallyanalytic-memoir}.
\end{proof}

\begin{lm}\label{lm:analytic-vecs-tensor}
Suppose that for $i = 1,2$ we are given locally $K$-analytic groups $G_i$, as well as admissible continuous representations $\sigma_i$ of $G_i$ on Banach spaces \(\cite[\S3]{MR1900706}\). Then $\sigma_1 \whotimes \sigma_2$ is an admissible continuous representation of $G_1 \times G_2$ on a Banach space, and we have a natural isomorphism
\begin{equation}\label{eq:25}
(\sigma_1 \whotimes \sigma_2)\lan[(G_1 \times G_2,K)] \cong (\sigma_1)\lan[(G_1,K)] \whotimes (\sigma_2)\lan[(G_2,K)]
\end{equation}
in $\replac (G_1 \times G_2)$.
\end{lm}
\begin{proof}
We simplify notation ``$K\textrm{-an}=(G,K)\textrm{-an}$'' as in the proof of Lemma \ref{lm:an-vectors-in-parab-ind}. The admissibility claim follows from~\cite[Lemma A.3]{fondamental} (by passing to compact open subgroups that act unitarily). Now suppose that $H_i \subseteq G_i$ are analytic open subgroups which are $K$-points of affinoid rigid analytic group varieties $\H_i$ over $K$, so $H_1 \times H_2$ is an analytic open subgroup of $G_1 \times G_2$ ($=K$-points of $\H_1\times_K\H_2$). Then, by definition (\cite[\S 3.3]{locallyanalytic-memoir})
 \begin{equation}\label{eq:14}
 (\sigma_1 \whotimes \sigma_2)\lrig[\H_1 \times \H_2] \cong 
 (\C\rig(\H_1, \sigma_1) \whotimes \C\rig(\H_2, \sigma_2))^{\Delta_{1,2}(H_1)\times \Delta_{1,2}(H_2)}.
\end{equation}
We note that $\Delta_{1,2}(H_i)$ acts continuously on $\C\rig(\H_i, \sigma_i)$, because it acts continuously on $\C^0(H_i,\sigma_i)$ \cite[Prop.\ 3.1.5, 3.2.10]{locallyanalytic-memoir} and therefore on $\C^0(H_i,\sigma_i)\lrig[\H_i] \cong \C\rig(\H_i, \sigma_i)$ (see the comment after Def.\ 3.3.1 and also Prop.\ 3.3.7 in \cite{locallyanalytic-memoir}). (Alternatively, note that $\Delta_{1,2}(H_i)$ acts continuously on $\C\rig(\H_i,E) \otimes \sigma_i$ and hence on its completion $\C\rig(\H_i, \sigma_i)$.) As the $H_i$ are topologically finitely generated, we easily deduce using Corollary \ref{cor:exactness} and Corollary \ref{cor:kernel-complete-tens} (note that the strong dual of a Banach space is still a Banach space, hence is bornological)
\begin{multline}
 (\C\rig(\H_1, \sigma_1) \whotimes \C\rig(\H_2, \sigma_2))^{\Delta_{1,2}(H_1)\times \Delta_{1,2}(H_2)} \\
 \cong \C\rig(\H_1, \sigma_1)^{\Delta_{1,2}(H_1)} \whotimes \C\rig(\H_2, \sigma_2)^{\Delta_{1,2}(H_2)}
 \cong (\sigma_1)\lrig[\H_1] \whotimes (\sigma_2)\lrig[\H_2].\label{eq:15}
\end{multline}
Finally we choose cofinal descending sequences $H_i^{(j)}$ ($j \ge 1$) of analytic open subgroups of $G_i$ ($=K$-points of affinoid group varieties $\H_i^{(j)}$ over $K$) as in the proof of \cite[Prop.\ 6.1.3]{locallyanalytic-memoir}, so that the transition maps in $(\sigma_i)\lan[K] \cong \ilim_{j \ge 1} (\sigma_i)\lrig[\H_i^{(j)}]$ are injective and compact. By using the isomorphisms \eqref{eq:14} and \eqref{eq:15} together with \cite[Prop.\ 1.1.32]{locallyanalytic-memoir} we deduce that~\eqref{eq:25} holds.
\end{proof}

We refer to \cite[\S 3]{MR1990669} and \cite[\S 1.2]{locallyanalytic-memoir} for (weak) Fr\'echet--Stein algebras and coadmissible modules.

\begin{lm}\label{lm:tensor-product-frechet-stein}
Suppose that $A$, $B$ are Fr\'echet--Stein \($E$-\)algebras. Then $A \whotimes B$ is a weak Fr\'echet--Stein algebra, and if $M$ is a coadmissible $A$-module and $N$ is a coadmissible $B$-module, then $M \whotimes N$ is a coadmissible $A \whotimes B$-module.
\end{lm}
\begin{proof}
We start by choosing Fr\'echet--Stein structures $A \cong \plim A_n$, $B \cong \plim B_n$ in the sense of~\cite[Def.\ 1.2.10]{locallyanalytic-memoir}. We note that $A \whotimes B \cong \plim (A_n \whotimes B_n)$ by \cite[Prop.\ 1.1.29]{locallyanalytic-memoir} (and cofinality) and claim that this expression gives a weak Fr\'echet--Stein structure on $A \whotimes B$ in the sense of \cite[Def.\ 1.2.6]{locallyanalytic-memoir}. The first two conditions of that definition are verified, as $A_n \whotimes B_n$ is a Banach algebra. To check the third condition, it remains to show that the natural map $A \whotimes B \to A_n \whotimes B_n$ has dense image. We will show more generally that if $V_i$, $V_i'$ ($i = 1,2$) are locally convex vector spaces and $f_i : V_i \to V_i'$ are continuous maps with dense image, then $f_1 \whotimes_\pi f_2 : V_1 \whotimes_\pi V_2 \to V_1' \whotimes_\pi V_2'$ has dense image. It suffices to show that $f_1 \otimes_\pi f_2 : V_1 \otimes_\pi V_2 \to V_1' \otimes_\pi V_2'$ has dense image. By factoring $f_1 \otimes_\pi f_2 = (f_1 \otimes_\pi 1)\circ (1 \otimes_\pi f_2)$ we moreover reduce to the case where $V_2 = V_2'$ and $f_2 = \id$ (as the composition of continuous maps with dense images has dense image). It is enough to show that for any open lattices $\Lambda_1'$ of $V_1'$ and $\Lambda_2$ of $V_2$ we have $\im(f_1)\otimes V_2 + \Lambda_1' \otimes \Lambda_2 = V_1' \otimes V_2$. Given $v_1' \otimes v_2 \in V_1' \otimes V_2$, choose $a \in E\s$ such that $a v_2 \in \Lambda_2$ and write $a^{-1} v_1' = f_1(v_1) + x_1$ with $v_1 \in V_1$, $x_1 \in \Lambda_1'$ (which is possible by assumption). Then $v_1' \otimes v_2 = f_1(v_1) \otimes a v_2 + x_1 \otimes a v_2
 \in \im(f_1)\otimes V_2 + \Lambda_1' \otimes \Lambda_2$, as required.
 
To check that $M \whotimes N$ is coadmissible, we write $M_n := A_n \otimes_A M(\cong A_n \whotimes_A M)$, $N_n := B_n \otimes_B N(\cong B_n \whotimes_B N)$ (which are Banach spaces by \cite[Cor.~3.1]{MR1990669}) and note that $M \whotimes N \cong \plim(M_n \whotimes N_n)$ by \cite[Prop.\ 1.1.29]{locallyanalytic-memoir}. By assumption, $M_n$ (resp.\ $N_n$) is a finitely generated (locally convex) topological module over $A_n$ (resp.\ $B_n$) in the sense of \cite[\S1.2]{locallyanalytic-memoir}. It follows from Corollary~\ref{cor:exactness} that $M_n \whotimes N_n$ is a finitely generated topological $A_n \whotimes B_n$-module. We have
\begin{equation}\label{eq:29}
(A_n \whotimes B_n) \whotimes_{A_{n+1} \whotimes B_{n+1}} (M_{n+1} \whotimes N_{n+1}) \cong (A_n \whotimes_{A_{n+1}} M_{n+1}) \whotimes (B_n \whotimes_{B_{n+1}} N_{n+1})
\end{equation}
as topological $A_n \whotimes B_n$-modules, because all spaces in question are Fr\'echet (even Banach), so using for instance Corollary~\ref{cor:exactness} one easily checks that either side represents the Hausdorff quotient of $A_n\whotimes M_{n+1} \whotimes B_n \whotimes N_{n+1}$ on which the two natural actions of $A_{n+1}$, $B_{n+1}$ agree (acting on $A_n$, $B_n$ on the right and $M_{n+1}$, $N_{n+1}$ on the left). Thus the module in~\eqref{eq:29} is isomorphic to $M_n \whotimes N_n$. This completes the proof of coadmissibility.
\end{proof}

\begin{lm}\label{lm:tensor-product-admissible}
Suppose that for $i = 1,2$ we are given locally $K$-analytic groups $G_i$, as well as admissible locally $K$-analytic representations $\sigma_i$ of $G_i$. Then the locally $K$-analytic representation $\sigma_1 \whotimes \sigma_2$ of $G_1 \times G_2$ is admissible.
\end{lm}
\begin{proof}
For a locally $K$-analytic group $G$ let $D\kan(G) := \C\kan(G,E)'_b$ denote the locally $K$-analytic distribution algebra of $G$. Let now $G := G_1 \times G_2$. We may assume without loss of generality that both $G_i$ are compact, and we need to show that $(\sigma_1 \whotimes \sigma_2)'_b \cong (\sigma_1)'_b \whotimes (\sigma_2)'_b$ is a coadmissible $D\kan(G)$-module (the isomorphism follows from \cite[Prop.\ 1.1.32(ii)]{locallyanalytic-memoir}). By the argument in Step 2 of the proof of \cite[Thm.\ 5.1]{MR1990669} we have a quotient map $D\qpan(G) \onto D\kan(G)$ and it follows that a locally $K$-analytic representation of $G$ is admissible if and only if it is admissible as locally $\qp$-analytic representation (see the proof of \cite[Prop.~3.7]{MR1990669}). Thus we are reduced to the case where moreover $K = \qp$.

From Lemma~\ref{lm:product-spaces} and \cite[Prop.\ 1.1.32(ii)]{locallyanalytic-memoir} we have an isomorphism of Fr\'echet--Stein algebras $D\qpan(G) \cong D\qpan(G_1) \whotimes D\qpan(G_2)$, identifying a Dirac distribution $\delta_{(g_1,g_2)}$ with
$\delta_{g_1} \whotimes\delta_{g_2}$. The result then follows from Lemma~\ref{lm:tensor-product-frechet-stein}.
\end{proof}

The following technical lemma will be needed later. We denote by $\C^0_G(X,Y)$ the space of continuous $G$-equivariant maps from $X$ to $Y$.

\begin{lm}\label{lm:socle-for-product-group}
Suppose that $G_i$ \($i = 1,2$\) are topological groups and that $\pi_i$ and $\Pi_i$ are locally convex vector spaces equipped with a topological action of $G_i$ such that the $\Pi_i$ are of compact type. If $\C^0_{G_1}(\pi_1,\Pi_1) = 0$, then $\C^0_{G_1 \times G_2}(\pi_1 \whotimes_\pi \pi_2,\Pi_1 \whotimes \Pi_2) = 0$.
\end{lm}
\begin{proof}
Suppose that $f \in \C^0_{G_1 \times G_2}(\pi_1 \whotimes_\pi \pi_2,\Pi_1 \whotimes \Pi_2)$. By assumption, for any $y \in \pi_2$ and any $\lambda \in (\Pi_2)'$ we have $((1 \whotimes \lambda)\circ f)(\pi_1 \whotimes y) = 0$, so also $((1 \whotimes \lambda)\circ f)(\pi_1 \whotimes_\pi \pi_2) = 0$. It thus suffices to show that if $x \in \Pi_1 \whotimes \Pi_2$ and $(1 \whotimes \lambda)(x) = 0$ for all $\lambda \in (\Pi_2)'$, then $x = 0$. By \cite[Cor.\ 18.8]{nfa} we have $\Pi_1 \whotimes \Pi_2 \cong \cL_b((\Pi_2)'_b,\Pi_1)$ (see {\it loc.cit.} for the notation), and it is easily checked that this isomorphism is compatible with evaluation at any $\lambda \in (\Pi_2)'$. The lemma follows. Note that the same argument shows in fact that $\C^0_{G_1 \times G_2}(\pi_1 \whotimes_\iota \pi_2,\Pi_1 \whotimes \Pi_2) = 0$, where $\pi_1 \otimes_\iota \pi_2$ is the tensor product equipped with the inductive topology (\cite[\S17]{nfa}).
\end{proof}

\section{Universal unitary completions}\label{sec:universal}

We compute the universal unitary completion of certain locally analytic parabolic inductions (Proposition \ref{prop:univ-unit-compl}).\\

We still denote by $K$ a subfield of $E$ containing $\qp$ and by $G$ the group of $K$-points of a connected reductive group over $K$. We refer to \cite[\S1]{MR2181093} for the definition of universal unitary completions.

\begin{prop}\label{prop:univ-unit-compl}
Suppose that $P$ is a parabolic subgroup of $G$ with Levi subgroup $M$ and that $\sigma$ is a locally analytic representation of $M$ of compact type satisfying the following assumptions:
\begin{enumerate}
\item \label{item:1} $\sigma$ admits a central character $\chi_\sigma$;
\item \label{item:2} there exists a BH-subspace $\sigma_0$ of $\sigma$ such that $\sigma = \sum_{m \in M} m \sigma_0$.
\end{enumerate}
Then $\sigma$ has a universal unitary completion $\wh\sigma$ and the locally analytic representation $\pi := (\pInd \sigma)\kan$ satisfies the same hypotheses as $\sigma$ \(for $G$ instead of $M$\). If moreover $\chi_\sigma$ is unitary, then the universal unitary completion $\wh\pi$ of $\pi$ is given by $\wh\pi = (\pInd \wh\sigma)\cts$, together with the evident canonical map $\pi \to \wh\pi$.
\end{prop}

Note that if $\chi_\sigma$ is non-unitary, then $\wh\sigma = 0$, whereas $\wh\pi$ may be non-zero. For an example, see the representation $\pi$ in the proof of Proposition~\ref{prop:univ-completion-ordinary} below. Proposition \ref{prop:univ-unit-compl} has the following immediate corollary.

\begin{coroll}\label{cor:univ-unit-compl}
Suppose that $P$ is a parabolic subgroup of $G$ with Levi subgroup $M$ and that $\sigma$ is a finite-dimensional locally analytic representation of $M$ such that $\sigma$ has a unitary central character. Then $(\pInd \sigma)\kan$ has universal unitary completion $(\pInd \wh\sigma)\cts$.
\end{coroll}

Note that in the situation of Corollary \ref{cor:univ-unit-compl}, $\wh\sigma$ is the largest unitary quotient of the Banach representation $\sigma$ and thus an isomorphism if and only if the $M$-action on $\sigma$ is unitary. For example, if $G$ is quasi-split, $B = TU$ is a Borel subgroup, and $\chi$ a unitary character of $T$, then the universal unitary completion of ${(\bInd \chi)\kan}$ is $(\bInd \chi)\cts$.\\

To prepare for the proof of Proposition \ref{prop:univ-unit-compl}, we first need some preliminary results.\\

Let $\o P$ denote the parabolic subgroup opposite to $P$ with common Levi subgroup $M$, and let $\o N$ denote its unipotent radical. Choose $z \in Z_M$, the centre of $M$. Choose an analytic open subgroup $\o N_0 \congto \o\N_0(K)$ of $\o N$ (where as usual $\o\N_0$ is an affinoid rigid analytic group over $K$) such that the function $\o N_0 \to \o N_0$, $\o n \mapsto z\o nz^{-1}$ is rigid analytic, i.e.\ lifts to $\o\N_0 \to \o\N_0$. 
Recall \cite[\S0.1]{locallyanalytic-memoir} that a chart of a locally $K$-analytic manifold $X$ consists of
an open subset $X' \subseteq X$, a $K$-affinoid closed ball $\X'$, together with a locally analytic
isomorphism $X' \congto \X'(K)$.
For any chart $\phi : X \congto \X(K)$ of $\o N$ and any $m \in M$, $\o n \in \o N$ we denote for short by $m \o n \X m^{-1}$ the chart $m \o n X m^{-1} \congto X \xrightarrow{\phi} \X(K)$, where the first map sends $m \o n x m^{-1}$ to $x$. Recall \cite[\S 2.1]{locallyanalytic-memoir}
that an analytic partition of a locally $K$-analytic manifold $X$ is a partition of $X$ into a disjoint union of charts.
The following lemma is clear.

\begin{lm}\label{lm:cofinal}
If $\bigcap_{n \ge 0} z^n \o N_0 z^{-n} = 1$, then for any fixed $m \ge 0$ the set $\{z^n \o\nu \o\N_0 z^{-n} : \o\nu \in z^{-n-m}\o N_0 z^{n+m}/\o N_0 \}_{n \ge 0}$ is cofinal among all analytic partitions of $z^{-m} \o N_0 z^m$.
\end{lm}

\begin{lm}\label{lm:universal-tensor}\
\begin{enumerate}
 \item
 Suppose for each $1 \le i \le r$ that $G_i$ is the group of $K$-points of a connected reductive group over $K$ and that $\sigma_i$ is a locally convex vector space equipped with a continuous action of $G_i$ such that $\sigma_i$ has a universal unitary completion $\wh \sigma_i$. Then the representation $\sigma_1 \otimes_\pi \cdots \otimes_\pi \sigma_r$ of $G_1 \times \cdots \times G_r$ has universal unitary completion $\wh\sigma_1 \whotimes \cdots \whotimes \wh\sigma_r$.
 \item
 Suppose that $H$ is the group of $K$-points of a connected reductive group over $K$ and that $\sigma$ is a locally convex vector space equipped with a continuous action of $H$ such that $\sigma$ has a universal unitary completion $\wh \sigma$. Then the \(usual\) completion of $\sigma$ has the same universal unitary completion.
\end{enumerate}
\end{lm}
\begin{proof}
(i) It is clear that $G := G_1 \times \cdots \times G_r$ acts continuously on the locally convex vector space $\sigma := \sigma_1 \otimes_\pi \cdots \otimes_\pi \sigma_r$. We first show that for any $G$-invariant open lattice $\Lambda$ in $\sigma$ there exist $G_i$-invariant open lattices $\Lambda_i$ in $\sigma_i$ such that $\Lambda_1 \otimes_{\O} \cdots \otimes_{\O}\Lambda_r \subseteq \Lambda$. By definition of the projective tensor product topology, there exist open lattices $\Lambda_i' \subseteq \sigma_i$ such that $\Lambda_1' \otimes_{\O} \cdots \otimes_{\O}\Lambda_r' \subseteq \Lambda$, so we can take $\Lambda_i := \sum_{g_i \in G_i} g_i \Lambda_i'$. Hence, if $\Lambda_i$ for each $i$ is a minimal (up to commensurability) $G_i$-stable open lattice in $\sigma_i$, then $\Lambda := \Lambda_1 \otimes_{\O} \cdots \otimes_{\O}\Lambda_r$ is a minimal (up to commensurability) $G$-stable open lattice in $\sigma$. Hence by \cite[Lemma~1.3]{MR2181093} $\sigma$ has universal unitary completion $(\sigma_1 \otimes \cdots \otimes \sigma_r)_\Lambda\sphat$ (with the notation in \cite[\S19]{nfa}), which is isomorphic to $\wh\sigma_1 \whotimes \cdots \whotimes \wh\sigma_r$ by \cite[Lemma 19.10(ii)]{nfa}.

(ii) Let us write here $\sigma^c$ for the usual completion of $\sigma$. Since $H$ is locally compact it follows that $H$ acts continuously on $\sigma^c$ (use that a compact open subgroup of $H$ acts equicontinuously on $\sigma$ \cite[\S 3.1]{locallyanalytic-memoir}). If $\Lambda$ is a minimal (up to commensurability) $H$-stable open lattice in $\sigma$, then its closure $\o\Lambda$ is easily checked to be a minimal (up to commensurability) $H$-stable open lattice in $\sigma^c$. Finally, observe that the natural map $\sigma_\Lambda\sphat \to (\sigma^c)_{\o\Lambda}\sphat$ (again with the notation in \cite[\S19]{nfa}) is an isomorphism.
\end{proof}

Suppose that $X$ is a locally $K$-analytic manifold and that $\sigma$ is a Hausdorff locally convex
vector space.  We denote by $\C^0_c(X,\sigma)$ the vector space of functions in $\C^0(X,\sigma)$
that have compact support in $X$ (recall the support is the closure of the non-vanishing locus) and
by $\C_c\kan(X,\sigma)$ the linear subspace of $\C\kan(X,\sigma)$ consisting of compactly supported
locally analytic functions, i.e.\ $\C_c\kan(X,\sigma) = \C\kan(X,\sigma) \cap \C^0_c(X,\sigma)$.  As
in \cite[\S1]{MR2392361} we give $\C_c\kan(X,\sigma)$ the locally convex inductive limit topology
according to the isomorphism
\begin{equation}\label{eq:31}
  \C_c\kan(X,\sigma) \cong \ilim_{\{X_i,V_i\}_{i \in I}} \bigoplus_{i \in I} \C\rig(\X_i,\o V_i),
\end{equation}
where $\{X_i \congto \X_i(K)\}_{i \in I}$ runs through all analytic partitions of $X$ and the $V_i$
run through all BH-subspaces of $\sigma$. The inclusion $\C_c\kan(X,\sigma) \into \C\kan(X,\sigma)$ is continuous.
Moreover, the natural map
\begin{equation}
  \ilim_V \C_c\kan(X, \o V) \to \C_c\kan(X,\sigma),\label{eq:32}
\end{equation}
where $V$ runs through all BH-subspaces of $\sigma$, is a continuous bijection. 
If $X$ is compact we have $\C_c\kan(X,\sigma) = \C\kan(X,\sigma)$ and
the map~(\ref{eq:32}) is a topological isomorphism.

Suppose that $\X$ is a reduced affinoid rigid analytic space over $K$. Let $\C\rig(\X,E)$ denote the Banach algebra of $E$-valued rigid analytic functions on $\X$, i.e.\ $\C\rig(\X,E)=A \otimes_{K} E$, where $A$ is the affinoid algebra of $\X$. Then $\C\rig(\X,E)$ is a reduced affinoid algebra over $E$ by \cite[Lemma 3.3.1]{conrad-rigid}. Since it is reduced, the usual supremum norm defines the affinoid Banach topology. Recall also that if $V$ is any Banach space (over $E$), then $\C\rig(\X,V)$ denotes the Banach space $\C\rig(\X,E) \whotimes V$.

\begin{df}\label{df:rigid-lattices}\
\begin{enumerate}
\item We denote by $\C\rig(\X,\O)$ the open subring of power-bounded functions in $\C\rig(\X,E)$ (equivalently, the functions of supremum norm at most 1).
\item If $V^0$ is a unit ball of a Banach space $V$, we denote by $\C\rig(\X,V^0)$ the closure of $\C\rig(\X,\O) \otimes_{\O} V^0$ inside $\C\rig(\X,V)$.
\end{enumerate}
\end{df}

Note that $\C\rig(\X,V^0)$ is a unit ball in $\C\rig(\X,V)$. For example, if $\X$ is the $n$-dimensional closed unit ball, we can think of $\C\rig(\X,V)$ as the Banach space of all power series $\sum_{\u i \in \N^n} \u T^{\u i} v_{\u i}$ with $\u T^{\u i}:=T_1^{i_1}\cdots T_n^{i_n}$ and $v_{\u i} \in V$ tending to zero in $V$ as $\u i \to \infty$, and then $\C\rig(\X,V^0)$ consists of all such power series with $v_{\u i} \in V^0$. If $X$ is a topological space we remark that $\C^0(X,V^0)$ is the space of functions in $\C^0(X,V)$ with image in $V^0$. From Definition \ref{df:rigid-lattices} we obtain the following lemma.

\begin{lm}\label{lm:rigid-lattices}\
\begin{enumerate}
\item The evaluation map $\C\rig(\X,V) \to \C^0(\X(K),V)$ restricts to $\C\rig(\X,V^0) \to \C^0(\X(K),V^0)$.
\item Suppose that $V \to W$ is any map of Banach spaces sending a unit ball $V^0 \subseteq V$ to a unit ball $W^0 \subseteq W$. Then the natural map $\C\rig(\X,V) \to \C\rig(\X,W)$ restricts to $\C\rig(\X,V^0) \to \C\rig(\X,W^0)$.
\item Suppose that $\X \to \Y$ is a map of reduced affinoid rigid analytic spaces over $K$. Then the natural map $\C\rig(\Y,V) \to \C\rig(\X,V)$ restricts to $\C\rig(\Y,V^0) \to \C\rig(\X,V^0)$.
\end{enumerate}
\end{lm}

Suppose now that $X$ is a locally $K$-analytic manifold and $V$ a Banach space. Using Lemma \ref{lm:rigid-lattices} we can make the following definition. We note from~(\ref{eq:31}) that 
$$\C_c\kan(X,V) \cong \ilim_{\{X_i\}_{i \in I}} \bigoplus_{i \in I} \C\rig(\X_i,V),$$
where the inductive limit runs over all analytic partitions $\{X_i\}_{i \in I}$ of $X$ and $\C_c\kan(X,V)$ is equipped with the locally convex inductive limit topology.

\begin{df}\label{df:analytic-lattices}
We let
$$\C_c\kan(X,V^0) := \ilim_{\{X_i\}_{i \in I}} \bigoplus_{i \in I} \C\rig(\X_i,V^0),$$
where the inductive limit runs over all analytic partitions $\{X_i\}_{i \in I}$ of $X$. If $X$ is compact we write $\C\kan(X,V^0)$ for $\C_c\kan(X,V^0)$.
\end{df}

If $X' \congto \X'(K)$ is any chart of $X$, then the natural map $\C\rig(\X',V) \to \C_c\an(X,V)$ is injective,
since $\X'(K)$ is Zariski-dense inside $\X'$.
From the definitions it follows that $\C_c\kan(X,V)=\sum_{X'} \C\rig(\X',V)$, where the sum runs over all charts $X' \congto \X'(K)$. Note also that $\C_c\kan(X,V^0)$ is an open lattice in $\C_c\kan(X,V)$ and that $\C_c\kan(X,V^0) = \sum_{X'} \C\rig(\X',V^0)$ inside $\C_c\kan(X,V)$. If $X = \coprod_{i\in I} X_i$ is a partition by open subsets, then the natural map $\bigoplus_{i\in I} \C_c\kan(X_i,V) \to \C_c\kan(X,V)$ is a topological isomorphism, identifying $\bigoplus_{i\in I} \C_c\kan(X_i,V^0)$ with $\C_c\kan(X,V^0)$. (There are continuous maps in both directions, using a cofinality argument.) \\

From \ Definition \ \ref{df:analytic-lattices} \ and \ Lemma \ \ref{lm:rigid-lattices} \ we \ obtain \ the \ following \ lemma, \ where $\C^0_c(X,V^0):=\C^0_c(X,V)\cap \C^0(X,V^0)$.

\begin{lm}\label{lm:analytic-lattices}\
\begin{enumerate}
\item The natural  injection $\C_c\kan(X,V) \hookrightarrow \C^0_c(X,V)$ restricts to $\C_c\kan(X,V^0) \hookrightarrow \C_c^0(X,V^0)$.
\item Suppose that $V \to W$ is any map of Banach spaces sending a unit ball $V^0 \subseteq V$ to a unit ball $W^0 \subseteq W$. Then the natural map $\C_c\kan(X,V) \to \C_c\kan(X,W)$ restricts to $\C_c\kan(X,V^0) \to \C_c\kan(X,W^0)$.
\item Suppose that $X \to Y$ is a proper map of locally $K$-analytic manifolds \(i.e., inverse \ images \ of \ compact \ sets \ are \ compact\). \ Then \ the \ natural \ map $\C_c\kan(Y,V) \to \C_c\kan(X,V)$ restricts to $\C_c\kan(Y,V^0) \to \C_c\kan(X,V^0)$.
\end{enumerate} 
\end{lm}

We next establish a compatibility of the lattices $\C\rig(\X,V^0)$ and $\C_c\kan(X,V^0)$ with respect to addition.

\begin{lm}\label{lm:additivity-rigid-lattices}
Suppose that $\X$ is a reduced affinoid rigid analytic space over $K$.
\begin{enumerate}
\item Suppose that $V \onto W$ is a continuous surjection of Banach spaces, sending a unit ball $V^0 \subseteq V$ onto a unit ball $W^0 \subseteq W$. Then the natural map $\C\rig(\X,V^0) \to \C\rig(\X,W^0)$ is surjective.
\item Suppose that $\sigma$ is a Hausdorff locally convex vector space and that $V_1$, $V_2$ are BH-subspaces of $\sigma$ with unit balls $\Lambda_i \subseteq \o{V_i}$. Then $\Lambda_1 + \Lambda_2$ is a unit ball in $\o{V_1+V_2}$ and we have $\C\rig(\X,\Lambda_1) + \C\rig(\X,\Lambda_2) = \C\rig(\X,\Lambda_1+\Lambda_2)$ inside $\C\rig(\X,\o{V_1+V_2})$.
\end{enumerate}
\end{lm}
\begin{proof}
(i) We equip $\C\rig(\X,E)$ (resp.\ $V$) with the gauge norm of the open lattice $\C\rig(\X,\O)$ (resp.\ $V^0$). Then $\C\rig(\X,V^0)$ is the unit ball in $\C\rig(\X,V) = \C\rig(\X,E) \whotimes V$ with respect to the tensor product norm. Since the gauge norm of $\C\rig(\X,E)$ by definition takes values in $|E|_E$, it follows that $\C\rig(\X,E)$ is isometric to the Banach space $c_0(I,E)$ of functions $I \to E$ that tend to zero, for some set $I$ \cite[I.1.5]{colmez-une-variable}. By \cite[I.1.8]{colmez-une-variable} we have isometric isomorphisms $\C\rig(\X,V) \cong c_0(I,E) \whotimes V \cong c_0(I,V)$, hence giving a isomorphism $\C\rig(\X,V^0) \cong c_0(I,V^0)$ of unit balls. The claim now follows since $c_0(I,V^0) \to c_0(I,W^0)$ is clearly surjective.

(ii) By \cite[Prop.~1.1.5]{locallyanalytic-memoir} note first that $V_1+V_2$ is also a BH-subspace of $\sigma$. By \cite[Prop.~1.1.2(ii)]{locallyanalytic-memoir} we have a continuous surjection of Banach spaces $\o{V_1}\oplus \o{V_2} \onto \o{V_1+V_2}$. It is an open map by the open mapping theorem and sends $V^0 := \Lambda_1 \oplus \Lambda_2$ surjectively onto $W^0 := \Lambda_1 + \Lambda_2$, which is thus a unit ball in $\o{V_1+V_2}$. We conclude by (i).
\end{proof}

\begin{coroll}\label{cor:additivity-analytic-lattices}
Suppose that $X$ is a locally $K$-analytic manifold. Suppose that $\sigma$ is a Hausdorff locally convex vector space and that $V_1$, $V_2$ are BH-subspaces of $\sigma$ with unit balls $\Lambda_i \subseteq \o{V_i}$. Then we have the following equalities
inside $\C_c\kan(X,\o{V_1+V_2})$:
\begin{align*}
 \C_c\kan(X,\Lambda_1) + \C_c\kan(X,\Lambda_2) &= \C_c\kan(X,\Lambda_1+\Lambda_2),\\
 \C_c\kan(X,\o{V_1}) + \C_c\kan(X,\o{V_2}) &= \C_c\kan(X,\o{V_1+V_2}).
\end{align*}
\end{coroll}
\begin{proof}
The first equality follows by summing the equality in Lemma~\ref{lm:additivity-rigid-lattices}(ii) over all charts (or all analytic partitions) of $X$. The second follows by inverting $p$.
\end{proof}

\begin{prop}\label{prop:amice}
Suppose that $X$ is a locally $K$-analytic manifold and $V$ a Banach space with unit ball $V^0$. Then we have inside $\C_c^0(X,V)$:
\begin{equation}
 \C_c\kan(X,V) \cap \C_c^0(X,V^0) = \C_c\kan(X,V^0).\label{eq:1}
 \end{equation}
\end{prop}
\begin{proof}
It is clear that the right-hand side is contained in the left-hand side. To show that the other inclusion holds, choose any analytic partition $\{X_i\}_{i\in I}$ of $X$. Then each term in~(\ref{eq:1}) decomposes as a direct sum over $i$, and so we may assume without loss of generality that $X \cong \O_K^d$ for some $d \ge 1$.

We first make some preliminary remarks. Let $\varpi$ denote a uniformiser of $K$, $k$ the residue field of $K$, $q$ the cardinality of $k$, and $\varpi_E$ a uniformiser of $E$. We equip the Banach space $V$ with the gauge norm of $V^0$. Then any lift $(v_i)_{i \in I}$ of an algebraic basis of $V^0/\varpi_E V^0$ is an orthonormal basis of $V$. Note that if $X$, $Y$ are compact topological spaces, then we have isomorphisms of Banach spaces (see e.g. \cite[\S 17]{nfa})
$$\C^0(X,V) \cong \C^0(X,E) \whotimes V\ {\rm and}\ \C^0(X \times Y,E) \cong \C^0(X,E) \whotimes \C^0(Y,E).$$
It is easy to verify that they are all isometries, where each function space carries the supremum norm. On the other hand, for $r \in |K^\times|_K$ and $\u a \in K^n$ let $\B_r(\u a)$ denote the $K$-affinoid closed ball of radius $r$ and centre $\u a$. Then for $\u a \in K^n$, $\u b \in K^m$ the canonical isomorphism $\B_r(\u a) \times \B_r(\u b) \cong \B_r(\u a, \u b)$ induces an isometric isomorphism $\C\rig(\B_r(\u a, \u b),E) \cong \C\rig(\B_r(\u a),E) \whotimes \C\rig(\B_r(\u b),E)$. Define the Banach space $\LA_h(\O_K^d,E) := \bigoplus \C\rig(\B_{|\varpi|_K^h}(\u a), E)$, where the sum runs over representatives $\u a$ of $\O_K^d/\varpi^h \O_K^d$ and each summand is equipped with the supremum norm. It follows from what is before that $\LA_h(\O_K^d,E) \cong \wh\bigotimes_{1 \le i \le d} \LA_h(\O_K,E)$ is an isometric isomorphism. Define more generally the Banach space $\LA_h(\O_K^d,V) := \bigoplus \C\rig(\B_{|\varpi|_K^h}(\u a), V)$, where each summand is equipped with the tensor product norm. Hence its unit ball is $\LA_h(\O_K^d,V^0) := \bigoplus \C\rig(\B_{|\varpi|_K^h}(\u a), V^0)$, and we have a topological isomorphism
$$\C\kan(\O_K^d,V) \cong \ilim_h \LA_h(\O_K^d,V)\cong \ilim_h (\LA_h(\O_K^d,E)\whotimes V)$$
which restricts to $\C\kan(\O_K^d,V^0) \cong \ilim_h \LA_h(\O_K^d,V^0)$.

We now recall a particular Mahler basis of $\C^0(\O_K,E)$ from \cite[\S 1]{deshalit-mahler}. For $m \ge 1$ let $g_{q^m}(z) := \varpi^{-(q^m-1)/(q-1)} \prod_{r \in \RR_m} (z-r)$, where $\RR_m := \{ \sum_{i = 0}^{m-1} \wt{a_i} \varpi^i : a_i \in k \}$ and tilde denotes the Teichm\"uller lift. Any integer $n \ge 0$ can be written $n = \sum_{i=0}^{m-1} b_i q^i$ for some $m\geq 1$ and some $b_i\in \{0,\dots,q-1\}$, and we define $g_n := \prod_{i=0}^{m-1} g_{q^i}^{b_i}$ and, for $h \ge 0$, $g_{n,h} := \varpi^{\sum_{s > h} \lfloor \frac n{q^s}\rfloor} g_n$. Then the $\{g_n\}_{n \ge 0}$ form an orthonormal basis of $\C^0(\O_K,E)$, and for any $h \ge 0$, the $\{g_{n,h}\}_{n \ge 0}$ form an orthonormal basis of $\LA_h(\O_K,E)$ (see Thm.\ 1.1 and Prop.\ 4.2 in \cite{deshalit-mahler}).

Suppose now that $f \in \C\kan(\O_K^d,V) \cap \C^0(\O_K^d,V^0)$. Then $f \in \LA_h(\O_K^d,V)$ for some $h \ge 0$ and by what is above we can uniquely write
\begin{equation*}
 f(z_1,\dots,z_d) = \sum_{\u n, i} \lambda_{\u n, i} g_{n_1}(z_1) \cdots g_{n_d}(z_d) v_i,
\end{equation*}
where $\u n=(n_1,\dots,n_d)\in \Z^d_{\ge 0}$, $i\in I$, $\lambda_{\u n, i} \in \O$ and $\mu_{\u n, i, h} := \lambda_{\u n, i} \varpi^{-\sum_{j \le d; s > h} \lfloor \frac{n_j}{q^s}\rfloor} \to 0$ in $E$. Then for all but finitely many pairs $(\u n,i) \in \Z^d_{\ge 0} \times I$ we have $\mu_{\u n, i} \in \O$. We can therefore find an $h' \ge h$ such that $\mu_{\u n, i, h'}= \lambda_{\u n, i} \varpi^{-\sum_{j \le d; s > h'} \lfloor \frac{n_j}{q^s}\rfloor} \in \O$ for all pairs $(\u n,i)$. It follows that $f \in \LA_{h'}(\O_K^d,V^0) \subseteq \C\kan(\O_K^d,V^0)$, as required.
\end{proof}

\begin{lm}\label{lm:existence-universal-unitary}
Suppose that $\sigma$ is a locally analytic representation of $G$ of compact type such that there exist a BH-subspace $\sigma_0$ of $\sigma$ and elements $g_i \in G$ for $i \ge 1$ satisfying $\sigma = \sum_{i = 1}^\infty g_i \sigma_0$. Then $\sigma$ has a universal unitary completion. Explicitly, if $\Lambda_0$ is any unit ball in $\o{\sigma_0}$, then $\Lambda := \sum_{g \in G} g \Lambda_0$ is a minimal \(up to commensurability\) $G$-stable open lattice in $\sigma$.
\end{lm}

\begin{proof}
The proof only uses that $\sigma$ is an LF-space \cite[Def.\ 1.1.16]{locallyanalytic-memoir}. 
Since $\sigma$ is an LF-space we can write $\sigma$ as an increasing union of FH-subspaces $\sigma_i$ ($i \ge 1$) such that $\sigma \cong \ilim \o{\sigma_i}$ is a topological isomorphism. It is clear that $\Lambda$ as in the statement is a $G$-stable lattice in $\sigma$. We now verify that $\Lambda$ is open. Let $\sigma'_i := \sum_{j \le i} g_j \sigma_0$. Then by assumption, $\sigma$ is an increasing union of the FH-subspaces $\sigma'_i$ (using \cite[Prop.\ 1.1.5]{locallyanalytic-memoir}). By \cite[Prop.\ 1.1.10]{locallyanalytic-memoir} the sequences $\sigma_i$ and $\sigma'_i$ are mutually cofinal, and we obtain that $\sigma \cong \ilim \o{\sigma'_i}$ is a topological isomorphism. Therefore, to check that $\Lambda$ is open, it suffices to show that its preimage in $\o{\sigma'_i}$ is open for all $i$ (see \cite[\S 5.E]{nfa}). This is true since that preimage contains the lattice $\sum_{j \le i} g_j \Lambda_0$, which is easily checked to be open using the open mapping theorem.

Suppose now that $\Xi$ is any $G$-stable open lattice in $\sigma$. Then $\Xi \cap \o{\sigma_0}$ is open, hence after scaling $\Xi$ we may assume that $\Xi \cap \o{\sigma_0} \supseteq \Lambda_0$, so $\Xi$ contains $\sum_{g \in G} g \Lambda_0 = \Lambda$. This shows that $\Lambda$ is a minimal (up to commensurability) $G$-stable open lattice in $\sigma$, and hence that $\sigma$ has a universal unitary completion, namely $\sigma\sphat_{\Lambda}$.
\end{proof}

\begin{proof}[Proof of Proposition~\ref{prop:univ-unit-compl}]
Choose a cofinal sequence $M_1 \supseteq M_2 \supseteq \dots$ of analytic open subgroups of $M$ and for each $i$ let $\M_i$ be the affinoid rigid analytic $K$-group such that $M_i=\M_i(K)$.

\emph{Step 1:} We show that we can enlarge our given $\sigma_0$, if necessary, to ensure that there exists an analytic open subgroup $M_0'$ of $M$ such that $\sigma_0$ is $M_0'$-stable and the induced action of $M_0'$ on $\o{\sigma_0}$ is $\M_0'$-analytic (with obvious notation).

By assumption, $\sigma \cong \ilim_{i \ge 1} \o{\tau_i}$, where $\o{\tau_i}$ is a Banach space and the transition maps are injective (and compact). Let $\tau_i$ denote the image of $\o{\tau_i}$ in $\sigma$. By Propositions 3.2.15 and 1.1.10 in \cite{locallyanalytic-memoir} (using the completeness of $\sigma$) and an easy induction we may assume that each $\tau_i$ is $M_1$-stable. Thus, as $\sigma$ is a locally analytic representation of $M$, by \cite[Thm.~3.6.12]{locallyanalytic-memoir} we have topological isomorphisms
\begin{equation*}
 \sigma \cong \ilim_{j} \sigma_{\text{$\M_j$-an}}\cong \ilim_{j, i} (\o{\tau_i})_{\text{$\M_j$-an}} \cong \ilim_{i} (\o{\tau_i})_{\text{$\M_i$-an}}.
\end{equation*}
Let $\o{\sigma_i} := (\o{\tau_i})_{\text{$\M_i$-an}}$ and let $\sigma_i$ denote its image in $\sigma$, a BH-subspace.
By \cite[Prop.\ 1.1.10]{locallyanalytic-memoir} there exists an $i \ge 1$ such that $\sigma_0 \subseteq \sigma_i$.
Then we can enlarge $\sigma_0$ to $\sigma_i$ and take $M_0' := M_i$.

\emph{Step 2:} We show that $\sigma$ admits a universal unitary completion and identify a minimal (up to commensurability) $G$-stable open lattice in $\sigma$.

By \cite[Prop.\ 4.1.6]{jacquet1} we may choose an analytic open subgroup $H=\H(K)$ of $G$ that has a rigid analytic Iwahori decomposition with respect to $P$ and $\o P$ in the sense of \cite[Def.\ 4.1.3]{jacquet1}. This means in particular that there are closed rigid analytic subgroups $\M_0$, $\o\N_0$, $\N_0$ of $\H$ satisfying $\M_0(K) = H \cap M$, $\o\N_0(K) = H \cap \o N$, $\N_0(K) = H \cap N$, such that the multiplication map induces an isomorphism $\N_0 \times \M_0 \times \o\N_0\congto \H$ of rigid spaces. We let $M_0 := H \cap M$, $\o N_0 := H \cap \o N$, $N_0 := H \cap N$. Let $A_M\subseteq M$ be the split part of the connected centre of $M$ and choose $z \in A_M$ that strictly contracts each root subspace of the Lie algebra of $\o N$. By shrinking $H$, we may assume that $M_0 \subseteq M_0'$ (as analytic open subgroups) and that the conjugation map $\o N_0 \to \o N_0$, $\o n \mapsto z\o nz^{-1}$ is rigid analytic.

Note that $\o{\sigma_0}$ is an $\M_0$-analytic representation of $M_0$, as $M_0 \subseteq M_0'$. In particular, $M_0$ acts continuously on $\o{\sigma_0}$, so we can find an $M_0$-stable unit ball $\Lambda_0 \subseteq \o{\sigma_0}$ \cite[Lemma 6.5.4]{locallyanalytic-memoir}. Since $M/M_0$ is countable, Lemma~\ref{lm:existence-universal-unitary} implies that $\Lambda := \sum_{m \in M} m \Lambda_0$ is a minimal (up to commensurability) $M$-stable open lattice in $\sigma$, and that $\wh\sigma \cong \sigma\sphat_{\Lambda}$ exists.

Pick coset representatives $m_i$ ($i \ge 1$) for $M/M_0$. In particular, $\Lambda = \sum m_i \Lambda_0$. We define the increasing sequence of BH-subspaces $\sigma'_i := \sum_{j \le i} m_j \sigma_0$, so that $\sigma \cong \ilim \o{\sigma'_i}$ (see the proof of Lemma~\ref{lm:existence-universal-unitary}).

\emph{Step 3:} We set up some notation. Let $\pi := (\pInd \sigma)\kan$. Recall that $\C_c\kan(\o N,\sigma)$ is naturally a closed $\o P$-invariant subspace of $\pi$, with image consisting of all functions whose support is contained in $P \o N$ \cite[Lemma 2.3.6]{jacquet2}. We normalise the closed embedding $\C_c\kan(\o N,\sigma)\hookrightarrow \pi$ so that the inverse map (on its image) is $f\in \pi \longmapsto (\o n\in \o N\mapsto f(\o n^{-1}))$. The $\o P$-action on $\C_c\kan(\o N,\sigma)$ is then given explicitly by $(m \o n f)(\o n') = m f({\o n}^{-1} m^{-1} \o n' m)$. In particular, for any chart $X \congto \X(K)$ of $\o N$ and any BH-subspace $V$ of $\sigma$, $\C\rig(\X, \o V)$ embeds as a BH-subspace of $\C_c\kan(\o N,\sigma)$ and hence of $\pi$. For $m \in M$, $\o n \in \o N$ we see that $m \o n \C\rig(\X, \o V) = \C\rig(m \o n \X m^{-1}, \o{mV})$ and in the limit that $m \o n \C_c\kan(X, \o V) = \C_c\kan(m \o n X m^{-1}, \o{mV})$ (as subspaces of $\C_c\kan(\o N,\sigma)$). Similarly, if $V^0$ is a unit ball in $\o V$, then $m \o n \C\rig(\X, V^0) = \C\rig(m \o n \X m^{-1}, mV^0)$, where $mV^0$ denotes the image of $V^0$ inside $\o{mV}$ and $m \o n \C_c\kan(X, V^0) = \C_c\kan(m \o n X m^{-1}, mV^0)$.

\emph{Step 4:} We show that $\pi$ satisfies hypotheses \ref{item:1} and \ref{item:2} (for $G$ instead of $M$). It is obvious that $\pi$ has a central character. To verify hypothesis~\ref{item:2}, we take as BH-subspace $\C\rig(\o\N_0, \o{\sigma_0})$.
By compactness of $P\backslash G$ we can find a finite partition $P\backslash G = \bigsqcup_{j = 1}^{s} C_j$ into
compact open subsets and elements $g_j \in G$ such that $C_j \subseteq P\backslash P \o N g_j^{-1}$.
It is enough to show that
\begin{equation}\label{eq:30}
  \pi = \sum_{1 \le j \le s,i \ge 1,n \ge 0, \o\nu \in \o N/\o N_0} g_j m_i z^n \o\nu \C\rig(\o\N_0, \o{\sigma_0})
\end{equation}
with $\o N_0$, $\sigma_0$ and $(m_i)_{i\geq 1}$ as in Step 2.
Inside $\C_c\kan(\o N,\sigma)$ we have by the equalities in Step 3
$$\sum_{n\ge 0,\o\nu \in \o N/\o N_0} z^n \o\nu \C\rig(\o\N_0, \o{\sigma_0}) = \sum_{n,\o\nu} \C\rig(z^n \o\nu \o\N_0 z^{-n}, \o{\sigma_0})= \C_c\kan(\o N,\o{\sigma_0}),$$
as $\sigma$ has a central character and using Lemma~\ref{lm:cofinal}. From Corollary \ref{cor:additivity-analytic-lattices} and the bijectivity of~(\ref{eq:32}) we then deduce
\begin{equation*}
\sum_{i \ge 1} m_i \C_c\kan(\o N,\o{\sigma_0}) = \sum_{i\ge 1} \C_c\kan(\o N,\o{m_i \sigma_0}) = \bigcup_{\ell\ge 1} \C_c\kan(\o N,\o{\sigma_\ell'}) = \C_c\kan(\o N,\sigma)
\end{equation*}
with $\sigma_\ell'$ as in Step 2. (For the last equality note once again that the $\sigma_\ell'$ are cofinal among all BH-subspaces of $\sigma$.) By partitioning the support of functions in $\pi = (\pInd \sigma)\kan$ and translating, it follows that
$\pi = \sum_{j=1}^{s} g_j \C_c\kan(\o N,\sigma)$. We obtain~(\ref{eq:30}) by combining these equalities.

\emph{Step 5:} We now assume that the central character $\chi_\sigma$ of $\sigma$ is unitary. Define $\wt\Lambda := \sum_{g \in G} g \C\rig(\o\N_0,\Lambda_0)$ and $\wt\Lambda' := \pi \cap (\pInd \Lambda)\cts$, where $\Lambda=\sum_{m \in M} m \Lambda_0$ (see Step 2) and the intersection is taken inside $(\pInd \sigma)\cts$. We will show that $\wt\Lambda = \wt\Lambda'$, and that it is a minimal (up to commensurability) $G$-stable open lattice in $\pi$.

According to equation~(\ref{eq:30}) in Step 4 we may apply Lemma~\ref{lm:existence-universal-unitary} to the BH-subspace $\C\rig(\o\N_0,\o{\sigma_0})$ of $\pi$ and see that $\wt\Lambda$ is a minimal (up to commensurability) $G$-stable open lattice in $\pi$. Since $\C\rig(\o\N_0,\Lambda_0) \subseteq \C^0(\o N_0,\Lambda_0) \subseteq (\pInd \Lambda)\cts$ (noting that $M \Lambda_0 \subseteq \Lambda$), it follows that $\wt\Lambda \subseteq \wt\Lambda'$, as $\wt\Lambda'$ is $G$-stable by definition.

For the reverse inclusion, suppose that $f \in \wt\Lambda'$. By partitioning the support of $f$ and by translating (as in Step 4) we may assume, without loss of generality, that $\supp(f) \subseteq P \o N_0$, so $f \in \C\kan(\o N_0,\sigma) \cap \C^0(\o N_0,\Lambda)$. By Propositions 2.1.6 and 1.1.10 in \cite{locallyanalytic-memoir} and recalling that $\sigma \cong \ilim_i \o{\sigma_i'}$, there exists $i$ such that $f \in \C^0(\o N_0,\o{\sigma_i'})$. It follows that $\im(f) \subseteq \o{\sigma_i'}$ is compact. Recall that $\Lambda_k' := \sum_{j \le k} m_j \Lambda_0$ is an open lattice in $\o{\sigma_k'}$ (see the proof of Lemma~\ref{lm:existence-universal-unitary}), and note that $\Lambda = \bigcup_k \Lambda_k'$. Then for the value of $i$ above, $\Lambda \cap \o{\sigma_i'}$ is an increasing union of the open lattices $\Lambda_k' \cap \o{\sigma_i'}$ for $k \ge i$. By compactness, $\im(f) \subseteq \Lambda_k' \cap \o{\sigma_i'} \subseteq \Lambda_k'$ for some $k \ge i$. By increasing $i$, we conclude that $f \in \C\kan(\o N_0,\o{\sigma_i'}) \cap \C^0(\o N_0,\Lambda_i')$.

By Proposition~\ref{prop:amice}, we have $f \in \C\kan(\o N_0,\Lambda_i') \subseteq \C_c\kan(\o N,\Lambda_i')$. By Corollary~\ref{cor:additivity-analytic-lattices}, Lemma~\ref{lm:cofinal} and since $\sigma$ has a unitary central character we have (as in Step 4)
\begin{align*}
f \in \sum_{j \le i} \C_c\kan(\o N,m_j \Lambda_0) \ = \ \sum_{j \le i} m_j \C_c\kan(\o N,\Lambda_0) &= &\sum_{\o\nu \in \o N, j \ge 1, n \ge 0} m_j \C\rig(z^n \o\nu \o\N_0 z^{-n},\Lambda_0)\\
&= &\sum_{\o\nu \in \o N, j \ge 1, n \ge 0} m_j z^n \o\nu \C\rig(\o\N_0,\Lambda_0),
\end{align*}
and this is contained in $\wt\Lambda$, as required.

\emph{Step 6:} We show that $\wh \pi$ is canonically isomorphic to $(\pInd \wh\sigma)\cts$.

Consider the natural map $\theta : \pi = (\pInd \sigma)\kan \to (\pInd \sigma\sphat_{\Lambda})\cts$ (with $\Lambda$ as in Step 2 and Step 5). By choosing a locally analytic section of the map $G \to P \backslash G$, the map $\theta$ is identified with $\C\kan(P \backslash G,\sigma) \to \C^0(P \backslash G,\sigma\sphat_{\Lambda})$. Since the locally constant, $\sigma$-valued functions are dense in $\C^0(P \backslash G,\sigma\sphat_{\Lambda})$, we deduce that $\theta$ has dense image. Moreover, if $\wh\Lambda$ denotes the closure of the image of $\Lambda$ in $\sigma\sphat_{\Lambda}$ (whose preimage in $\sigma$ is $\Lambda$), then $\wt\Lambda'$ is the preimage under $\theta$ of the unit ball $(\pInd \wh\Lambda)\cts$. It follows that $\theta$ induces a topological isomorphism $\pi\sphat_{\wt\Lambda'} \congto (\pInd \sigma\sphat_{\Lambda})\cts$. By Step 2 and Step 5 we know that $\wh\sigma \cong \sigma\sphat_{\Lambda}$ and $\wh\pi \cong \pi\sphat_{\wt\Lambda'}$. This finishes the proof of Proposition \ref{prop:univ-unit-compl}.
\end{proof}

\section{A partial adjunction}\label{sec:partial-adjunction}

We prove an adjunction result (Proposition~\ref{prop:adjunct-max-subrep-ps-abstract}) which will be crucially used in the proof of the main result.\\

We first need to discuss some technical preliminaries. We let $G := \prod_{i=1}^d \GL_n(L_i)$ with $L_i/\qp$ finite for all $i$. We let $\o B$ (resp., $\o N$, resp., $B$, resp., $T$) denote the subgroup of upper-triangular (resp., upper-triangular unipotent, resp., lower-triangular, resp., diagonal) matrices in $G$. Let $\o\b, \b, \t$ denote the $\qp$-Lie algebras of $\o B$, $B$, $T$, respectively. Let $\o N_0$ denote any compact open subgroup of $\o N$ and let $T^+ := \{t \in T : t \o N_0 t^{-1} \subseteq \o N_0\}$. Then for any locally analytic representation $\pi$ of $\o B$ we have a Hecke action of $T^+$ on $\pi^{\o N_0}$ by $t\cdot v := \sum_{\o N_0/t \o N_0 t^{-1}} ntv$ for $t \in T^+$ and $v \in \pi^{\o N_0}$ (cf.\ \cite[\S3.4]{jacquet1}, but note that we do not normalise). Let $K := \prod_{i=1}^d \GL_n(\O_{L_i})$.

\begin{prop}[Emerton]\label{prop:breuil-633}
Suppose that $\Pi$ is an admissible continuous representation of $G$ on a Banach space and that $0 \to \Pi\lan[\qp] \to \pi \to \pi_1 \to 0$ is an exact sequence of admissible locally $\qp$-analytic representations of $G$, where $\Pi\lan[\qp]\subseteq \Pi$ is the subspace of locally $\qp$-analytic vectors for the action of $G$ \(an admissible locally $\qp$-analytic representation of $G$ by \cite[Thm.~7.1]{MR1990669}\). Suppose that $\chi : T \to E\s$ is a locally $\qp$-analytic character and $\eta : \t \to E$ its derivative. If $H^1(\o\b, \Pi\lan[\qp] \otimes (-\eta)) = 0$ then we have short exact sequences
\begin{gather}
 0 \to (\Pi\lan[\qp])^{\o N_0}[\t = \eta] \to \pi^{\o N_0}[\t = \eta] \to \pi_1^{\o N_0}[\t = \eta] \to 0,\label{eq:9} \\
 0 \to (\Pi\lan[\qp])^{\o N_0}[\t = \eta]_\chi \to \pi^{\o N_0}[\t = \eta]_\chi \to \pi_1^{\o N_0}[\t = \eta]_\chi \to 0,\label{eq:10}
\end{gather}
where $V[\t = \eta] \subseteq V$ \(resp.\ $V_\chi \subseteq V$\) denotes the largest subspace on which $\t$ acts via $\eta$ \(resp.\ the
generalised $\chi$-eigenspace under the action of $T^+$\), and all vector spaces in the last exact sequence are finite dimensional.
\end{prop}
\begin{proof}
This is the statement of \cite[Prop.\ 6.3.3]{breuil-ext1}, and we (finally!) give here full details of its proof. Let $T^0 := T \cap K$. The exactness of \eqref{eq:9} follows as in the proof of {\it loc.cit.} from $H^1(\o\b, \Pi\lan[\qp] \otimes (-\eta)) = 0$. To see the finite-dimensionality, we first show that
\begin{equation}
\pi^{\o N_0}[\t = \eta]_\chi = J_{\o B}(\pi)[\t = \eta]_\chi,\label{eq:22}
\end{equation}
where $J_{\o B}$ denotes Emerton's locally analytic Jacquet functor for the subgroup $\o B$ of $G$ (\cite{jacquet1}). Note that the $T^+$-action on the left extends uniquely to a $T$-action. We have
\begin{equation}\label{eq:21}
\pi^{\o N_0}[\t = \eta][(t-\chi(t))^s : t \in T^+] = J_{\o B}(\pi)[\t = \eta][(t-\chi(t))^s : t \in T^+]
\end{equation}
for all $s \ge 1$ by an analogue of \cite[Prop.\ 3.2.12]{jacquet1}, noting that the $T$-action on the left is locally $\qp$-analytic and $T$-finite (as $T^0$ acts semisimply and $T/T^0$ is finitely generated). Taking the union of~(\ref{eq:21}) over all $s \ge 1$, we obtain (\ref{eq:22}), since again $T^0$ acts semisimply and $T/T^0$ is finitely generated.

Next, \cite[Prop.\ 4.2.33]{jacquet1} shows that the strong dual $J_{\o B}(\pi)[\t = \eta]'_b$ is the space of global sections of a coherent sheaf $\F$ on $\wh T$ (the rigid analytic variety of locally $\qp$-analytic characters of $T$) with discrete support. We write $\wh T = \bigcup_{n\ge 1}\wh T_n$ as increasing union of an admissible cover of affinoid subdomains. Then $\F(\wh T_n)\cong \oplus_{x \in \wh T_n} \F(\wh T_n)_{\m_x}$, where the direct sum is over all closed points $x$ of $\wh T_n$ and $\m_x$ is the corresponding maximal ideal. It is a finite decomposition, as $\supp \F$ is discrete; moreover, $\m_x^{s(x)}$ annihilates $\F(\wh T_n)_{\m_x}$ for some $s(x) \ge 1$ and $\F(\wh T_n)$ is finite-dimensional. Also, $\F(\wh T_n)_{\m_x} \cong \F_{x}$ by \cite[Prop.\ 9.4.2/6]{MR746961}. By passing to the limit over $n$ we get a topological isomorphism $J_{\o B}(\pi)[\t = \eta]'_b \cong \prod_{x \in \wh T} \F_{x}$ and dually
\begin{equation*}
 J_{\o B}(\pi)[\t = \eta] \cong \bigoplus_{x \in \wh T} (\F_{x})'_b.
\end{equation*}
For any $t \in T$ the function $t - \chi(t)$ on $\wh T$ is contained in the maximal ideal $\m_\chi$ of $\chi \in \wh T$, hence by above $(t-\chi(t))^{s(\chi)}$ annihilates $\F_\chi$ for all $t \in T$. Therefore, $J_{\o B}(\pi)[\t = \eta]_\chi \cong (\F_{\chi})'_b$ is finite-dimensional.

To see the exactness of \eqref{eq:10}, we first note that by the above we can identify it with the sequence $0 \to J_{\o B}(\Pi\lan[\qp])[\t = \eta]_\chi \to J_{\o B}(\pi)[\t = \eta]_\chi \to J_{\o B}(\pi_1)[\t = \eta]_\chi \to 0$, and thus it is independent of the choice of $\o N_0$.

Next, let $H := \prod_{i=1}^d \{ g \in \GL_n(\O_{L_i}) : g \equiv 1 \pmod {p^r} \}$ for some large $r$ and $\o N_0 := H \cap \o N$ (these are uniform pro-$p$ groups). Then $T^0$ normalises $\o N_0$, so the Hecke action of $T^0$ on $(\Pi\lan[\qp])^{\o N_0}$ coincides with its natural action. Using \cite[Prop.\ 3.3.2]{jacquet1} we choose $z_1,\dots,z_s \in T^+$ such that $z_i \o N_0 z_i^{-1} \subseteq \o N_0^p$ for all $i$ and such that $T^0$, $z_1,\dots,z_s$ generate $T$ as group. By increasing $r$ we may assume that moreover $z_i^{-1} H z_i \subseteq K$ for all $i$. Then by \cite[Lemma~5.3]{BHS2} we know that $\Pi\lan[\qp]$ is an increasing union of $H$-stable BH-subspaces $\Pi^{(h)}$ for $h\ge 1$ and, moreover, for each $h \ge 1$ the elements $z_i$ preserve the BH-subspace $\Pi_h := \Pi^{(h)}\cap (\Pi\lan[\qp])^{\o N_0}[T^0 = \chi]$ and induce a compact operator on the Banach space $\o{\Pi_h}$ (a closed subspace of $\o{\Pi^{(h)} \cap (\Pi\lan[\qp])^{\o N_0}[\t = \eta]}$), see \S\ref{sec:notation} for the notation.

We observe that if $Y \subseteq T^+$ is any submonoid that contains $T^0$ and generates $T$ as group, then we can also compute $J_{\o B}(\pi)$ using $Y$ instead of $T^+$ by \cite[Lemma 3.2.19]{jacquet1}. In fact, all results of \cite[\S3.2]{jacquet1} apply and we deduce in particular that $\pi^{\o N_0}[\t = \eta]_{Y=\chi} = J_{\o B}(\pi)[\t =\eta]_{Y=\chi}$ by the argument at the beginning of our proof, where the subscript $Y=\chi$ denotes the generalised eigenspace for $\chi|_Y$. As $Y$ generates $T$ as group we also have $J_{\o B}(\pi)[\t =\eta]_{Y=\chi} = J_{\o B}(\pi)[\t = \eta]_\chi$. Thus the space $\pi^{\o N_0}[\t = \eta]_{Y=\chi} = \pi^{\o N_0}[T^0 = \chi]_{Y=\chi}$ is finite-dimensional and independent of the choice of $Y$.

Let $Y$ denote the submonoid of $T^+$ generated by $T^0$, $z_1,\dots,z_s$. From \eqref{eq:9} we deduce an exact sequence $0 \to (\Pi\lan[\qp])^{\o N_0}[T^0 = \chi] \to \pi^{\o N_0} [T^0 = \chi] \to \pi_1^{\o N_0}[T^0 = \chi] \to 0$ and it remains to show that the last map is surjective on generalised $\chi|_Y$-eigenspaces. Let $M_1 := \pi_1^{\o N_0}[T^0 =\chi]_{Y=\chi}=\pi_1^{\o N_0}[\t = \eta]_{\chi}$ and let $M$ be the preimage of $M_1$ in $\pi^{\o N_0}[T^0 = \chi]$. Let $V$ be any subspace of $\pi^{\o N_0}[T^0 = \chi]$ mapping isomorphically onto the finite-dimensional space $M_1$. Now choose $h \ge 1$ such that $z_i V \subseteq V + \Pi_h$ for all $i$ (recall $(\Pi\lan[\qp])^{\o N_0}[T^0 = \chi]$ is the increasing union of the $\Pi_h$), and let $V' := \Pi_h\oplus V$ (a BH-subspace of $\pi^{\o N_0} [T^0 = \chi]$). We have a $Y$-equivariant exact sequence of Banach spaces $0 \to
\o{\Pi_h} \to \o{V'} \to M_1 \to 0$. As the $z_i$ are compact on $\o{\Pi_h}$ and $\dim_E M_1 <\infty$, we deduce that the $z_i$ are compact on $\o{V'}$. It follows that the generalised eigenspace $V'_{Y = \chi}$ is finite-dimensional and surjects onto $M_1$. This finishes the proof of \eqref{eq:10}.
\end{proof}

\begin{ex}\label{ex:cohomology}
Suppose that $\Pi|_K \cong \C^0(K,E)^{\oplus r}$ for some $r$, or just that $\Pi\lan[\qp]|_K \cong \C\qpan(K,E)^{\oplus r}$, then $H^1(\o\b, \Pi\lan[\qp] \otimes (-\eta)) = 0$ for all linear maps $\eta :\t \to E$ by the proof of \cite[Prop.\ 6.3.3]{breuil-ext1}.
\end{ex}

We will need the following corollary for future reference.

\begin{coroll}[Emerton]\label{cor:breuil-633}
Suppose that $\pi$ is an admissible locally $\qp$-analytic representations of $G$ and that $\chi : T \to E\s$ is a locally $\qp$-analytic character. Then $\Hom_{T^+}(\chi,\pi^{\o N_0}) = \Hom_T(\chi,J_{\o B}(\pi))$ is finite-dimensional.
\end{coroll}
\begin{proof}
The two spaces are equal by \cite[Prop.\ 3.2.12]{jacquet1}. Moreover, the image of any $T^+$-equivariant linear map $\chi \to \pi^{\o N_0}$ has to be contained in $\pi^{\o N_0}[\t = \eta]_\chi$, which is finite-dimensional by Proposition~\ref{prop:breuil-633}.
\end{proof}

We recall that for $M \in \O_{\alg}^{\b}$ and $\pi^\infty$ an admissible smooth representation of $T$, Orlik--Strauch
have defined an admissible locally $\qp$-analytic representation $\F_B^G(M,\pi^\infty)$
(see the paragraph before Lemma~\ref{lm:OS-product-groups}). Let $X(\u T)$ denote the $\qp$-algebraic characters $T \to E\s$. The derivative
of any $\mu \in X(\u T)$ gives a $\qp$-linear map $\t \to E$ and hence a Verma module $M(\mu) \in \O_{\alg}^{\b}$ whose simple cosocle 
we denote by $L(\mu)$ (cf.\ \cite[\S2]{socle1}).

\begin{prop}\label{prop:adjunct-max-subrep-ps-abstract}
Suppose that $\Pi$ is an admissible continuous representation of $G$ on a Banach space together with an action of a commutative algebra $\T$ by continuous $G$-linear endomorphisms. Let $\m \lhd \T$ be a maximal ideal with residue field $E$ and $\sigma$ an admissible locally $\qp$-analytic representation of $G$ of finite length such that:
\begin{enumerate}
\item we have $H^1(\o\b, \Pi\lan[\qp] \otimes (-\eta)) = 0$ for all linear maps $\eta : \t \to E$;
\item we are given a \($G$-equivariant continuous\) injection $j : \soc_G \sigma \into \Pi\lan[\qp][\m]$;
\item any irreducible constituent of $\sigma/\soc_G \sigma$ is of the form $\F_B^G(L(\mu),\chi)$ for some $\mu \in X(\u T)$ and some smooth character $\chi : T\to E\s$ such that moreover no irreducible subquotient of $\F_B^G(M(\mu),\chi)$ injects into $\Pi\lan[\qp][\m]$.
\end{enumerate}
Then $j$ extends uniquely to a \($G$-equivariant continuous\) map $\sigma \to \Pi\lan[\qp][\m]$ that is moreover injective.
\end{prop}
\begin{proof}
The uniqueness of the extension is clear by the last assumption. To show existence we follow the strategy of \cite[\S6.4]{breuil-ext1}. For existence, we may assume by induction that we have already extended $j$ to an injection $j_0 : \sigma' \into \Pi\lan[\qp][\m]$, where $\soc_G \sigma \subseteq \sigma' \subseteq \sigma$ and $C := \sigma/\sigma'$ is irreducible. We define the amalgam $S := \Pi\lan[\qp] \oplus_{j_0,\sigma'} \sigma$ and let $\T$ act on $S$ by declaring that it acts through $\T/\m \cong E$ on $\sigma$. Then $S[\m] = \Pi\lan[\qp][\m] \oplus_{j_0,\sigma'} \sigma$, and so we have an exact sequence $0 \to \Pi\lan[\qp][\m] \to S[\m] \to C \to 0$. By assumption, we can write $C \cong \F_B^G(L(\mu),\chi)$, and hence by \cite[Thm.\ 4.3]{socle2},
\begin{equation}\label{eq:6}
\Hom_{T^+}(t^{-\mu} \chi,C^{\o N_0}) \ne 0.
\end{equation}
(We note that \cite{socle2} works with the group of $L$-points of a split reductive group over $L$. However, the proofs work unchanged for our group $G$. See also \cite[Rk.\ 5.1.2]{BHS3}.) Let $\eta : \t \to E$ denote the derivative of $t^{-\mu} \chi$, i.e.\ $\eta = -\mu$. By Proposition~\ref{prop:breuil-633} we have an exact sequence of finite-dimensional vector spaces
\begin{equation}\label{eq:8}
0 \to (\Pi\lan[\qp])^{\o N_0}[\t = \eta]_{t^{-\mu}\chi} \to S^{\o N_0}[\t = \eta]_{t^{-\mu}\chi} \to C^{\o N_0}[\t = \eta]_{t^{-\mu}\chi} \to 0.
\end{equation}
If we had $((\Pi\lan[\qp])^{\o N_0}[\t = \eta]_{t^{-\mu}\chi})_\m \ne 0$ then $\Hom_{T^+}(t^{-\mu} \chi,(\Pi\lan[\qp][\m])^{\o N_0}) \ne 0$, so by \cite[Thm.\ 4.3]{socle2} we would get a non-zero $G$-linear map $\F_B^G(M(\mu)\dual,\chi)\to \Pi\lan[\qp][\m]$ (see {\it loc.cit.} for the notation) and hence some irreducible constituent of $\F_B^G(M(\mu)\dual,\chi)$, or equivalently of $\F_B^G(M(\mu),\chi)$, would inject into $\Pi\lan[\qp][\m]$, contradicting our last assumption. Therefore the second and third terms in the sequence~\eqref{eq:8} become isomorphic after $\m$-localisation, and thus also on taking $\m$-torsion. In particular, we have an isomorphism
\begin{equation}\label{eq:7}
\Hom_{T^+}(t^{-\mu} \chi,S[\m]^{\o N_0}) \congto \Hom_{T^+}(t^{-\mu} \chi,C^{\o N_0}),
\end{equation}
and they are non-zero by \eqref{eq:6}. (In fact they are one-dimensional, but we will not need that.)

As in \cite[\S6.4]{breuil-ext1} a non-zero element of $\Hom_{T^+}(t^{-\mu} \chi,S[\m]^{\o N_0})$ corresponds to a non-zero $(\g,\o B)$-linear map
\begin{equation}\label{map}
U(\g) \otimes_{U(\o\b)} (-\mu) \otimes_E \C_c^\infty(\o N,\chi) \to S[\m],
\end{equation}
where $\C_c^\infty(\o N,\chi)$ is the $\o B$-representation defined in {\it loc.cit.} twisted by the character $\chi$. The map (\ref{map}) factors through a map
\begin{equation*}
\psi_1 : C\lp(\o N) = \o L(-\mu) \otimes_E \C_c^\infty(\o N,\chi) \to S[\m]
\end{equation*}
by our last assumption, where $\o L(-\mu)$ denotes the simple module of highest weight $-\mu$ in $\O_{\alg}^{\o\b}$. (For $C\lp(\o N)\subseteq C$ we use the notation of \cite[\S2.7]{jacquet2} and for the identification with $\o L(-\mu) \otimes_E \C_c^\infty(\o N,\chi)$, see the proof of \cite[Prop.\ 4.2]{socle2}.) By the isomorphism~\eqref{eq:7} we may assume that the composite of $\psi_1$ with the natural map $S[\m] \to C$ is the inclusion $C\lp(\o N) \to C$.

Let us write $\sigma = (\sigma' \text{ --- } C)$ and let $\sigma' \text{ --- }C\lp(\o N)$ be the pullback of $\sigma$ along the inclusion $C\lp(\o N) \hookrightarrow C$. Let $\psi_2$ be the composite $\sigma' \text{ --- }C\lp(\o N) \into \sigma \into S[\m]$, so $\psi_2|_{\sigma'} = j_0$ and let $s$ denote the projection $(\sigma' \text{ --- }C\lp(\o N)) \onto C\lp(\o N)$. Then by construction, $\psi := \psi_2-\psi_1\circ s : (\sigma' \text{ --- }C\lp(\o N)) \to \Pi\lan[\qp][\m]$ and it restricts to $j_0$ on $\sigma'$. By \cite[Thm.\ 7.1.1]{breuil-ext1}, the map $\psi$ extends uniquely to a $G$-equivariant map $\sigma \to \Pi\lan[\qp][\m]$. This latter map is injective, since it is non-zero on $\soc_G \sigma$ (as $j_0$ is).
\end{proof}

\begin{rk}\label{rk:hypothesis-partial-adj}
At least in our global application, the condition that no irreducible subquotient of $\F_B^G(M(\mu),\chi)$ injects into $\Pi\lan[\qp][\m]$ in Proposition \ref{prop:adjunct-max-subrep-ps-abstract}(iii) will be equivalent to demanding that $\F_B^G(L(\mu),\chi)$ does not inject into $\Pi\lan[\qp][\m]$ by Conjecture~\ref{conj:loc-an-socle} (which holds in many cases by \cite[Thm.~1.3]{BHS3}), see the proof of Theorem~\ref{thm:adjunct-max-subrep-ps-coho}.
\end{rk}

\begin{rk}
We note that Proposition \ref{prop:adjunct-max-subrep-ps-abstract} generalises \cite[Thm.~B]{BC-adjunction} when the parabolic subgroup $P$ of {\it loc.cit.} is the Borel subgroup and the character is locally algebraic.
\end{rk}

\section{The finite slope space in the generic crystabelline case}\label{sec:fin-slope-gen-crystab}

We prove our main result (Corollary \ref{cor:adjunct-max-subrep-ps-coho}).

\subsection{Local setup and results}\label{sec:local-setup-crystab}

We define and study the ``finite slope'' representation $\Pi(\rho)\fs$.\\

Let $L/\qp$ be a finite extension and set $\cS_L := \Hom_{\qp}(L,E)$. We fix a crystabelline representation $\rho : \Gal(\lb/L) \to \GL_n(E)$ satisfying the following genericity hypothesis.

\begin{hyp}\label{hyp:generic-crystabelline}
We assume that $\rho$ is potentially crystalline with $\WD(\rho) = \oplus_{i=1}^n \chi_i$ for some smooth characters $\chi_i : W_L \to E\s$ (and $N =0$), that $\chi_i \chi_j^{-1} \not\in \{ 1, |\cdot|_L^{\pm 1}\}$ for all $i \ne j$, and that for each $\sigma \in \cS_L$ the $\sigma$-Hodge--Tate weights of $\rho$ are distinct.
\end{hyp}

Let $L'/L$ be finite abelian such that $\rho|_{\Gal(\lb/L')}$ is crystalline. Then, in particular, the Deligne--Fontaine module $D := D\cris(\rho|_{\Gal(\lb/L')})$ covariantly associated to $\rho$ satisfies Hypotheses 5.1 and 5.2 in \cite{socle1}. For each $\sigma \in \cS_L$, let $h_{1,\sigma} < \cdots < h_{n,\sigma}$ denote the jumps in the Hodge filtration on $D_{L',\sigma} := D_{L'} \otimes_{L \otimes_{\qp} E,\sigma} E$ (or equivalently, in the $n$-dimensional $E$-vector space $D_{L',\sigma}^{\Gal(L'/L)}$), and let $\lambda_{i,\sigma} := -h_{i,\sigma}-(n-i)$ for all $1 \le i \le n$, so that $\lambda_{1,\sigma}\geq \cdots \geq \lambda_{n,\sigma}$.

Let $G := \GL_{n/L}$ and $\u G := \Res_{L/\qp} G \times_{\qp} E \cong \prod_{\cS_L} \GL_{n/E}$. Let $B$ (resp.\ $\o B$) denote the {\it lower}-triangular (resp.\ {\it upper}-triangular) Borel subgroup of $G$, let $T$ denote the diagonal maximal torus and $W \cong S_n$ the Weyl group of $(G,T)$ with Bruhat order $\le$, and let $w_0$ denote the longest element of $W$. We let $\Delta \subseteq X(T)$ denote the simple roots of $G$ with respect to $\o B$. Similarly we define $\u B, \u{\o B}, \u T, \u W, \u
w_0, \u \Delta$ for $\u G$. Note that we can naturally identify $W$ as a subgroup of $\u W$. Also, we can and will think of $\lambda=(\lambda_{i,\sigma})_{i,\sigma}$ as an element of $X(\u T)$, which is dominant with respect to ${\u {\o B}}$.

We say that a \emph{refinement} of $D$ is a complete flag of Deligne--Fontaine submodules of $D$. Equivalently, it is a complete flag of Weil--Deligne subrepresentations of $\WD(D) = \WD(\rho)\cong \oplus_{i=1}^n \chi_i$. Since the characters $\chi_i$ are distinct by assumption, every refinement of $D$ is of the form
\begin{equation}
\F : 0 \subseteq D_{\F,1} \subseteq D_{\F,1} \oplus D_{\F,2} \subseteq \cdots \subseteq D,\label{eq:24}
\end{equation}
where the $D_{\F,i}$ are pairwise distinct Deligne--Fontaine submodules of $D$ of rank 1. Let $\chi_{\F,i} := \WD(D_{\F,i})$, which is a smooth character of $W_L\ab$ or equivalently of $L\s$. Thus we see that a refinement of $D$ can also be thought of as an ordering $\chi_{\F,1}$, \dots, $\chi_{\F,n}$ of the characters $\chi_i$. Moreover,
we have a simply transitive action of $W$ on the set $\Ref(D)$ of refinements of $D$, given by
\[w\F : 0 \subseteq D_{\F,w^{-1}(1)} \subseteq D_{\F,w^{-1}(1)} \oplus D_{\F,w^{-1}(2)} \subseteq \cdots \subseteq D,\] 
where $w \in W$ and $\F \in \Ref(D)$ is as in~(\ref{eq:24}), i.e.\ $\chi_{w\F,i} = \chi_{\F,w^{-1}(i)}$.

Let $\W := {\u W} \times \Ref(D)$. For $(w\alg,\F) \in \W$ we define a locally $\qp$-algebraic character $\eta(w\alg,\F): T(L)\rightarrow E^\times$ as follows:
\begin{equation}\label{eq:23}
{\left( \begin{smallmatrix} x_1& & \\ & \ddots &\\ && x_n \end{smallmatrix} \right)}\in T(L)\longmapsto 
\prod_{i=1}^n \Big(\prod_{\sigma \in \cS_L} \sigma(x_i)^{-h_{(w\alg)^{-1}(i),\sigma}}\Big) \chi_{\F,i}(x_i) \ve^{-(n-i)}(x_i)
\end{equation}
and we let $\pi_{B,\F} := \chi_{\F,1}|\cdot|_L^{-(n-1)} \otimes \chi_{\F,2}|\cdot|_L^{-(n-2)}\otimes \cdots \otimes \chi_{\F,n}$ denote the smooth part of $\eta(w\alg,\F)$. We define the locally $\qp$-analytic principal series representation of $G(L)$
\begin{equation}\label{PS}
\PS(w\alg,\F):=\Big(\Ind_{B(L)}^{G(L)} \eta(w\alg,\F)\Big)\qpan
\end{equation}
(which is of compact type and admissible). For $w\alg \in \u W$ let $W(w\alg)$ denote the subgroup of $W$ generated by all reflections $s_\alpha$ for $\alpha \in \Delta$ with $s_\alpha w\alg > w\alg$ (where we see $W$ as subgroup of $\u W$) and $P(w\alg)$ the parabolic subgroup of $G$ containing $B$ with Weyl group $W(w\alg)$.

\begin{lm}\label{lm:princ-series}\
\begin{enumerate}
\item The principal series $\PS(w\alg,\F)$ has finite length and irreducible socle, which we denote by $C(w\alg,\F)$.
\item We have $C(w\alg_1,\F_1) \cong C(w\alg_2,\F_2)$ if and only if $w\alg_1 = w\alg_2$ and $W(w\alg_1) \F_1 = W(w\alg_2) \F_2$.
\item The irreducible constituents of $\PS(w\alg,\F)$ are the representations $C(\tau,\F)$ for $\tau \in \u W$, $\tau \ge w\alg$. They occur only once in case $\tau \in \{ w\alg, \u w_0 \}$.
\end{enumerate}
\end{lm}
\begin{proof}
Part (i) is a special case of the results in \cite[\S6]{socle1}: using the theory of Orlik--Strauch \cite{OS}, as extended in \cite{socle1}, we can write $\PS(w\alg,\F) = \cF^{G}_{B}(M(w\alg \cdot (-\lambda)),\pi_{B,\F})$, so by \cite[Cor.\ 2.5]{socle1} we have that $C(w\alg,w) = \cF^{G}_{B}(L(w\alg \cdot (-\lambda)),\pi_{B,\F})$ is the (irreducible) socle of $\PS(w\alg,\F)$. (We remark that for us, the construction of Orlik--Strauch, as well as the dot action of $\u W$ on $X(\u T)$, are defined relative to our choice of lower-triangular Borel subgroup $\u B$, just as in \cite{socle2}. In \cite{socle1} the dot action was defined relative to $\u{\o B}$.) Part (ii) now follows from \cite[Lem.~6.2]{socle1}, the sentence before \cite[Lem.~6.3]{socle1}, and \cite[Lem.~4.2]{socle1}. Since $-\lambda$ is dominant with respect to  $\u B$, the Verma module $M(w\alg \cdot (-\lambda))$ has finite length and constituents $L(\tau \cdot (-\lambda))$ for $\tau \ge w\alg$, which occur only once in case $\tau \in \{ w\alg, \u w_0 \}$ (see for example the bottom of p.~155 in \cite{humphreys-bgg} for the latter case). As the smooth induction $(\Ind_B^G \pi_{B,\F})^\infty$ is irreducible by genericity, parts (i) and (iii) follow from the main results of \cite{OS}.
\end{proof}

By base change a refinement $\F \in \Ref(D)$ gives rise to a $\Gal(L'/L)$-stable flag $\F_{L'}$ of $D_{L'} := L' \otimes_{L_0'} D$. On the other hand, forgetting the indexation, the Hodge filtration $\Fil^*D_{L'}$ gives rise to another $\Gal(L'/L)$-stable flag of $D_{L'}$. By Galois descent we obtain two complete flags $(\F_{L'})^{\Gal(L'/L)}$, $(\Fil^*D_{L'})^{\Gal(L'/L)}$ on the free rank-$n$ $L \otimes_{\qp} E$-module $(D_{L'})^{\Gal(L'/L)}$. Their relative position is given by an element of $\u W$, which we denote by $w\alg(\F)\u{w_0}$ (thus defining $w\alg(\F)$). Explicitly, if $\alpha : (L \otimes_{\qp} E)^n \congto (D_{L'})^{\Gal(L'/L)}$ is any isomorphism of $L \otimes_{\qp} E$-modules, the flags $\alpha^{-1}((\F_{L'})^{\Gal(L'/L)})$, $\alpha^{-1}((\Fil^*D_{L'})^{\Gal(L'/L)})$ are described by an element of
$$\u G(E)\backslash \big((\u G(E)/\u{\o B}(E)) \times (\u G(E)/\u{\o B}(E))\big)$$
which is independent of the choice of $\alpha$. We mean that this element is in the same coset as $(1,w\alg(\F)\u{w_0})$.

\begin{rk}\label{rk:comparing-C-walg-F-with-socle-papers}
We relate $C(w\alg,\F)$ and $w\alg(\F)$ to the notions introduced in \cite[\S6]{socle1}, \cite[\S6]{socle2}. In those references, a refinement $\F$ is fixed at the outset to define $C(w\alg,w)$ and $w\alg(w)$. To indicate the dependence on $\F$, we write $C_\F(w\alg,w)$ and $w\alg_\F(w)$ in this remark. With this convention, we have $C_\F(w\alg,w) = C(w\alg,w\F)$ and $w\alg_\F(w) = w\alg(w\F)$. (The latter equality follows from $w\alg_\F(w) = w\alg_{w\F}(1)$, which holds by an elementary argument just as in the proof of \cite[Prop.\ 6.4(i)]{socle1}.)
\end{rk}

Now, we let $\Ws(\rho) := \{ (w\alg, \F) \in \W : w\alg \le w\alg(\F) \}$ and
\begin{equation*}
\Cs(\rho) := \{ \text{$C(w\alg, \F)$ (up to isomorphism) such that\ } (w\alg, \F) \in \Ws(\rho) \}.
\end{equation*}
We recall that if $C(w_1\alg, \F_1) \cong C(w_2\alg, \F_2)$, then $w_1\alg = w_2\alg=:w\alg$ and we have $w\alg \le w\alg(\F_1)$ if and only if $w\alg \le w\alg(\F_2)$ (as follows from \cite[Lem.\ 6.3]{socle2}). The following construction takes place in the abelian category of admissible locally $\qp$-analytic representations (in fact, all representations that are involved are of finite length).

\begin{df}\label{our-amalgam}\
\begin{enumerate}
\item For any $(w\alg, \F) \in \Ws(\rho)$, let $M(w\alg,\F)$ be the largest (non-zero) subrepresentation of the principal series $\PS(w\alg,\F)$ such that none of the irreducible constituents of $M(w\alg,\F)/C(w\alg,\F)$ is contained in $\Cs(\rho)$.
\item For any $C \in \Cs(\rho)$ let $M(\rho)_C$ denote the amalgam
\begin{equation*}
 \bigoplus_{C} \{ M(w\alg,\F) : (w\alg,\F) \in \Ws(\rho),\ C(w\alg,\F) \cong C\}
\end{equation*}
over the common socle $C$.
\item Let $\Pi(\rho)\fs_C$ denote the unique quotient of $M(\rho)_C$ whose socle is isomorphic to $C$.
\item Let $\Pi(\rho)\fs := \bigoplus_{C \in \Cs(\rho)} \Pi(\rho)\fs_C$.
\end{enumerate}
\end{df}

\begin{rk}\label{rk:our-amalgam}\
 \begin{enumerate}
\item By construction, $C$ is contained in the socle of $M(\rho)_C$, but equality does not hold in general (see the examples in \S\ref{sec:examples}). Moreover, the quotient $\Pi(\rho)\fs_C$ is well-defined, since $C$ occurs just once in $M(\rho)_C$.
\item Note that $M(w\alg,\F) = \PS(w\alg,\F)$ if and only if $w\alg = w\alg(\F)$ (use Lemma \ref{lm:princ-series}(iii)).
\item Note that $M(w\alg,\F)$ injects into $\Pi(\rho)\fs_{C(w\alg,\F)}$ for each $(w\alg,\F) \in \Ws(\rho)$.
\end{enumerate}
\end{rk}

\begin{prop}\label{prop:constituents-of-Qrho}
Suppose that $\F \in \Ref(D)$.
\begin{enumerate}
\item If $\alpha \in \u \Delta$, then $C(s_\alpha,\F)$ occurs precisely once in $\Pi(\rho)\fs$.
\item The representation $C(\u w_0,\F)$ occurs precisely once in $\Pi(\rho)\fs$.
\end{enumerate}
\end{prop}
\begin{proof}
(i) If $(s_\alpha,\F) \in \Ws(\rho)$, then $C := C(s_\alpha,\F)$ is a subrepresentation of $\Pi(\rho)\fs$ and occurs only once. Otherwise, if $C$ occurs in $M(w\alg,\F')$ with $(w\alg,\F') \in \Ws(\rho)$, then $w\alg = 1$ and $W(s_\alpha) \F = W(s_\alpha) \F'$ by Lemma~\ref{lm:princ-series}. By the Kazhdan--Lusztig conjectures (or Jantzen's multiplicity 1 criterion), the Verma module $M(-\lambda)$ contains the constituent $L(s_\alpha\cdot(-\lambda))$ with multiplicity one, and in the second radical layer. Therefore $\PS(1,\F')$, having an irreducible socle, contains a non-split extension $\E$ of the form $C(1,1) \text{---} C$ as subrepresentation, and hence so does $M(1,\F')$. More precisely, we have $\E \cong \F_B^G(M, \pi_{B,\F'})$, where $M$ is the length two quotient of $M(-\lambda)$ with socle $L(s_\alpha\cdot(-\lambda))$. As $s_\alpha\cdot(-\lambda)$ and $-\lambda$ are dominant with respect to the Borel $B \cap L_{P(s_\alpha)}$ of $L_{P(s_\alpha)}$, it follows that $\E \cong \F_{P(s_\alpha)}^G(M, (\Ind_{B(L)}^{P(s_\alpha)(L)} \pi_{B,\F'})^\infty)$, and hence is independent of the choice of $\F'$, as $(\Ind_{B(L)}^{P(s_\alpha)(L)} \pi_{B,\F'})^\infty$ is by the genericity conditions. Therefore, all occurrences of $C$ inside $M(\rho)_{C(1,1)}$ are contained in an amalgam of $r := \# W(s_\alpha)$ copies of the extension $\E$ over the common socle $C(1,1)$. But this amalgam is easily seen to be isomorphic to $\E \oplus C^{r-1}$, hence there is only one copy of $C$ in $\Pi(\rho)\fs_{C(1,1)}$.

(ii) On the one hand, $M(w\alg(\F),\F) = \PS(w\alg(\F),\F)$ contains $C := C(\u w_0,\F)$ precisely once as constituent by Lemma~\ref{lm:princ-series}. On the other hand, if $C$ occurs in $M(w\alg,\F')$ with $(w\alg,\F') \in\Ws(\rho)$, then $\F = \F'$ as $W(\u w_0) = 1$, so $w\alg \le w\alg(\F)$. If $w\alg < w\alg(\F)$, then we have a surjection $\PS(w\alg,\F) \onto \PS(w\alg(\F),\F)$ (because of a corresponding injection of Verma modules). It has to send $M(w\alg,\F)$ to zero, as otherwise $M(w\alg,\F)$ would contain $C(w\alg(\F),\F) \in \Cs(\rho)$ as constituent (not in its socle), so $C$ does not occur in $M(w\alg,\F)$. Hence $w\alg = w\alg(\F)$ and we are done, since $M(w\alg(\F),\F)$ injects into $\Pi(\rho)\fs_{C(w\alg(\F),\F)}$.
\end{proof}

\subsection{Global applications}\label{sec:global-app-crystab}

We prove our main global results (Theorem \ref{thm:adjunct-max-subrep-ps-coho} and Corollary \ref{cor:adjunct-max-subrep-ps-coho}).\\

We first explain our global setup, which is essentially the same as that of \cite[\S5--\S6]{socle2} (except we do not assume that $p$ splits completely in our totally real field), to which we refer the reader for further details and references. We fix a totally real number field $F^+ \ne \Q$ and a totally imaginary quadratic extension $F/F^+$. We let $c$ denote the unique complex conjugation of $F$ and suppose that every place $v|p$ of $F^+$ splits in $F$. We let $\u G$ be a unitary group over $F^+$ defined by a hermitian form of rank $n$ over $F$, so we have an isomorphism $\iota_G : \u G \times_{F^+} F \congto \GL_n$. We assume moreover that the hermitian form is totally definite, i.e.\ $\u G(F^+ \otimes_\Q \R)$ is compact. We fix a compact open subgroup $U^p$ of $\u G(\A_{F^+}^{\infty,p})$ of the form $U^p = \prod_{v \nmid p} U_v$ for compact open subgroups $U_v \subseteq \u G(F^+_v)$ for $v \nmid p$. We define
\begin{equation*}
 \wh S(U^p, E) := \C^0(\u G(F^+) \backslash \u G(\A_{F^+}^\infty) / U^p, E),
\end{equation*}
which is a Banach space for the supremum norm on the (profinite) compact topological space $\u G(F^+) \backslash \u G(\A_{F^+}^\infty) / U^p$. We will sometimes just write $\wh S$ for $\wh S(U^p,E)$. For any place $v$ of $F^+$ that splits as $v = w w^c$ in $F$, we choose an isomorphism $\iota_w : \u G(F_v^+) \congto \GL_n(F_w)$ that is conjugate to the isomorphism induced by $\iota_G$ on $F_w$-points (the choice of which won't matter). If in addition $U_v$ is a maximal compact subgroup of $\u G(F_v^+)$, then we demand moreover that $\iota_w(U_v) = \GL_n(\O_{F_w})$. We let $\Sigma(U^p)$ denote the (finite) set of places $v$ of $F^+$ that split in $F$ and are such that $U_v$ is not maximal compact. We let $\T(U^p) = E[T_w^{(j)}]$ denote the polynomial algebra over $E$ generated by all $T_w^{(j)}$ for $w$ a place of $F$ lying over a place $v$ of $F^+$ that splits in $F$ and such that $v \not\in \Sigma(U^p)$ and $v \nmid p$. Then $\T(U^p)$ acts topologically on $\wh
S(U^p,E)$ by letting $T_w^{(j)}$ act as the double coset operator $[U_v\, \iota_w^{-1} \big(\begin{smallmatrix} 1_{n-j} \\ &\varpi_w 1_j\end{smallmatrix}\big) U_v]$, where $\varpi_w$ is a uniformiser of $F_w$. This action commutes with the unitary left action of $\u G(F_p^+)$ on $\wh S(U^p,E)$ by right translation of functions, where $F_p^+ := F^+ \otimes_{\Q} \qp$. For each place $v|p$ of $F^+$ we choose a place $\tv|v$ of $F$ (the choice of which won't matter). We note that $F_p^+ \cong \prod_{v|p} F_\tv$ and that $\u G(F_p^+)$ is identified with $\prod_{v|p} \GL_n(F_\tv)$ via the $\iota_\tv$.

We let $r : \Gal(\o F/F) \to \GL_n(E)$ be a continuous representation and we assume:
\begin{enumerate}
\item $r$ is unramified at the places of $F^+$ that split in $F$ and are not in $\Sigma(U^p)$;
\item $r^c\cong r^\vee\otimes\varepsilon^{1-n}$ (where $r^c(g):=r(cgc)$, $g\in \Gal(\o F/F)$ and $r^\vee$ is the dual of $r$);
\item $r$ is an absolutely irreducible representation of $\Gal(\o F/F)$.
\end{enumerate}
We associate to $r$ and $U^p$ the maximal ideal $\m_r$ in $\T(U^p)$ generated by all elements
$$\Big((-1)^j{\rm Norm}(w)^{j(j-1)/2}T_w^{(j)}-a^{(j)}_w\Big)_{j,w},$$
where $j\in \{1,\dots,n\}$, $w$ is a place of $F$ lying over a place $v$ of $F^+$ that splits in $F$ and such that $v \not\in \Sigma(U^p)$ and $v \nmid p$, ${\rm Norm}(w)$ is the cardinality of the residue field at $w$, and where $X^n+a_w^{(1)}X^{n-1}+\cdots + a_w^{(n-1)}X+a_w^{(n)}$ is the characteristic polynomial of $r(\Frob_w)$ (an element of $\O[X]$, $\Frob_w$ is a geometric Frobenius element at $w$). We assume that $\wh S(U^p,E)\lan[\qp][\m_r]$ has non-zero locally $\qp$-{\it algebraic} vectors (where $(-)\lan[\qp]$ is the locally convex subspace of locally $\qp$-analytic vectors for the action of $\u G(F_p^+)$), i.e.\ $r$ is automorphic of ``level'' $U^p$.

We now assume in the following that for each place $v | p$ of $F^+$ the representation $r_\tv := r|_{\Gal(\o {F_\tv}/F_\tv)}$ is crystabelline and satisfies Hypothesis \ref{hyp:generic-crystabelline}. We use the notation of \S\ref{sec:local-setup-crystab} modified in a trivial way as follows: we let $\W := \prod_{v|p} \W_\tv$, $\Ws(r) := \prod_{v|p} \Ws_\tv(r_{\tv})$ and we write $C_\tv(w\alg_\tv, \F_\tv)$, $M_\tv(w\alg_\tv, \F_\tv)$, etc.\ for the representations of $\GL_n(F_\tv)$ defined in \S\ref{sec:local-setup-crystab}. Usually we will omit the subscript $\tv$ on the outside, which should not lead to any confusion.

Recall the following conjecture of the first author (slightly generalised, since we do not assume that $p$ splits completely in $F^+$).

\begin{conj}\cite[Conj.~6.1]{socle2}\label{conj:loc-an-socle}
Suppose that $(w\alg_\tv, \F_\tv)_v \in \W$. We have
$$\Hom_{\u G(F^+_p)}(\whotimes_{v|p}C(w\alg_\tv, \F_\tv)(\ve^{n-1}), \wh S(U^p,E)\lan[\qp][\m_r]) \ne 0$$
if and only if $(w\alg_\tv, \F_\tv)_v \in \Ws(r)$ \(where $*(\chi):=*\otimes (\chi\circ{\det})$\).
\end{conj}

Recall this is known to be true in many cases where $r_\tv$ is crystalline for all $v|p$.

\begin{thm}\cite[Thm.~1.3]{BHS3}\label{BHS3}
Assume that $U^p$ is sufficiently small \(\cite[\S3.3]{CHT}\), the residual representation $\overline r$ is absolutely irreducible, $r_\tv$ is crystalline for all $v|p$ \(and satisfies Hypothesis \ref{hyp:generic-crystabelline}\), and that the following assumptions hold:
\begin{enumerate}
\item[(i)]$p>2$;
\item[(ii)]$F/F^+$ is unramified and $\u G$ is quasi-split at all finite places of $F^+$;
\item[(iii)]$U_v$ is hyperspecial when the finite place $v$ of $F^+$ is inert in $F$;
\item[(iv)]${\overline r}({\rm Gal}(\overline F/F(\sqrt[p]{1}))$ is adequate \(\cite[Def.~2.3]{Thor}\).
\end{enumerate}
Then Conjecture \ref{conj:loc-an-socle} is true.
\end{thm}

Before stating our theorem, we observe that for any $(w\alg_\tv, \F_\tv)_v \in \W$ the representation
$\whotimes_{v|p} C(w\alg_\tv, \F_\tv)(\ve^{n-1})$ of $\u G(F^+_p) \cong \prod_{v|p} \GL_n(F_\tv)$ is
admissible and topologically irreducible. This follows from Lemma~\ref{lm:OS-product-groups} and
the theory of Orlik--Strauch. Below we will tacitly use the exactness of $\whotimes$ for compact type spaces (Corollary~\ref{cor:exactness}) and that the external tensor product of admissible locally analytic representations is admissible
(Lemma~\ref{lm:tensor-product-admissible}).

\begin{thm}\label{thm:adjunct-max-subrep-ps-coho}
Suppose that Conjecture~\ref{conj:loc-an-socle} holds. If $C_\tv \in \Cs(r_\tv)$ for all $v | p$, then any injective $\u G(F^+_p)$-equivariant homomorphism $\whotimes_{v|p} C_\tv(\ve^{n-1}) \to \wh S(U^p,E)\lan[\qp][\m_r]$ extends uniquely to a $\u G(F^+_p)$-equivariant homomorphism
$$\whotimes_{v|p} \Pi(r_\tv)\fs_{C_\tv}(\ve^{n-1}) \longrightarrow \wh S(U^p,E)\lan[\qp][\m_r]$$
that is moreover injective.
\end{thm}
\begin{proof}
We let $f$ denote the given injection $\whotimes_{v|p} C_\tv(\ve^{n-1}) \hookrightarrow \wh S(U^p,E)\lan[\qp][\m_r]$ and we use the notation of \S\ref{sec:local-setup-crystab}.

\emph{Step 1:} We show that $f$ extends uniquely to a $\u G(F^+_p)$-equivariant homomorphism
\begin{equation}\label{eq:19}
 \whotimes_{v|p} M(r_\tv)_{C_\tv}(\ve^{n-1}) \longrightarrow \wh S(U^p,E)\lan[\qp][\m_r].
\end{equation}
By construction (see Definition \ref{our-amalgam}(ii)), and by Lemma~\ref{lm:compl-tensor-finite-colimits}, it suffices to show that for any fixed $(w\alg_\tv, \F_\tv)_v \in \Ws(r)$ such that $C_\tv \cong C(w\alg_\tv, \F_\tv)$ for all $v|p$ the map $f$ extends uniquely to a $\u G(F^+_p)$-equivariant homomorphism
$$\whotimes_{v|p} M(w\alg_\tv, \F_\tv)(\ve^{n-1}) \longrightarrow \wh S(U^p,E)\lan[\qp][\m_r].$$
This follows by applying Proposition \ref{prop:adjunct-max-subrep-ps-abstract} for the group $\u G(F^+_p)$, $\Pi = \wh S(U^p,E)$, $\sigma = \whotimes_{v|p} M(w\alg_\tv, \F_\tv)(\ve^{n-1})$, $\T = \T(U^p)$, $\m = \m_r$, as we now explain. Note first that assumption (i) in the proposition is satisfied by the proof of~\cite[Prop.\ 6.3.3]{breuil-ext1} (where it is reduced to Example~\ref{ex:cohomology}). Next, $\whotimes_{v|p} M(w\alg_\tv, \F_\tv)$ is of finite length (by the remark preceding the theorem). As $C_\tv$ only occurs once in $M(w\alg_\tv, \F_\tv)$ for all $v | p$, it easily follows from Lemma~\ref{lm:socle-for-product-group} that the $\u G(F^+_p)$-socle of $\wh \otimes_{v|p} M(w\alg_\tv, \F_\tv)$ is isomorphic to $\whotimes_{v|p} C_\tv$. Any irreducible constituent of $\whotimes_{v|p} M(w\alg_\tv, \F_\tv)/\whotimes_{v|p} C_\tv$ has the form $\whotimes_{v|p} C(\tau_\tv, \F_\tv)$, where $\tau_\tv \ge w\alg_\tv$ and for at least one $v'|p$ we have $\tau_{\tv'} \nleq w\alg_{\tv'}(\F_{\tv'})$. Using Lemma~\ref{lm:OS-product-groups} we have $\whotimes_{v|p} C(\tau_\tv, \F_\tv) \cong \F_{\u B(F_p^+)}^{\u G(F_p^+)}(L(\mu),\chi)$ with $\mu = (\tau_\tv \cdot (-\lambda_{\tv}))_{v|p}$ and $\chi = \prod_{v|p} \pi_{B_\tv,\F_\tv}$. By Lemma \ref{lm:OS-product-groups} again, $\F_{\u B(F_p^+)}^{\u G(F_p^+)}(M(\mu),\chi) \cong \whotimes_{v|p} \PS(\tau_\tv,\F_\tv)$ which has irreducible constituents $\whotimes_{v|p} C(\tau'_\tv, \F_\tv)$ with $\tau'_\tv \ge \tau_\tv$ for all $v | p$, so $\tau'_{\tv'} \nleq w\alg_{\tv'}(\F_{\tv'})$. Therefore, by Conjecture~\ref{conj:loc-an-socle}, we may indeed apply Proposition \ref{prop:adjunct-max-subrep-ps-abstract} and the claim follows. 

\emph{Step 2:} We show that the map~(\ref{eq:19}) factors uniquely through a map
\begin{equation}
\whotimes_{v|p} \Pi(r_\tv)\fs_{C_\tv}(\ve^{n-1}) \longrightarrow \wh S(U^p,E)\lan[\qp][\m_r]\label{eq:26}
\end{equation}
which is moreover injective.
By \ Corollary~\ref{cor:exactness} \ we \ have \ a \ surjection \ $\whotimes_{v|p} M(r_\tv)_{C_\tv} \onto \whotimes_{v|p} \Pi(r_\tv)\fs_{C_\tv}$ \ with \ kernel $\sum_{v|p} K_\tv \whotimes (\whotimes_{v'\ne v} M(r_{\tv'})_{C_{\tv'}})$, where $K_\tv := \ker(M(r_\tv)_{C_\tv} \onto \Pi(r_\tv)\fs_{C_\tv})$. By definition, no irreducible constituent of $K_\tv$ is contained in $\Cs(r_\tv)$, hence the claim follows from Conjecture~\ref{conj:loc-an-socle}. The resulting map~\eqref{eq:26} is injective, since it is non-zero, the left-hand side has socle $\whotimes_{v|p} C_\tv$ by Lemma~\ref{lm:socle-for-product-group}, and no other irreducible constituent injects into the right-hand side (by Conjecture~\ref{conj:loc-an-socle}).
\end{proof}

\begin{lm}\label{lm:fd}
For $(w\alg_\tv, \F_\tv)_v \in \W$ the $E$-vector space
\begin{equation*}
 \Hom_{\u G(F_p^+)}\big(\whotimes_{v|p} C(w\alg_\tv, \F_\tv)(\ve^{n-1}),\wh S(U^p,E)\lan[\qp][\m_r]\big)
\end{equation*}
is finite-dimensional.
\end{lm}
\begin{proof}
Let $\Pi := \wh S(U^p,E)\lan[\qp][\m_r]$, which is an admissible locally $\qp$-analytic $\u G(F_p^+)$-representation. As in (\ref{eq:6}) and using the notation of \S\ref{sec:partial-adjunction} (recall $\u G(F_p^+)\cong \prod_{v\vert p}\GL_n(F_\tv)$), we can find a locally $\qp$-analytic character $\chi : T \to E\s$ and a non-zero $T^+$-equivariant homomorphism $f : \chi \to C^{\o N_0}$, where $C:=\whotimes_{v|p} C(w\alg_\tv, \F_\tv)$. Then restriction to $\o N_0$-invariants and composition with $f$ gives a map
\begin{equation*}
\Hom_{\u G(F_p^+)}\big(\whotimes_{v|p} C(w\alg_\tv, \F_\tv)(\ve^{n-1}),\Pi\big) \longrightarrow \Hom_{T^+}(\theta,\Pi^{\o N_0}),
\end{equation*}
(where $\theta$ is the relevant twist of $\chi$) which is injective, as $\whotimes_{v|p} C(w\alg_\tv, \F_\tv)$ is irreducible. The latter space is finite-dimensional by Corollary~\ref{cor:breuil-633}.
\end{proof}

\begin{coroll}\label{cor:adjunct-max-subrep-ps-coho}
We keep the hypotheses of Theorem~\ref{thm:adjunct-max-subrep-ps-coho}. For each $C = \wh \otimes_{v|p} C_\tv$ with $C_\tv \in \Cs(r_\tv)$ let
$$n_C := \dim_E \Hom_{\u G(F_p^+)}\big(\wh \otimes_{v|p} C_\tv(\ve^{n-1}),\wh S(U^p,E)\lan[\qp][\m_r]\big) \in \Z_{> 0}.$$
Then there exists an injective $\u G(F_p^+)$-equivariant linear map
\begin{equation}\label{eq:20}
\bigoplus_{C = \whotimes C_\tv} \big(\whotimes_{v|p} \Pi(r_\tv)\fs_{C_\tv}(\ve^{n-1})\big)^{\oplus n_C} \longrightarrow \wh S(U^p,E)\lan[\qp][\m_r].
\end{equation}
\end{coroll}
\begin{proof}
By assumption we have an injection
\begin{equation}\label{injsocle}
\bigoplus_{C = \whotimes C_\tv} \big(\whotimes_{v|p} C_\tv(\ve^{n-1})\big)^{\oplus n_C} \hookrightarrow \wh S(U^p,E)\lan[\qp][\m_r].
\end{equation}
Applying Theorem~\ref{thm:adjunct-max-subrep-ps-coho} to each irreducible direct summand, we see that the given map extends uniquely to a map as in~(\ref{eq:20}). The resulting map is injective because from (\ref{injsocle}) it is injective on the socle.
\end{proof}

Combining Corollary \ref{cor:adjunct-max-subrep-ps-coho} with Theorem \ref{BHS3}, we obtain the result in the introduction.

\subsection{Special cases and examples}\label{sec:examples}

We give explicit examples for the representations $\Pi(\rho)\fs$ and also relate $\Pi(\rho)\fs$ to previous results or conjectures. For simplicity, we only consider here crystalline representations.\\

We first give two examples in the crystalline case for $\GL_3(\qp)$. A refinement is here an ordering of the (distinct) eigenvalues $\{\varphi_1,\varphi_2,\varphi_3\}$ of the crystalline Frobenius. We denote by $s_\alpha$, $s_\beta$ the two simple reflections, which generate the Weyl group $W={S}_3$.

We start with the {\it noncritical} case, by which we mean $w\alg(\F)=1$ for all refinements $\F$. We fix an arbitrary refinement $\F_0:=(\varphi_1,\varphi_2,\varphi_3)$ (the choice of which won't matter) and recall that $w\F_0=(\varphi_{w^{-1}(1)},\varphi_{w^{-1}(2)},\varphi_{w^{-1}(3)})$ for $w\in {S}_3$. One can then check that $\Pi(\rho)\fs$ has the following explicit form, where the constituent $C(w\alg,w\F_0)$ is just denoted $C_{w\alg,w}$ below, where the (irreducible) socle is the constituent $C(1,\F_0)=C_{1,1}$ in the {\it middle}, where we use without comment the intertwinings provided by Lemma \ref{lm:princ-series}(ii), where a line between two constituents means as usual a non-split extension as subquotient and where the constituent further away from the centre is the quotient:
$$
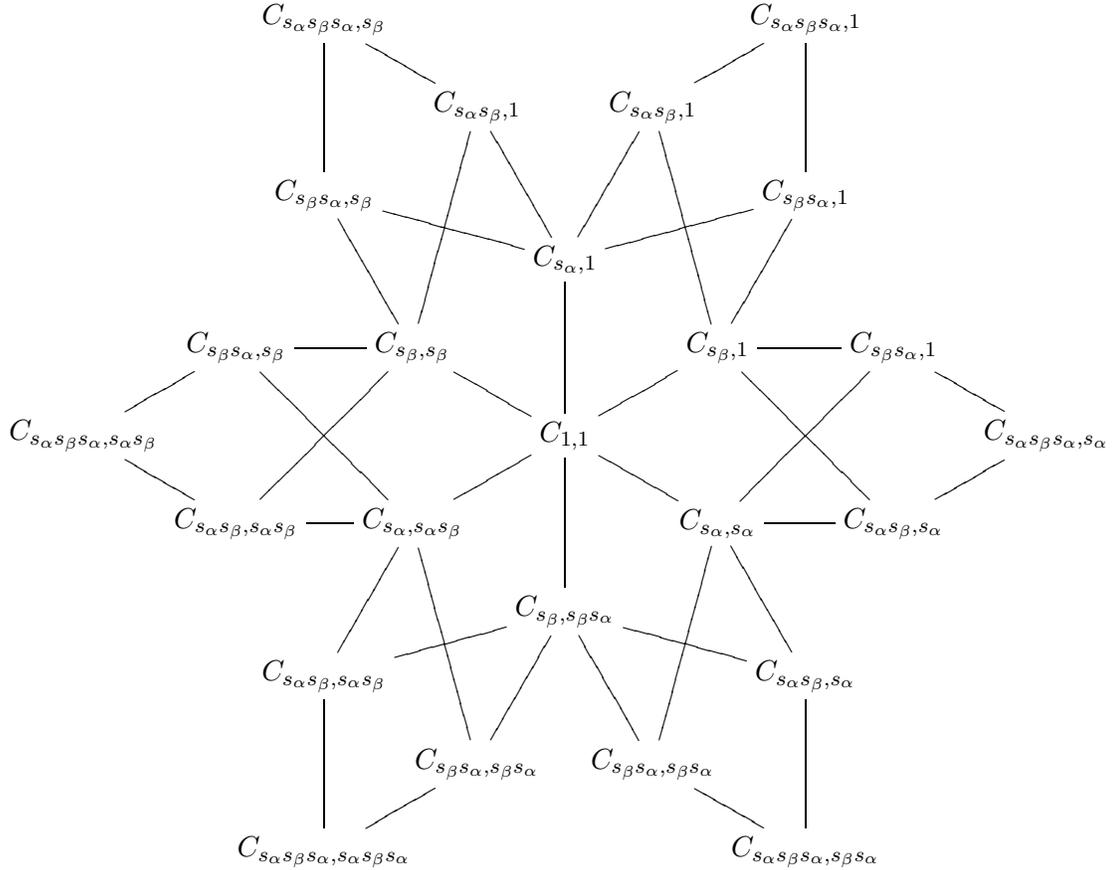
\begin{xy}
(0,0)*+{C_{1,1}}="a"; (0,23.4)*+{C_{s_\alpha,1}}="b"; (20.3,11.7)*+{C_{s_\beta,1}}="c"; (20.3,-11.7)*+{C_{s_\alpha,s_\alpha}}="d";(0,-23.4)*+{C_{s_\beta,s_\beta s_\alpha}}="e";(-20.3,11.7)*+{C_{s_\beta,s_\beta}}="f"; (-20.3,-11.7)*+{C_{s_\alpha,s_\alpha s_\beta}}="g"; 
(11.7,43.7)*+{C_{s_\alpha s_\beta,1}}="h"; (32,32)*+{C_{s_\beta s_\alpha,1}}="i";(32,55.4)*+{C_{s_\alpha s_\beta s_\alpha,1}}="j";
(-11.7,43.7)*+{C_{s_\alpha s_\beta,1}}="k"; (-32,32)*+{C_{s_\beta s_\alpha,s_\beta}}="l";(-32,55.4)*+{C_{s_\alpha s_\beta s_\alpha,s_\beta}}="m";
(43.7,11.7)*+{C_{s_\beta s_\alpha,1}}="n";(43.7,-11.7)*+{C_{s_\alpha s_\beta,s_\alpha}}="o";(64,0)*+{C_{s_\alpha s_\beta s_\alpha,s_\alpha}}="p";
(11.7,-43.7)*+{C_{s_\beta s_\alpha,s_\beta s_\alpha}}="q"; (32,-32)*+{C_{s_\alpha s_\beta,s_\alpha}}="r";(32,-55.4)*+{C_{s_\alpha s_\beta s_\alpha,s_\beta s_\alpha}}="s";
(-11.7,-43.7)*+{C_{s_\beta s_\alpha,s_\beta s_\alpha}}="t"; (-32,-32)*+{C_{s_\alpha s_\beta,s_\alpha s_\beta}}="u";(-32,-55.4)*+{C_{s_\alpha s_\beta s_\alpha,s_\alpha s_\beta s_\alpha}}="v";
(-43.7,-11.7)*+{C_{s_\alpha s_\beta,s_\alpha s_\beta}}="w";(-43.7,11.7)*+{C_{s_\beta s_\alpha,s_\beta}}="x";(-64,0)*+{C_{s_\alpha s_\beta s_\alpha,s_\alpha s_\beta}}="y";
{\ar@{-}"a";"b"}; {\ar@{-}"a";"e"}; {\ar@{-}"a";"c"}; {\ar@{-}"a";"d"}; {\ar@{-}"a";"f"}; {\ar@{-}"a";"g"};{\ar@{-}"b";"h"} ;{\ar@{-}"c";"i"};{\ar@{-}"b";"i"} ;{\ar@{-}"c";"h"};{\ar@{-}"h";"j"};{\ar@{-}"i";"j"};{\ar@{-}"b";"k"};{\ar@{-}"f";"l"};{\ar@{-}"f";"k"};{\ar@{-}"b";"l"};{\ar@{-}"l";"m"};{\ar@{-}"k";"m"};{\ar@{-}"c";"n"};{\ar@{-}"d";"o"};{\ar@{-}"c";"o"};{\ar@{-}"d";"n"};{\ar@{-}"o";"p"};{\ar@{-}"n";"p"};{\ar@{-}"d";"r"};{\ar@{-}"d";"q"};{\ar@{-}"e";"r"};{\ar@{-}"e";"q"};{\ar@{-}"q";"s"};{\ar@{-}"r";"s"};{\ar@{-}"e";"t"};{\ar@{-}"e";"u"};{\ar@{-}"g";"t"};{\ar@{-}"g";"u"};{\ar@{-}"t";"v"};{\ar@{-}"u";"v"};{\ar@{-}"f";"w"};{\ar@{-}"g";"w"};{\ar@{-}"f";"x"};{\ar@{-}"g";"x"};{\ar@{-}"x";"y"};{\ar@{-}"w";"y"}
\end{xy}$$
(Note with Lemma \ref{lm:princ-series}(ii) that $\Pi(\rho)\fs$ is {\it not} multiplicity free: the $6$ distinct constituents of the form $C(s_\alpha s_\beta,\F)$ or $C(s_\beta s_\alpha,\F)$ all appear with multiplicity $2$.)

We go on with an example in the critical ordinary case. Here we have a canonical refinement $\F_0=\F_\rho$ (see  the beginning of \S\ref{sec:local-setup-ord}) due to the fact that $\rho$ is upper triangular with distinct Hodge--Tate weights. The possible locally analytic socles are worked out in \cite[\S6.2]{Breuil-peking}, we only give here $\Pi(\rho)\fs$ when its socle is $C(1,\F_\rho)\oplus C(s_\alpha,s_\alpha s_\beta\F_\rho)=C_{1,1}\oplus C_{s_\alpha,s_\alpha s_\beta}$ (the interested reader can easily work out the other cases). We get the following form (same notation as before, the socle in each summand being now on the {\it left})
$$\!\!\!
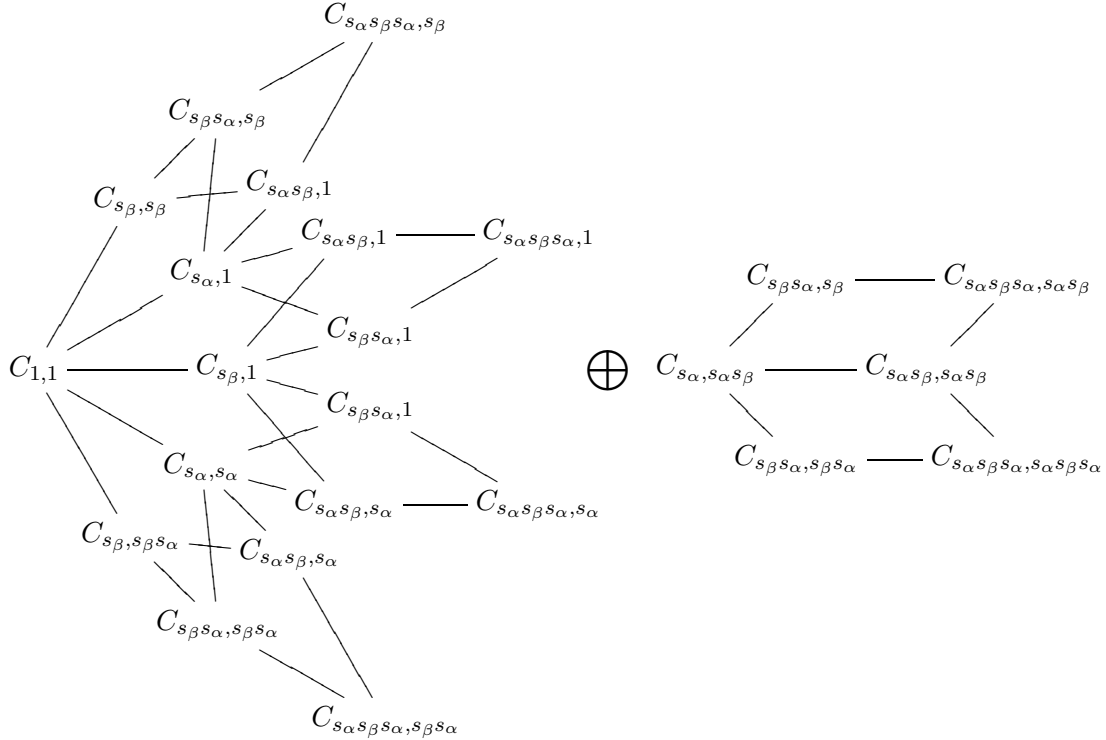
\begin{xy}
(0,0)*+{C_{1,1}}="a"; 
(12.9,22.3)*+{C_{s_\beta,s_\beta}}="e"; (22.3,12.9)*+{C_{s_\alpha,1}}="b"; 
(24.5,34)*+{C_{s_\beta s_\alpha,s_\beta}}="g"; (34,24.5)*+{C_{s_\alpha s_\beta ,1}}="h"; (46.8,46.8)*+{C_{s_\alpha s_\beta s_\alpha, s_\beta}}="i";
(25.7,0)*+{C_{s_\beta,1}}="c"; 
(41.3,18)*+{C_{s_\alpha s_\beta ,1}}="j"; (44.7,5)*+{C_{s_\beta s_\alpha,1}}="k"; 
(67.1,18)*+{C_{s_\alpha s_\beta s_\alpha,1}}="l";
(22.3,-12.9)*+{C_{s_\alpha,s_\alpha}}="d"; (12.9,-22.3)*+{C_{s_\beta,s_\beta s_\alpha}}="f";
(41.3,-18)*+{C_{s_\alpha s_\beta ,s_\alpha}}="m"; (44.7,-5)*+{C_{s_\beta s_\alpha,1}}="n"; 
(67.1,-18)*+{C_{s_\alpha s_\beta s_\alpha,s_\alpha}}="o";
(24.5,-34)*+{C_{s_\beta s_\alpha,s_\beta s_\alpha}}="p"; (34,-24.5)*+{C_{s_\alpha s_\beta ,s_\alpha}}="q"; (46.8,-46.8)*+{C_{s_\alpha s_\beta s_\alpha, s_\beta s_\alpha}}="r";
{\ar@{-}"a";"b"}; {\ar@{-}"a";"c"}; {\ar@{-}"a";"c"}; {\ar@{-}"a";"d"}; {\ar@{-}"a";"e"}; {\ar@{-}"a";"f"} ; {\ar@{-}"e";"g"}; {\ar@{-}"e";"h"};{\ar@{-}"b";"g"}; {\ar@{-}"b";"h"};{\ar@{-}"g";"i"} ; {\ar@{-}"h";"i"};
{\ar@{-}"b";"j"};{\ar@{-}"b";"k"}; {\ar@{-}"c";"j"};{\ar@{-}"c";"k"};
{\ar@{-}"j";"l"};{\ar@{-}"k";"l"};
{\ar@{-}"d";"m"};{\ar@{-}"d";"n"}; {\ar@{-}"c";"m"};{\ar@{-}"c";"n"};
{\ar@{-}"m";"o"};{\ar@{-}"n";"o"};
{\ar@{-}"d";"p"};{\ar@{-}"d";"q"}; {\ar@{-}"f";"p"};{\ar@{-}"f";"q"};
{\ar@{-}"p";"r"};{\ar@{-}"q";"r"}
\end{xy}
\!\!\!\!\!\!\bigoplus\ \ 
\begin{xy}
(0,0)*+{C_{s_\alpha,s_\alpha s_\beta}}="s";(12,12)*+{C_{s_\beta s_\alpha,s_\beta}}="t";(29.4,0)*+{C_{s_\alpha s_\beta,s_\alpha s_\beta}}="u";(41.4,12)*+{C_{s_\alpha s_\beta s_\alpha,s_\alpha s_\beta}}="v";(12,-12)*+{C_{s_\beta s_\alpha,s_\beta s_\alpha}}="w";(41.4,-12)*+{C_{s_\alpha s_\beta s_\alpha,s_\alpha s_\beta s_\alpha}}="x";
{\ar@{-}"s";"t"}; {\ar@{-}"s";"u"};{\ar@{-}"s";"w"};
{\ar@{-}"t";"v"};{\ar@{-}"u";"v"};{\ar@{-}"u";"x"};{\ar@{-}"w";"x"}
\end{xy}
$$
One can check that $\Pi(\rho)\fs$ again fails to be multiplicity free ($4$ constituents appear with multiplicity $2$) and that $(\Pi(\rho)\ord)\lan[\qp]$ (see Proposition \ref{prop:Qrho-ordinary}) is the direct summand on the left.

In the crystalline case for $\GL_2(L)$, the representation $\Pi(\rho)\fs$ is easily checked to be exactly the representation $\Pi(D_{\rm cris}(\rho))$ in \cite[\S4(9)]{socle1}. In particular, in this case Theorem \ref{thm:adjunct-max-subrep-ps-coho} was already proven by Ding in the setting of the completed $H^1$ of unitary Shimura curves (see \cite[Thm.~6.3.7]{Ding}). The proof of {\it loc.cit.} however is different from that of Theorem \ref{thm:adjunct-max-subrep-ps-coho} (e.g.\ it doesn't use \cite[Thm.\ 7.1.1]{breuil-ext1}). Note that here $\Pi(\rho)\fs$ is multiplicity free.\\

Finally, Theorem \ref{thm:adjunct-max-subrep-ps-coho} (assuming Conjecture \ref{conj:loc-an-socle}) together with Proposition \ref{prop:constituents-of-Qrho}(i) imply that any constituent of the form
$$\big(\whotimes_{v'|p,v'\ne v}C(1,\F_\tvprime)(\ve^{n-1})\big)\whotimes C(s_{\alpha_\tv}, \F_\tv)(\ve^{n-1}),$$
where $v\vert p$ and $s_{\alpha_\tv}$ is a {\it simple} reflection in $\Res_{F_\tv/\qp} \GL_n \times_{\qp} E \cong \prod_{\cS_{F_\tv}} \GL_{n/E}$ (see \S\ref{sec:local-setup-crystab}), that {\it does not} inject into $\wh S(U^p,E)\lan[\qp][\m_r]$ is such that there is a non-split extension
$$\big(\whotimes_{v'|p,v'\ne v}C(1,\F_\tvprime)(\ve^{n-1})\big)\whotimes \big(C(1,\F_\tv) \text{---} C(s_{\alpha_\tv}, \F_\tv)\big)(\ve^{n-1})$$
that {\it does} inject into $\wh S(U^p,E)\lan[\qp][\m_r]$. By \cite[\S3.3]{breuil-ext1} and together with Theorem \ref{BHS3}, this gives further evidence to \cite[Conj.~6.1.1]{breuil-ext1} in the crystalline case (note that in {\it loc.cit.} it is assumed that there is only one $v$ dividing $p$ in $F^+$, in which case the factor $\whotimes_{v'|p,v'\ne v}C(1,\F_\tvprime)(\ve^{n-1})$ disappears).

\section{Ordinary representations}\label{sec:ordin-repr}

For $L=\qp$ and $\rho$ crystabelline upper triangular satisfying Hypothesis~\ref{hyp:generic-crystabelline} we prove that the locally analytic vectors of the representation $\Pi(\rho)\ord$ of \cite[\S3.3]{fondamental} is a subrepresentation of $\Pi(\rho)\fs$, and then deduce strong evidence to \cite[Conj.~4.2.2]{fondamental} in the crystalline case using Theorem \ref{cor:adjunct-max-subrep-ps-coho} (and Theorem \ref{BHS3}).

\subsection{Local setup and results}\label{sec:local-setup-ord}

For $L=\qp$ and $\rho$ crystabelline upper triangular satisfying Hypothesis~\ref{hyp:generic-crystabelline} we prove (among other results) that $(\Pi(\rho)\ord)\lan[\qp]$ is a subrepresentation of $\Pi(\rho)\fs$ (Proposition \ref{prop:Qrho-ordinary}) and that $\Pi(\rho)\ord$ is its universal unitary completion (Proposition \ref{prop:univ-completion-of-analytic-vecs}).\\

We keep the notation of \S\ref{sec:fin-slope-gen-crystab} and specialise to the case where $L = \qp$ and $\rho$ is crystabelline ordinary, that is $\rho : \Gal(\qpb/\qp) \to \GL_n(E)$ is of the form
\begin{equation}\label{eq:rho}
 \rho \sim
 \left(\begin{matrix}
 \psi_1 & * & \dots & * \\ & \psi_2 & \dots & * \\ && \ddots & \vdots \\ &&& \psi_n
 \end{matrix}\right)
\end{equation}
and satisfies Hypothesis~\ref{hyp:generic-crystabelline}. We write $\psi_i(x) = x^{-h_i} \chi_i(x)$ (for $x\in \qp\s$) with $\chi_i = 
\WD(\psi_i)$ smooth and $h_i \in \Z$. We remark that Hypothesis~\ref{hyp:generic-crystabelline} implies that $\rho$ is generic in the sense of \cite[Def.\ 3.3.1]{fondamental}.

As $\rho$ is regular de Rham we may assume without loss of generality that $h_1 < \cdots < h_n$. Thus $D := D\cris(\rho)$ has a canonical refinement $\F_\rho$ with $\chi_{\F_\rho,i} = \chi_i$ for all $i$, and for $w\alg,w \in W$ we write $\eta(w\alg,w) := \eta(w\alg,w\F_\rho)$, $\PS(w\alg,w) := \PS(w\alg,w\F_\rho)$, $C(w\alg,w) := C(w\alg,w\F_\rho)$, and $w\alg(w) := w\alg(w
\F_\rho)$.

We fix a representative homomorphism in the conjugacy class $\rho$ that is a \emph{good conjugate} in the sense of \cite[Def.\ 3.2.4]{fondamental}, and we will also denote it by $\rho$. This is possible after conjugating by a suitable element of $\u B(E)$ by \cite[Prop.\ 3.2.3]{fondamental}. We emphasise that the following definition depends on our choice of good conjugate.

\begin{df}\label{compareBH}
We let $W_\rho = \{ w \in W : \text{$\dot w \rho \dot w^{-1}$ is upper-triangular} \}$, where $\dot w$ is a representative of $w$ (this is the inverse of the subset $W_{C_\rho}$ defined in \cite[\S3.2]{fondamental} and is in general {\it different} from the subset denoted by $W_\rho$ in \cite[(14)]{fondamental}). For each $w \in W_\rho$ we let $\Sigma_w \subseteq W$ consist of all (commuting!)\ products $s_{\alpha_1} \cdots s_{\alpha_r}$ with $r \ge 0$, $\alpha_i \in \Delta$ pairwise orthogonal, and $s_{\alpha_i} w \not \in W_\rho$ for all $i$.
\end{df}

Note that if $\sigma \in \Sigma_w$ and $\sigma' \le \sigma$, then $\sigma' \in \Sigma_w$. (In fact, $\sigma' = \prod_{i \in I} s_{\alpha_i}$ for a unique subset $I \subseteq \{1,\dots,r\}$.) Note also that in this case we can uniquely write $\sigma = \sigma'' \sigma'$ with $\ell(\sigma) = \ell(\sigma'')+\ell(\sigma')$. (Namely, $\sigma'' = \prod_{i \not\in I} s_{\alpha_i}$.)

\begin{lm}\label{lm:walg}
For any $w \in W_\rho$, $\sigma \in \Sigma_w$ we have $w\alg(\sigma w) = w$.
\end{lm}
\begin{proof}
For any $0 \le i \le n$, let $\F_{w,i}$ denote the member of the flag $w\F_{\rho}$ that has rank $i$. The element $w\alg(w)w_0$ gives the relative position of the flags $(w\F_{\rho})_{L'}^{\Gal(L'/\qp)}$ and $(\Fil^*D_{L'})^{\Gal(L'/\qp)}$. A calculation shows that this means that the filtration $\Fil^*D_{L'} \cap (\F_{w,i})_{L'}$ on $(\F_{w,i})_{L'}$ (omitting Galois invariants for simplicity) jumps precisely at the integers $\{h_{w\alg(w)^{-1}(j)} : 1 \le j \le i\}$.

Next we claim that for $w \in W$ the flag $w\F_{\rho}$ (with induced structures) consists of weakly admissible subobjects if and only if $w \in W_\rho$. The first condition is equivalent to saying that $\rho$ has a filtration with subquotients $\psi_{w^{-1}(i)}$, $1 \le i \le n$ in this order (where $i = 1$ corresponds to the subobject), i.e.\
\begin{equation}\label{eq:5}
\rho \sim
 \left(\begin{matrix}
 \psi_{w^{-1}(1)} & \dots & * \\ & \ddots & \vdots \\ && \psi_{w^{-1}(n)}
\end{matrix}\right).
\end{equation}
By \cite[Prop.\ 3.2.3]{fondamental} we may assume, after further conjugation by $\o B(E)$, that the right-hand side of~\eqref{eq:5} is a good conjugate. By \cite[Prop.\ 3.2.6]{fondamental} it then follows that $w \in W_\rho$. The converse is clear.

Let us now go back to our given $w \in W_\rho$ and $\sigma \in \Sigma_w$. We can write $\sigma =\prod_{\ell=1}^r s_\ell$ with $s_\ell$ corresponding to the simple root $\epsilon_{n_\ell}-\epsilon_{n_\ell+1}$ and such that $n_\ell+ 1 < n_{\ell+1}$ for all $1 \le \ell< r$. Assume first $\sigma=1$. By the previous paragraphs each subobject $\cF_{w,i}$ is weakly admissible and the Hodge filtration of $\cF_{w,i}$ jumps at $\{h_{w\alg(w)^{-1}(j)} : 1 \le j \le i\}$, from which we easily deduce $w\alg(w) = w$ by induction. For general $\sigma\in \Sigma_w$ as above, note that $\cF_{\sigma w,i} = \cF_{w,i}$ for $i \not \in \{n_1,\dots,n_r\}$ and that $\cF_{\sigma w,n_\ell} = \cF_{s_\ell w,n_\ell}$, which is not weakly admissible by the previous paragraph, as $s_\ell w \not\in W_\rho$. As moreover $\cF_{\sigma w,i} = \cF_{w,i}$ for $i = n_\ell \pm 1$, the only possibility is that the Hodge filtration of $\cF_{\sigma w,n_\ell}$ jumps at $\{ h_{w\alg(w)^{-1}(j)} : 1 \le j \le n_\ell \}$. Hence again $w\alg(\sigma w) = w$.
\end{proof}

We recall the following result of Breuil--Emerton (\cite[Thm.~2.2.2]{BE}).

\begin{prop}\label{prop:breuil-emerton}
Suppose that $n = 2$, that $k_1 > k_2$ are integers, and that $\theta_1,\theta_2 : \qp\s \to E\s$ are smooth characters such that $x\mapsto \theta_i(x) (x)^{k_i}$ is unitary for $i = 1, 2$. If $k_1 = k_2+1$, we further assume that $\theta_1|\cdot|_{\qp}^{-1} \ne \theta_2$. Then the universal unitary completion $\wh\pi$ of
\begin{equation*}
\pi := \big(\Ind_{B(\qp)}^{G(\qp)} \theta_2 (-)^{k_1} \ve^{-1} \otimes \theta_1 (-)^{k_2} \big)\qpan
\end{equation*}
\(where $(-)^{k}$ means the character $x\in \qp^\times\mapsto x^k$\) is an admissible representation that is a non-split extension of $\big(\Ind_{B(\qp)}^{G(\qp)} \theta_2 (-)^{k_2} \ve^{-1} \otimes \theta_1 (-)^{k_1} \big)\cts$ by $\big(\Ind_{B(\qp)}^{G(\qp)} \theta_1 (-)^{k_1} \ve^{-1} \otimes \theta_2 (-)^{k_2} \big)\cts$, each of which is topologically irreducible. Moreover if $\theta_1 \theta_2^{\pm 1} \not\in \{ 1, |\cdot|_{\qp}^{\pm 1}\}$, then the canonical map $\pi \to \wh\pi$ is injective.
\end{prop}
\begin{proof}
We first reduce by twisting by a power of the unitary character $\ve$ to the case where $k_1 = 1$. Then the first result follows from \cite[Thm.\ 2.2.2]{BE}: in their notation we need to take $\chi_2 =\theta_2 (-)^{k_1} \ve^{-1}$ smooth, $\chi_1 = \theta_1 (-)^{k_2}$, and $k = k_1-k_2+1$. Note also that \cite{BE} work with the upper-triangular Borel $\o B$.

For \ the second \ result, \ recall \ that \ $\pi$ \ is \ a \ non-split \ extension \ of $\pi'' := \big(\Ind_{B(\qp)}^{G(\qp)} \theta_2 (-)^{k_2} \ve^{-1} \otimes \theta_1 (-)^{k_1} \big)\qpan$ by $\soc_{G(\qp)} \big(\Ind_{B(\qp)}^{G(\qp)} \theta_1 (-)^{k_1} \ve^{-1} \otimes \theta_2(-)^{k_2} \big)\qpan$, both of which are irreducible, as $\theta_1 \theta_2^{\pm 1} \not\in \{ 1,|\cdot|_{\qp}^{\pm 1}\}$ (see e.g.\ \cite[\S3.2]{Breuil-peking}). If the map $\pi \to \wh\pi$ isn't injective, it thus has to factor through the quotient $\pi''$. From the definition of universal unitary completions it would then follow that $\wh\pi \cong (\pi'')\sphat\,$. However, we know that $(\pi'')\sphat\, \cong (\Ind_{B(\qp)}^{G(\qp)} \theta_2 (-)^{k_2} \ve^{-1} \otimes \theta_1(-)^{k_1})\cts$ (for example by Proposition \ref{prop:univ-unit-compl}, but see also the proof of \cite[Thm.\ 2.2.2]{BE}). This contradiction shows that the map $\pi \to \wh\pi$ is indeed injective.
\end{proof}

\begin{prop}\label{prop:univ-completion-ordinary}
For any $w \in W_\rho$, $\sigma \in \Sigma_w$, the locally analytic principal series $\PS(w,\sigma w)$ admits an \emph{admissible} universal unitary completion $\PS(w,\sigma w)\sphat\,$ that is isomorphic to the representation $\Pi(\rho)_{w^{-1}(J)}$ constructed in \cite[\S 3.3]{fondamental}, where $J = \{ \alpha \in \Delta : s_\alpha \in \supp(\sigma) \}$. Moreover, the canonical map $\PS(w,\sigma w) \to \PS(w,\sigma w)\sphat\,$ is injective.
\end{prop}
\begin{proof}
Let $P_J$ be the parabolic containing $B$ determined by the subset $J \subseteq \Delta$ of the statement, let $M_J$ be the Levi subgroup of $P_J$ that contains $T$ and let $N_J$ denote the unipotent radical of $P_J$. We let $\pi:=\big(\Ind^{M_J(\qp)}_{(B \cap M_J)(\qp)}\eta(w,\sigma w)\big)\qpan$, so that $\PS(w,\sigma w) \cong \big(\Ind_{P_J(\qp)}^{G(\qp)} \pi\big)\qpan$.

From Proposition~\ref{prop:univ-unit-compl} applied to the character $\eta(w,\sigma w)$ we deduce that there exists a BH-subspace $\pi_0 \subseteq \pi$ such that $\pi = \sum_{m \in M_J(\qp)} m \pi_0$ and that $\wh\pi$ exists. Moreover observe that $\sigma$ lies in the Weyl group of $M_J$ with respect to $T$. Hence the central character $\eta(w,\sigma w)|_{Z_{M_J}(\qp)}$ of $\pi$ is equal to $\eta(w,w)|_{Z_{M_J}(\qp)}$, which is unitary by~(\ref{eq:23}), as each $\psi_i$ is. Applying \ Proposition~\ref{prop:univ-unit-compl} \ again, \ this \ time \ to \ $\pi$, \ we \ deduce \ that $\PS(w,\sigma w) \cong \big(\Ind_{P_J(\qp)}^{G(\qp)} \pi\big)\qpan$ has universal unitary completion $\big(\Ind_{P_J(\qp)}^{G(\qp)} \wh\pi\big)\cts$.

We now determine $\wh\pi$ explicitly. There exist integers $1 \le n_1 < \dots < n_r < n$ such that $n_i + 1 < n_{i+1}$ for all $i$, and $J$ consists of the simple roots $\epsilon_{n_i}-\epsilon_{n_i+1}$. Then we can identify $M_J$ with $\GL_2^r \times \GL_1^{n-2r}$, where the $i$-th factor of $\GL_2$ corresponds to the simple root $\epsilon_{n_i}-\epsilon_{n_i+1}$. By Lemma~\ref{lm:ind-product-groups}, $\pi$ is the external completed tensor product of all $\psi_{w^{-1}(j)} \ve^{-(n-j)}$ for $\sigma(j) = j$ and
\begin{equation}\small
 \big(\Ind_{B_2(\qp)}^{\GL_2(\qp)} \chi_{w^{-1}(n_i+1)} (-)^{-h_{w^{-1}(n_i)}} \ve^{-(n-n_i)} \otimes
 \chi_{w^{-1}(n_i)} (-)^{-h_{w^{-1}(n_i+1)}} \ve^{-(n-n_i-1)} \big)\qpan\label{eq:16}
\end{equation}
for $1 \le i \le r$, where $B_2$ is the lower triangular Borel of $\GL_2$. For any $1 \le i \le r$, as $w \in W_\rho$ and $s_{\alpha_i} w \not\in W_\rho$, we note that the homomorphism $\dot w \rho \dot w^{-1}$ contains the $2 \times 2$-submatrix $\left(\begin{matrix}\psi_{w^{-1}(n_i)} & * \\ & \psi_{w^{-1}(n_i+1)}\end{matrix}\right)$ with $* \ne 0$. As $\rho$ is by assumption a good conjugate, we see that the extension $*$ is in fact non-split. It is moreover de Rham, as $\rho$ is, hence $-h_{w^{-1}(n_i)} > -h_{w^{-1}(n_i+1)}$.

By Lemma~\ref{lm:universal-tensor} (both parts), Proposition~\ref{prop:breuil-emerton}, and Hypothesis~\ref{hyp:generic-crystabelline} (and what is above), we deduce that $\wh\pi$ is the external completed tensor product of all $\psi_{w^{-1}(j)} \ve^{-(n-j)}$ for $\sigma(j) = j$ and of the unique non-split extension of $\big(\Ind_{B_2(\qp)}^{\GL_2(\qp)} \psi_{w^{-1}(n_i+1)} \ve^{-(n-n_i)} \otimes \psi_{w^{-1}(n_i)} \ve^{-(n-n_i-1)} \big)\cts$ by $\big(\Ind_{B_2(\qp)}^{\GL_2(\qp)} \psi_{w^{-1}(n_i)} \ve^{-(n-n_i)} \otimes \psi_{w^{-1}(n_i+1)}\ve^{-(n-n_i-1)} \big)\cts$ for $1 \le i \le r$ (cf.\ \cite[Prop.\ B.2]{fondamental} for the uniqueness). This is admissible and isomorphic to the representation $\wt\Pi(\rho)_{w^{-1}(J)}$ constructed in \cite[\S 3.3]{fondamental} (see in particular the construction in Step 2 of the proof of \cite[Prop.\ 3.3.3]{fondamental}). By the second paragraph it follows that $\PS(w,\sigma w)\sphat\, \cong \big(\Ind_{P_J(\qp)}^{G(\qp)} \wt\Pi(\rho)_{w^{-1}(J)}\big)\cts\cong \Pi(\rho)_{w^{-1}(J)}$, as desired.

We note that in the setting of Lemma~\ref{lm:universal-tensor}, it follows from its proof that if each canonical map $\sigma_i\to \wh\sigma_i$ is injective, then so is the canonical map $\sigma_1 \otimes_\pi \cdots \otimes_\pi\sigma_r \to \wh\sigma_1 \whotimes \cdots \whotimes \wh\sigma_r$. Similarly, in the context of Proposition~\ref{prop:univ-unit-compl}, if $\sigma \to\wh\sigma$ is injective, then so is $(\pInd \sigma)\qpan \to (\pInd \wh\sigma)\cts$. By the injectivity assertion of Proposition~\ref{prop:breuil-emerton} and by construction we thus deduce that the canonical map $\PS(w,\sigma w) \to \PS(w,\sigma w)\sphat\,$ is injective.
\end{proof}

For $w \in W_\rho$ and $\sigma, \sigma' \in \Sigma_w$ with $\sigma \le \sigma'$, by Proposition~\ref{prop:univ-completion-ordinary} and \cite[\S 3.3]{fondamental} there exists a $G(\qp)$-equivariant embedding $\PS(w,\sigma w)\sphat\, \hookrightarrow \PS(w,\sigma'w)\sphat\,$ that is
unique up to scalar. As in that reference we can fix a compatible system of injections $i_{w,\sigma,\sigma'} : \PS(w,\sigma w)\sphat\, \into \PS(w,\sigma'w)\sphat\,$ (the choice of which won't matter) and obtain that 
\begin{equation}\label{eq:4}
\Pi(\rho)_{C_\rho, w^{-1}} \cong \ilim_{\sigma \in \Sigma_w}
\PS(w,\sigma w)\sphat\,
\end{equation}
in the notation of that reference (see Definition \ref{compareBH}).\\

We will need the following lemmas below.

\begin{lm}\label{lm:sigma-w}
Suppose that $w \in W_\rho$ and $\sigma \in \Sigma_w$.
\begin{enumerate}
\item We have $\ell(\sigma w) = \ell(\sigma) + \ell(w)$.
\item If $w \le \tau \le \sigma w$, then $\tau = \sigma' w$ with $\sigma' \le \sigma$ \(hence $\sigma' \in \Sigma_w$\).
\end{enumerate}
\end{lm}
\begin{proof}
Write $\sigma = s_{\alpha_1} \cdots s_{\alpha_r}$ with $\alpha_i \in \Delta$ pairwise orthogonal and $s_{\alpha_i} w \not \in W_\rho$ for all $i$. The latter condition implies that $w^{-1}(\alpha_i) > 0$ for all $i$. Note that $\sigma$ is of length $r$, sending each $\alpha_i$ to $-\alpha_i$ and preserving $\Phi^+-\{\alpha_1,\dots,\alpha_r\}$, where $\Phi^+ \subseteq X(T)$ (resp.\ $\Phi^- \subseteq X(T)$) denotes the positive (resp.\ negative) roots of $G=\GL_n$ with respect to $\o B$. It follows that $\sigma w$ and $w$ send precisely the same elements of $\Phi^+$ to $\Phi^-$, except that $w^{-1}(\alpha_1)$, \dots, $w^{-1}(\alpha_r)$ are sent to $\Phi^-$ by $\sigma w$ and to $\Phi^+$ by $w$. This implies (i).

For (ii) we induct on the length of $\sigma$. If $\sigma = 1$ there is nothing to show. Trivially we have $s_{\alpha_1} \sigma < \sigma$, hence by (i) we have $s_{\alpha_1} \sigma w < \sigma w$.

If $\tau^{-1}(\alpha_1) \in \Phi^+$, then $\tau < s_{\alpha_1} \tau$. By the lifting property of Coxeter groups, from the previous two inequalities and $\tau \le \sigma w$ we deduce that $\tau \le s_{\alpha_1} \sigma w$. By the induction hypothesis applied to $s_{\alpha_1} \sigma$ we deduce the claim.

If, on the other hand, $\tau^{-1}(\alpha_1) \in \Phi^-$, then $s_{\alpha_1} \tau < \tau$ and $w < s_{\alpha_1} w$ by (i). As $w \le \tau$, the lifting property of Coxeter groups gives $w \le s_{\alpha_1} \tau$. On the other hand, using that $s_{\alpha_1} \sigma w < \sigma w$ and $\tau \le \sigma w$, the lifting property gives $s_{\alpha_1} \tau \le s_{\alpha_1} \sigma w$. By applying the induction hypothesis to $s_{\alpha_1} \sigma$ we see that $s_{\alpha_1} \tau = \sigma' w$ with $\sigma' \le s_{\alpha_1} \sigma$. Hence $\tau = s_{\alpha_1}\sigma' w$. Finally, $\sigma' < s_{\alpha_1} \sigma'$ (as $\sigma'$ is a product of some $s_{\alpha_i}$ with $i > 1$) and the lifting property shows that $s_{\alpha_1}\sigma' \le \sigma$, as required.
\end{proof}

\begin{lm}\label{lm:socle-of-PS-hat-an}
Suppose that $w \in W_\rho$, $\sigma \in \Sigma_w$.
\begin{enumerate}
\item For any $\sigma' \le \sigma$ the representation $C(\sigma'w,\sigma'w)$ occurs with multiplicity one in both $\PS(w,\sigma w)$ and $(\PS(w,\sigma w)\sphat\,)\lan[\qp]$.
\item The socle of $(\PS(w,\sigma w)\sphat\,)\lan[\qp]$ is isomorphic to $C(w,w)$.
\end{enumerate}
\end{lm}
\begin{proof}
(i) We note that by Proposition~\ref{prop:univ-completion-ordinary} and its proof the representation $\PS(w,\sigma w)\sphat\,$ has a filtration with graded pieces $\big(\Ind_{B(\qp)}^{G(\qp)} \eta(\sigma''w,\sigma''w)\big)\cts$ with $\sigma'' \le \sigma$. Hence $(\PS(w,\sigma w)\sphat\,)\lan[\qp]$ has a filtration with graded pieces $\PS(\sigma''w,\sigma''w)$ with $\sigma'' \le \sigma$. Moreover $\PS(\sigma''w,\sigma''w)$ has irreducible constituents $C(v,\sigma''w)$ with $v \ge \sigma''w$. Suppose that $C(\sigma'w,\sigma'w)$ occurs in $\PS(\sigma''w,\sigma''w)$ for some $\sigma'' \le \sigma$. Let $\tau := \sigma'w$ and $\tau' := \sigma''w$. By Lemma~\ref{lm:princ-series}(ii) we see that $\tau \ge \tau'$ and $W(\tau)\tau' = W(\tau)\tau$.

We claim that $\tau = \tau'$. By Lemma~\ref{lm:sigma-w}(ii) we deduce from $\tau \ge \tau'$ that $\sigma' \ge \sigma''$, so $\sigma' = u \sigma''$ with $\ell(\sigma')=\ell(u)+\ell(\sigma'')$. Hence $\tau = u \tau'$ and by Lemma~\ref{lm:sigma-w}(i) we deduce that $s_\alpha \tau < \tau$ for any simple reflection $s_\alpha$ in the support of $u$. On the other hand, by above, $u = \tau (\tau')^{-1} \in W(\tau)$, so $s_\alpha \tau > \tau$ for any simple reflection $s_\alpha$ in the support of $u$. This shows that $u = 1$, so indeed $\tau = \tau'$.

As $C(\tau,\tau)$ occurs with multiplicity one in $\PS(\tau,\tau)$ we have established part (i) for the representation $(\PS(w,\sigma w)\sphat\,)\lan[\qp]$.

For the representation $\PS(w,\sigma w)$, we first claim that $C(\tau,\tau) \cong C(\tau,\sigma w)$. Indeed, by Lemma~\ref{lm:princ-series}(ii) this is equivalent to showing that $\sigma (\sigma')^{-1} \in W(\tau)$. This is true by Lemma~\ref{lm:sigma-w}(i), which shows that $\ell(\sigma(\sigma')^{-1})+\ell(\tau) = \ell(\sigma (\sigma')^{-1}\tau) = \ell(\sigma w)$. It follows that $C(\tau,\tau)$ occurs in $\PS(w,\sigma w)$. As $\PS(w,\sigma w)$ injects into $(\PS(w,\sigma w)\sphat\,)\lan[\qp]$ by Proposition~\ref{prop:univ-completion-ordinary} the proof is complete. (Alternatively we could check directly that $P_{w,\tau}(1) = 1$ (Kazhdan--Lusztig polynomial) using Jantzen's criterion \cite[\S8.7]{humphreys-bgg}.)

(ii) By the filtration mentioned in (i), if $C$ is any irreducible closed subrepresentation of $(\PS(w,\sigma w)\sphat\,)\lan[\qp]$, then it has to inject into $\PS(\sigma'w,\sigma'w)$ for some $\sigma' \le \sigma$, hence $C \cong C(\sigma'w,\sigma'w)$ for some $\sigma' \le \sigma$. Using part (i) and that $\PS(w,\sigma w)$ injects into $(\PS(w,\sigma w)\sphat\,)\lan[\qp]$ we deduce that $C$ injects into $\PS(w,\sigma w)$, hence $C \cong C(w,\sigma w)$. By Lemma~\ref{lm:princ-series}(ii) we have $C(w,\sigma w) \cong C(w,w)$. Finally, the socle of $(\PS(w,\sigma w)\sphat\,)\lan[\qp]$ is irreducible, as $C(w,w)$ occurs with multiplicity one by part (i).
\end{proof}

We need to understand better the ordinary representations $\Pi(\rho)_{C_\rho,w^{-1}}$ constructed in \cite[\S3]{fondamental}. To do this, we introduce an abstract framework. Suppose that $(I,\le)$ is a finite poset, that $X_i$ ($i \in I$) are objects of some abelian category $\cA$, and that we have a compatible system of injections $X_i \into X_j$ for any $i \le j$. We say that a subset $J \subseteq I$ is a \emph{lower subset} if $i_1 \le i_2$ in $I$ and $i_2 \in J$ imply $i_1 \in J$. Consider the following condition.

\begin{cond}\label{cond:submodule-intersection} 
For any non-empty lower subset $J$ having upper bound $b \in I$ and for any maximal element $m$ of $J$ we have $(\sum_{J - \{m\}} X_j) \cap X_m = \sum_{j < m} X_j$ inside $X_b$.
\end{cond}

For any lower subset $J$ we define $L_J := \ilim_{j\in J} X_j$. If $J_1 \subseteq J_2$ are two lower subsets, then we have a canonical map $L_{J_1} \to L_{J_2}$. If $J = \{i \in I : i \le n\}$ for some $n \in I$, then we write $L_{\le n}$ for $L_J$ and $L_{< n}$ for $L_{J-\{n\}}$. Note that $L_{\le n}\cong X_n$. Also note that $L_\varnothing = 0$.

\begin{lm}\label{lm:injection-colimit}
The map $L_{J'} \to L_J$ is injective for all pairs of lower subsets $J' \subseteq J$ if and only if Condition~\ref{cond:submodule-intersection} holds. If this holds, then for any lower subsets $I_1$, $I_2$ we have $L_{I_1 \cap I_2} = L_{I_1} \cap L_{I_2}$ inside $L_I$.
\end{lm}
\begin{proof}
We first observe that if $I_1$, $I_2$ are lower subsets, then $L_{I_1 \cup I_2} \cong L_{I_1} \oplus_{L_{I_1 \cap I_2}} L_{I_2}$ (write down inverse isomorphisms).

To prove ``$\Leftarrow$'', we induct on $\# J$, the case $J = \varnothing$ being trivial. We may assume that $J' \ne J$. Pick a maximal element $m$ of $J$ which is moreover such that $J' \subseteq J-\{m\}$ (note that $m$ exists since $J'$ is a lower subset and that $J-\{m\}$ is still a lower subset). 

Assume first that $m$ is {\it not} a maximum of $J$, then the maps $L_{J'} \to L_{J-\{m\}}$, $L_{< m} \to L_{J-\{m\}}$, and $L_{< m} \to L_{\le m}$ are injective by induction hypothesis (which can be applied to the latter since $L_{\le m}\subsetneq J$). Hence so is $L_{J'} \to L_{J-\{m\}} \oplus_{L_{< m}} L_{\le m}= L_J$ and we are done. 

Assume now that $m$ is the maximum of $J$, so $L_J = X_m$. We now fix $J$ and induct on $\# J'$, the case $J' = \varnothing$ being trivial. Let $n$ denote a maximal element of $J'$. If $n$ is a maximum, then $X_n = L_{J'} \to L_J = X_m$ is injective by assumption. Otherwise, $L_{J'} = L_{J'-\{n\}} \oplus_{L_{<n}} L_{\le n}$. By induction the maps
\begin{equation*}
\xymatrix{L_{J'-\{n\}} \ar[r] & L_{J} \\ L_{<n} \ar[r] \ar[u] & L_{\le n} \ar[u]}
\end{equation*}
are all injective as previously. The images of $L_{J'-\{n\}}$, $L_{<n}$, $L_{\le n}$ inside $L_J = X_m$ are equal to $\sum_{J'-\{n\}} X_j$, $\sum_{j<n} X_j$, $X_n$, respectively. Thus the map $L_{J'} = L_{J'-\{n\}} \oplus_{L_{<n}} L_{\le n} \to L_J$ is injective if and only if $(\sum_{J'-\{n\}} X_j) \cap X_n = \sum_{j<n} X_j$ inside $X_m$, which holds by Condition~\ref{cond:submodule-intersection}.

As a consequence we know that whenever $I_1$, $I_2$ are lower subsets, the map $L_{I_1 \cup I_2} = L_{I_1} \oplus_{L_{I_1 \cap I_2}} L_{I_2} \to L_I$ is injective, so $L_{I_1 \cap I_2} = L_{I_1} \cap L_{I_2}$ inside $L_I$.

To prove ``$\Rightarrow$'', we just apply the identity $L_{I_1 \cap I_2} = L_{I_1} \cap L_{I_2}$ inside $L_{I_1 \cup I_2}\subseteq L_I$ with $I_1 = J-\{m\}$ and $I_2 = \{i \in I : i \le m\}$ (noting that $I_1 \cup I_2 \subseteq \{i \in I : i \le b\}$).
\end{proof}

\begin{ex}\label{ex:mult-free}
Suppose that $C_i$ ($i \in I$) are simple objects in $\cA$ that are pairwise non-isomorphic and that $X_i$ is a finite length object that is multiplicity-free with Jordan--H\"older factors $\{ C_j : j \le i\}$ such that the submodule structure (i.e.\ the Alperin diagram) of $X_i$ is described by the partial order $\le$. Then Condition~\ref{cond:submodule-intersection} holds.
\end{ex}

The following lemma is in fact already tacitly used in \cite{fondamental} (and should have been proved there!).

\begin{coroll}\label{cor:injection-pi-rho-ord}
Fix $w \in W_\rho$, $\sigma \in \Sigma_w$ and let $J_\sigma := \{ \alpha \in \Delta : s_\alpha \in \supp(\sigma) \}$. Then the map $\Pi(\rho)_{w^{-1}(J_\sigma)} \to \Pi(\rho)_{C_\rho, w^{-1}}$ is injective \(with the notation of \cite[\S3.3]{fondamental}, see \eqref{eq:4}\).
\end{coroll}
\begin{proof}
We apply the above formalism with $\cA$ the abelian category of admissible continuous representations of $G$ on Banach spaces, $I = \Sigma_w$ with respect to $\le$, and $X_\sigma = \Pi(\rho)_{w^{-1}(J_\sigma)}$ for $\sigma \in \Sigma_w$. Recall that we picked a compatible system of injections between the $X_\sigma$. To verify Condition~\ref{cond:submodule-intersection} we fix $\sigma \in \Sigma_w$, playing the role of the upper bound $b \in I$. For each $\sigma' \le \sigma$ we can write $X_{\sigma'} \cong \big(\Ind_{P_{J_\sigma}(\qp)}^{G(\qp)} Y_{\sigma'}\big)\cts$, where $Y_{\sigma'} := \big(\Ind_{(M_{J_\sigma}
\cap P_{J_{\sigma'}})(\qp)}^{M_{J_\sigma}(\qp)} \wt\Pi(\rho)_{w^{-1}(J_{\sigma'})}\big)\cts$. The functor $\cF : Y \mapsto \big(\Ind_{P_{J_\sigma}(\qp)}^{G(\qp)} Y\big)\cts$ from submodules of $\wt\Pi(\rho)_{w^{-1}(J_{\sigma})}$ (see Step 2 of the proof of \cite[Prop.\ 3.3.3]{fondamental} for the notation) to submodules of $\Pi(\rho)_{w^{-1}(J_{\sigma})}$ respects addition and intersections (for example, by choosing a continuous section and rewriting $\cF(Y) \cong \C^0(P_{J_\sigma}(\qp)\backslash G(\qp), Y)$). Then Condition~\ref{cond:submodule-intersection} follows from the corresponding condition on the Levi subgroup $M_{J_\sigma}(\qp)$ by Example~\ref{ex:mult-free} and \cite[Rk.\ 3.3.4(ii)]{fondamental}.
\end{proof}

For any $n \in I$ we let $Q_n := L_{\le n}/L_{< n}$.

\begin{lm}\label{lm:filtration}
Suppose that Condition~\ref{cond:submodule-intersection} holds. Then for any lower subset $J$, the object $L_J$ has a filtration with graded pieces isomorphic to $Q_j$ \($j \in J$\).
\end{lm}
\begin{proof}
By induction it suffices to show that if $J$ is a lower subset and $m \in J$ is a maximal element, then $L_J/L_{J-\{m\}} \cong Q_m$. To see this, note that the natural map $L_{\le m}/L_{< m} \to L_J/L_{J-\{m\}}$ is surjective by construction and injective by Lemma~\ref{lm:injection-colimit}.
\end{proof}

\begin{coroll}\label{cor:filtration-pi-rho-ord}
For any $w \in W_\rho$ the representation $\Pi(\rho)_{C_\rho, w^{-1}}$ has a filtration with graded \ \ pieces \ \ $\big(\Ind_{B(\qp)}^{G(\qp)} \eta(\sigma w,\sigma w)\big)\cts$ \ \ with \ \ $\sigma \in \Sigma_w$. \ \ Moreover \ \ we \ \ have \ \ $\soc_{G(\qp)} (\Pi(\rho)_{C_\rho, w^{-1}})\lan[\qp] \cong C(w,w)$.
\end{coroll}
\begin{proof}
As in the proof of Corollary~\ref{cor:injection-pi-rho-ord} we put ourselves in the context of the above formalism. Then the first claim follows from Lemma~\ref{lm:filtration}. Hence $(\Pi(\rho)_{C_\rho, w^{-1}})\lan[\qp]$ has a filtration with graded pieces $\PS(\sigma w,\sigma w)$ with $\sigma\in \Sigma_w$. We deduce that if $C$ is any irreducible closed subrepresentation of $(\Pi(\rho)_{C_\rho, w^{-1}})\lan[\qp]$, then it has to inject into $\PS(\sigma w,\sigma w)$ for some $\sigma \in \Sigma_w$, hence $C \cong C(\sigma w,\sigma w)$ for some $\sigma \in \Sigma_w$. We claim that $C(\sigma w,\sigma w)$ occurs in $(\Pi(\rho)_{C_\rho, w^{-1}})\lan[\qp]$ with multiplicity one, or equivalently that it occurs in $\PS(\sigma' w,\sigma' w)$ for $\sigma' \in \Sigma_w$ only when $\sigma' = \sigma$. If $C(\sigma w,\sigma w)$ occurs in $\PS(\sigma' w,\sigma' w)$, then by Lemma~\ref{lm:princ-series}(ii) we deduce that $\sigma w \ge \sigma' w$, hence $\sigma \ge \sigma'$ by Lemma~\ref{lm:sigma-w}(ii). It now follows from the proof of Lemma~\ref{lm:socle-of-PS-hat-an}(i) that $\sigma = \sigma'$, proving the claim. Therefore $C$ has to be contained in the subrepresentation $(\PS(w,\sigma w)\sphat\,)\lan[\qp]$ of $(\Pi(\rho)_{C_\rho, w^{-1}})\lan[\qp]$, and the claim follows from Lemma~\ref{lm:socle-of-PS-hat-an}(ii).
\end{proof}

\begin{rk}
Suppose that all representations $\pi_\sigma := \big(\Ind_{B(\qp)}^{G(\qp)} \eta(\sigma w,\sigma w)\big)\cts$ with $\sigma \in \Sigma_w$ are (topologically) irreducible. Then the above results show that $\Pi(\rho)_{C_\rho, w^{-1}}$ is a multiplicity-free representation with Jordan--H\"older factors $\pi_\sigma$ ($\sigma \in \Sigma_w$) such that the submodule structure is described by the poset $(\Sigma_w, \le)$. In particular, this establishes the existence part of Conjecture~3.5.1 in \cite{fondamental} (in case all $\pi_\sigma$ are irreducible). Hauseux \cite{Hauseux-BH} recently established the uniqueness part (under the same assumption).
\end{rk}

\begin{prop}\label{prop:Qrho-ordinary}
With the above assumptions we have that $(\Pi(\rho)\ord)\lan[\qp]$ is isomorphic to a subrepresentation of $\Pi(\rho)\fs$.
\end{prop}
\begin{proof}
Fix any $w\in W_\rho$ and let $C := C(w,w)$. It suffices to show that the representation $(\Pi(\rho)_{C_\rho,w^{-1}})\lan[\qp]$ injects into $\Pi(\rho)\fs_C$. By Corollary~\ref{cor:filtration-pi-rho-ord} we know that $(\Pi(\rho)_{C_\rho,w^{-1}})\lan[\qp]$ has socle $C$. By Proposition~\ref{prop:univ-completion-ordinary} and Corollary~\ref{cor:injection-pi-rho-ord}, for each $\sigma \in \Sigma_w$ we have an injection $\PS(w,\sigma w) \into (\Pi(\rho)_{C_\rho,w^{-1}})\lan[\qp]$, which is unique up to scalars. We also recall that $(w,\sigma w) \in \Ws(\rho)$ for $\sigma \in \Sigma_w$ by Lemma~\ref{lm:walg}.
 
\emph{Step 1:} We show that $(\Pi(\rho)_{C_\rho,w^{-1}})\lan[\qp] = \sum_{\sigma \in \Sigma_w} \PS(w,\sigma w)$. We first consider $n = 2$ (with arbitrary $\rho$ satisfying our assumptions) and note that $(\PS(1,s_\alpha)\sphat\;)\lan[\qp] = \PS(1,s_\alpha) + \PS(1,1)$: by Corollary~\ref{cor:filtration-pi-rho-ord} the left-hand side has irreducible constituents $C(1,1)$, $C(s_\alpha,1)$, $C(s_\alpha,s_\alpha)$, each occurring with multiplicity one, and these all occur in the right-hand side.

For general $n$, by (\ref{eq:4}) it suffices to show $(\PS(w,\sigma w)\sphat\;)\lan[\qp] = \sum_{\sigma' \le \sigma}\PS(w,\sigma' w)$ for any fixed $\sigma \in \Sigma_w$. We define $J$, $P_J$, $\pi$ as in Proposition~\ref{prop:univ-completion-ordinary} and its proof, so that $\PS(w,\sigma w)\sphat\; \cong \big(\Ind_{P_J(\qp)}^{G(\qp)} \wh\pi\big)\cts$, and by Lemma~\ref{lm:an-vectors-in-parab-ind} we have $(\PS(w,\sigma w)\sphat\;)\lan[\qp] \cong \big(\Ind_{P_J(\qp)}^{G(\qp)} \wh\pi\lan[\qp] \big)\qpan$. By the proof of Proposition~\ref{prop:univ-completion-ordinary} and by Lemma~\ref{lm:analytic-vecs-tensor}, we deduce that $\wh\pi\lan[\qp]$ is an external completed tensor product of all $\psi_{w^{-1}(j)} \ve^{-(n-j)}$ for $\sigma(j) = j$ and $(\wh \pi_i)\lan[\qp]$ for $1 \le i \le r$, where $\pi_i$ is the representation~\eqref{eq:16} and $\wh \pi_i$ its universal unitary completion (note that $\wh \pi_i$ is admissible by Proposition \ref{prop:breuil-emerton}). From the previous paragraph we deduce that $(\wh\pi_i)\lan[\qp]$ is the sum of the subrepresentations $\pi$ and $\big(\Ind_{B_2(\qp)}^{\GL_2(\qp)} \psi_{w^{-1}(n_i)} \ve^{-(n-n_i)} \otimes\psi_{w^{-1}(n_{i+1})} \ve^{-(n-n_i-1)} \big)\qpan$. The exactness of locally analytic parabolic induction then implies the claim.

\emph{Step 2:} We show that $(\Pi(\rho)_{C_\rho,w^{-1}})\lan[\qp]$ injects into $\Pi(\rho)\fs_C$. By Step 1 it follows that the amalgam $A := \bigoplus_{C} \{ \PS(w,\sigma w) : \sigma \in \Sigma_w \}$ over the common socle $C$ surjects onto $(\Pi(\rho)_{C_\rho,w^{-1}})\lan[\qp]$. As $C$ occurs precisely once in $A$, namely in the socle, we see that $(\Pi(\rho)_{C_\rho,w^{-1}})\lan[\qp]$ is the unique quotient of $A$ that has socle $C$. We also have that $M(w,\sigma w):=M(w,(\sigma w)\F_\rho) = \PS(w,\sigma w)$ by Lemma \ref{lm:walg} and Remark \ref{rk:our-amalgam}(ii). Now, consider the composition $A = \bigoplus_C \{M(w,\sigma w): \sigma \in \Sigma_w \} \into M(\rho)_C \onto \Pi(\rho)\fs_C$. As $\Pi(\rho)\fs_C$ has socle $C$, we deduce by what we showed at the beginning of Step $2$ that it factors through a map $(\Pi(\rho)_{C_\rho,w^{-1}})\lan[\qp] \to\Pi(\rho)\fs_C$. By considering socles we see that it is injective.
\end{proof}

\begin{rk}
The proof shows, in particular, that $(\Pi(\rho)_{C_\rho,w^{-1}})\lan[\qp]$ for $w \in W_\rho$ can be described more explicitly as the unique quotient of $\bigoplus_{C(w,w)} \{ \PS(w,\sigma w) : \sigma \in \Sigma_w \}$ that has socle $C(w,w)$.
\end{rk}

\begin{prop}\label{prop:univ-completion-of-analytic-vecs}
For any $w \in W_\rho$ the unitary representation $\Pi(\rho)_{C_\rho,w^{-1}}$ is the universal unitary completion of $(\Pi(\rho)_{C_\rho,w^{-1}})\lan[\qp]$. Also, $\Pi(\rho)\ord$ is the universal unitary completion of $(\Pi(\rho)\ord)\lan[\qp]$.
\end{prop}
\begin{proof}
\emph{Step 1:} We show that for any $w \in W_\rho$ and $\sigma \in \Sigma_w$ the unitary representation $\PS(w,\sigma w)\sphat\,$ is the universal unitary completion of $(\PS(w,\sigma w)\sphat\,)\lan[\qp]$.

We need to show that if $\Pi$ is a unitary continuous representation of $G(\qp)$ on a Banach space, then any continuous $G(\qp)$-equivariant map $\theta : (\PS(w,\sigma w)\sphat\,)\lan[\qp]\to \Pi$ extends uniquely to a continuous $G(\qp)$-equivariant map $\PS(w,\sigma w)\sphat\, \to \Pi$. The uniqueness is clear by the density of locally analytic vectors. Let $i : \PS(w,\sigma w) \to \PS(w,\sigma w)\sphat\,$, $i' : \PS(w,\sigma w) \to (\PS(w,\sigma w)\sphat\,)\lan[\qp]$ and $j : (\PS(w,\sigma w)\sphat\,)\lan[\qp] \to \PS(w,\sigma w)\sphat\,$ denote the canonical maps (all of which are injective), so $i = j\circ i'$. From the definition of $\PS(w,\sigma w)\sphat\;$ the map $\theta \circ i'$ extends uniquely to a map $\theta' : \PS(w,\sigma w)\sphat\, \to \Pi$, i.e.\ $\theta' \circ i = \theta\circ i'$. It follows that $(\theta' \circ j - \theta) \circ i' = 0$, i.e.\ $\theta' \circ j - \theta$ factors through the cokernel of $i'$.

We claim that no irreducible constituent of the cokernel of $i'$ admits a $G(\qp)$-invariant $\O$-lattice. This claim easily implies that $\theta' \circ j = \theta$, completing the proof. Suppose now that $C$ is any irreducible constituent of $\coker(i')$ that admits a $G(\qp)$-invariant $\O$-lattice. By the proof of Lemma~\ref{lm:socle-of-PS-hat-an}(i) we know that $C \cong C(\tau,\sigma' w)$ for some $\sigma' \le \sigma$ and $\tau \ge \sigma' w$. More generally, suppose that $C(\tau,w')$ admits a $G(\qp)$-invariant $\O$-lattice for any $(\tau,w') \in W^2$ with $\tau \ge w'$. Then the necessary condition of Emerton (cf.\ the proof of \cite[Cor.\ 7.7]{socle1}) shows that $(-\tau\cdot (-\lambda))(t) \pi_{B,w'}(t) \in \O$ for all $t \in T^+$, where $\pi_{B,w'}:= \pi_{B,w'\F_\rho}$, $T^+ := \{ \diag(t_1,\dots,t_n) \in T(\qp) : |t_i t_{i+1}^{-1}|_{\qp} \le 1\ \forall i\}$. (We note that the dot action in \cite{socle1} is defined relative to $\o B$.) 
By equation \cite[(8.8)]{socle1} and the line following \cite[Rq.\ 8.7]{socle1} we deduce that $p^{\sum_{i=1}^j(h_{w'^{-1}(i)}-h_{\tau^{-1}(i)})} \in \O$ for all $1\le j \le n$, which is easily seen to be equivalent to $\tau(-h) \ge w'(-h)$ by the dominance order on $X(T)$ relative to $\o B$. But $-h$ is dominant with respect to $\o B$ and $\tau \ge w'$ by assumption, so $\tau(-h) \le w'(-h)$. It follows that $\tau(-h) = w'(-h)$. As the $h_i$ are distinct, we deduce that $\tau = w'$. For our constituent $C$ above this means $C \cong C(\sigma' w, \sigma' w)$. But $C$ does not occur in $\coker(i')$ by Lemma~\ref{lm:socle-of-PS-hat-an}(i), contradiction.

\emph{Step 2:} We deduce the result. It is completely formal to see that universal unitary completions commute with finite colimits on the additive category of continuous representations of a $p$-adic reductive group on locally convex vector spaces. (In fact, finite colimits exist in this category, since finite direct sums and cokernels exist.) Similarly, the functor of passing to locally $\qp$-analytic vectors commutes with finite colimits on the abelian category of admissible continuous representations of $G(\qp)$ by \cite[Thm.~7.1]{MR1990669}. Hence from \eqref{eq:4} we get $(\Pi(\rho)_{C_\rho,w^{-1}})\lan[\qp] \cong \ilim_{\sigma \in \Sigma_w} (\PS(w,\sigma w)\sphat\,)\lan[\qp]$ and by Step 1 we deduce
\begin{equation*}
\big((\Pi(\rho)_{C_\rho,w^{-1}})\lan[\qp]\big)\sphat\; \cong \ilim_{\sigma \in \Sigma_w}
 \PS(w,\sigma w)\sphat\; \cong \Pi(\rho)_{C_\rho,w^{-1}},
\end{equation*}
completing the proof of the first statement. By passing to a finite direct sum over $W_\rho$ we deduce the second statement.
\end{proof}

\subsection{Global applications}\label{sec:global-app-ord}

We give strong evidence to \cite[Conj.~4.2.2]{fondamental} in the crystabelline case (Theorem \ref{thm:main}).\\

We keep the global setup and notation of \S\ref{sec:global-app-crystab}, but now assume in addition that $p$ splits completely in $F$ (or equivalently $F^+$). We assume in the following that for each place $v | p$ of $F^+$ the representation $r_\tv = r|_{\Gal(\o {F_\tv}/F_\tv)}$ is as in \S\ref{sec:ordin-repr}, namely that it is upper-triangular as in \eqref{eq:rho}, satisfies Hypothesis~\ref{hyp:generic-crystabelline}, and we choose a representative (still denoted by) $r_\tv$ of $r_\tv$ that is a good conjugate.

\begin{prop}\label{prop:princ-series-to-cohom}
Assume Conjecture~\ref{conj:loc-an-socle}. If $(w\alg_\tv, w_\tv)_v \in \W$ such that $w\alg_\tv = w\alg_\tv(w_\tv)$ for all $v$, then restriction to the socle induces an isomorphism
\begin{multline*}
\Hom_{\u G(F_p^+)}\big(\whotimes_{v|p} \PS(w\alg_\tv,w_\tv)(\ve^{n-1}), \wh S(U^p,E)\lan[\qp][\m_r]\big)\\
 \congto \Hom_{\u G(F_p^+)}\big(\whotimes_{v|p} C(w\alg_\tv,w_\tv)(\ve^{n-1}), \wh S(U^p,E)\lan[\qp][\m_r]\big).
\end{multline*}
Moreover, any non-zero element of the left-hand side is injective.
\end{prop}
\begin{proof}
By Remark~\ref{rk:our-amalgam}(ii) this is a special case of Step 1 of the proof of Theorem~\ref{thm:adjunct-max-subrep-ps-coho}.
\end{proof}

\begin{prop}\label{prop:unique-extension-property}
Assume Conjecture~\ref{conj:loc-an-socle}. For each place $v|p$ of $F^+$ suppose that $w_\tv \in W_{r_\tv}$ and $\sigma_\tv,\sigma'_\tv \in \Sigma_{w_\tv}$ with $\sigma_\tv \le \sigma'_\tv$. Then the restriction map
\begin{multline}\label{eq:2}
\Hom_{\u G(F_p^+)}\big(\whotimes_{v|p} \PS(w_\tv,\sigma_\tv w_\tv)\sphat\, (\ve^{n-1}), \wh S(U^p,E)[\m_r]\big) \\
\congto \Hom_{\u G(F_p^+)}\big(\whotimes_{v|p} \PS(w_\tv,\sigma'_\tv w_\tv)\sphat\, (\ve^{n-1}), \wh S(U^p,E)[\m_r]\big) 
\end{multline}
induced by the injections $i_{w_\tv,\sigma_\tv,\sigma'_\tv}$ of \S\ref{sec:local-setup-ord} \(see just above \eqref{eq:4}\) is an isomorphism of finite-dimensional \ \ vector \ \ spaces. \ \ Moreover, \ \ any \ \ non-zero \ \ element \ \ of \ $\Hom_{\u G(F_p^+)}(\whotimes_{v|p} \PS(w_\tv,\sigma_\tv w_\tv)\sphat\, (\ve^{n-1}), \wh S(U^p,E)[\m_r])$ is injective.
\end{prop}
\begin{proof}
We will first check the last assertion, by passing to locally $\qp$-analytic vectors, using Lemma~\ref{lm:analytic-vecs-tensor} and \cite[Thm.~7.1]{MR1990669}. Suppose that for each $v|p$ we are given an irreducible constituent $C_\tv$ of $(\PS(w_\tv,\sigma_\tv w_\tv)\sphat\,)\lan[\qp]$. Recall from the proof of Lemma~\ref{lm:socle-of-PS-hat-an} that the representation $(\PS(w_\tv,\sigma_\tv w_\tv)\sphat\,)\lan[\qp]$ has a filtration with graded pieces $\PS(\sigma'_\tv w_\tv,\sigma'_\tv w_\tv)$ with $\sigma'_\tv \le \sigma_\tv$, hence by Lemma~\ref{lm:princ-series}(iii) we have $C_\tv \cong C(\tau_\tv,\sigma'_\tv w_\tv)$ for some $\tau_\tv \ge \sigma'_\tv w_\tv$. If $\whotimes_{v|p} C_\tv(\ve^{n-1})$ injects into $\wh S(U^p,E)\lan[\qp][\m_r]$, then by Conjecture~\ref{conj:loc-an-socle}, Lemma~\ref{lm:walg} and Lemma \ref{lm:sigma-w}(i) we get that $\tau_\tv \le w\alg_\tv(\sigma'_\tv w_\tv) = w_\tv \le \sigma'_\tv w_\tv$ and hence $\tau_\tv = w_\tv = \sigma'_\tv w_\tv$ and $\sigma'_\tv = 1$, i.e.\ $C_\tv \cong C(w_\tv,w_\tv)$. From Lemma~\ref{lm:socle-for-product-group} and Lemma \ref{lm:socle-of-PS-hat-an}(ii) we deduce that $\whotimes_{v|p} C_\tv$ is the $\u G(F_p^+)$-socle of $\whotimes_{v|p} (\PS(w_\tv,\sigma_\tv w_\tv)\sphat\,)\lan[\qp]$. Therefore, any non-zero element of $\Hom_{\u G(F_p^+)}(\whotimes_{v|p}\PS(w_\tv,\sigma_\tv w_\tv)\sphat\, (\ve^{n-1}), \wh S(U^p,E)[\m_r])$ is injective.

Since $\whotimes_{v|p} C(w_\tv,w_\tv)$ is the $\u G(F_p^+)$-socle of $\whotimes_{v|p}
(\PS(w_\tv,\tau_\tv w_\tv)\sphat\,)\lan[\qp]$ for $\tau_\tv \in \{\sigma'_\tv,\sigma_\tv\}$, occurring
as constituent with multiplicity one (by Lemma~\ref{lm:socle-of-PS-hat-an}(i)), it follows that the map~\eqref{eq:2} is injective.
 
To complete the proof it suffices to show that the two sides of~\eqref{eq:2} have the same finite dimension. We note that
\begin{multline}\label{eq:3}
\Hom_{\u G(F_p^+)}\big(\whotimes_{v|p} \PS(w_\tv,\sigma_\tv w_\tv)\sphat\, (\ve^{n-1}), \wh S(U^p,E)[\m_r]\big) \\
 \cong \Hom_{\u G(F_p^+)}\big(\whotimes_{v|p} \PS(w_\tv,\sigma_\tv w_\tv)(\ve^{n-1}), \wh S(U^p,E)\lan[\qp][\m_r]\big)\\
\cong \Hom_{\u G(F_p^+)}\big(\whotimes_{v|p} C(w_\tv,\sigma_\tv w_\tv)(\ve^{n-1}), \wh S(U^p,E)\lan[\qp][\m_r]\big)
\end{multline}
by Lemma~\ref{lm:universal-tensor}, Proposition~\ref{prop:princ-series-to-cohom}, and Lemma~\ref{lm:walg}. As $C(w_\tv,\sigma_\tv w_\tv) \cong C(w_\tv,w_\tv)$ is independent of $\sigma_\tv$ by Lemma~\ref{lm:princ-series}(ii) and the vector space \eqref{eq:3} is finite-dimensional by Lemma~\ref{lm:fd} we complete the proof.
\end{proof}

The theorem that follows gives evidence for \cite[Conj.~4.2.2]{fondamental} (corrected as in \S\ref{sec:errata} below) in the crystabelline case.

\begin{thm}\label{thm:main}
Assume Conjecture~\ref{conj:loc-an-socle}. Then there exists an injection of admissible continuous representations $\whotimes_{v|p} \Pi(r_\tv)\ord(\ve^{n-1}) \into \wh S(U^p,E)[\m_r]$. More precisely, for any $w = (w_\tv)_v \in \prod_{v|p} W_{r_\tv}$ let
$$n_w := \dim_E \Hom_{\u G(F_p^+)}\big(\whotimes_{v|p} C(w_\tv,w_\tv)(\ve^{n-1}),\wh S(U^p,E)\lan[\qp][\m_r]\big) \in \Z_{> 0}.$$
Then we have an injection of admissible continuous representations
\begin{equation}\label{eq:28}
\bigoplus_{w = (w_\tv)_v} \big(\whotimes_{v|p} \Pi(r_\tv)_{C_{r_\tv}, w_\tv^{-1}}(\ve^{n-1})\big)^{\oplus n_w} \into \wh S(U^p,E)[\m_r].
\end{equation}
\end{thm}
\begin{proof}
By assumption we have an injection 
\begin{equation}
\bigoplus_{w = (w_\tv)_v} \big(\whotimes_{v|p} C(w_\tv,w_\tv)(\ve^{n-1})\big)^{\oplus n_w} \to \sanmr.\label{eq:27}
\end{equation}
Fix any $w = (w_\tv)_v \in \prod_{v|p} W_{r_\tv}$. By the isomorphism~\eqref{eq:4} and Lemma~\ref{lm:compl-tensor-finite-colimits} we have
\begin{multline*}
\Hom_{\u G(F_p^+)}\big(\whotimes_{v|p} \Pi(r_\tv)_{C_{r_\tv}, w_\tv^{-1}}(\ve^{n-1}), \wh S(U^p,E)[\m_r]\big) \\
\congto
\Hom_{\u G(F_p^+)}\big(\ilim_{\sigma_\tv \in \Sigma_{w_\tv}} \whotimes_{v|p} \PS(w_\tv,\sigma_\tv w_\tv)\sphat\, (\ve^{n-1}), \wh S(U^p,E)[\m_r]\big) \\
\cong \plim_{\sigma_\tv \in \Sigma_{w_\tv}} \Hom_{\u G(F_p^+)}\big(\whotimes_{v|p} \PS(w_\tv,\sigma_\tv w_\tv)\sphat\, (\ve^{n-1}), \wh S(U^p,E)[\m_r]\big).
\end{multline*}
By Proposition~\ref{prop:unique-extension-property} the projective limit is isomorphic to the final term where $\sigma_\tv = 1$ for all $v|p$ and hence by \eqref{eq:3} it is further isomorphic to
$$\Hom_{\u G(F_p^+)}\big(\whotimes_{v|p} C(w_\tv,w_\tv)(\ve^{n-1}), \wh S(U^p,E)\lan[\qp][\m_r]\big).$$
Thus we can extend the map~\eqref{eq:27} uniquely to a map as in~\eqref{eq:28}. The extended map is injective by the last statement of Corollary~\ref{cor:filtration-pi-rho-ord} (using Lemma~\ref{lm:socle-for-product-group}) and the injectivity of~\eqref{eq:27}.
\end{proof}

\begin{rk}
Alternatively we could prove Theorem \ref{thm:main} using Theorem~\ref{thm:adjunct-max-subrep-ps-coho}, Proposition~\ref{prop:Qrho-ordinary}, and Proposition~\ref{prop:univ-completion-of-analytic-vecs}. We also recall that Conjecture \ref{conj:loc-an-socle} is known in many cases (see Theorem \ref{BHS3}).
\end{rk}

\begin{rk}
  When $n = 3$, some cases of this theorem were claimed in an unpublished preprint \cite{BC-ord}.
\end{rk}

\section{Errata for \cite{fondamental}}\label{sec:errata}

The definition of $\Pi\ord$ just above \cite[Conj.\ 4.2.2]{fondamental} should be replaced by the following definition: $\Pi\ord$ is the closure (in the admissible continuous representation $\Pi$) of the sum of all its finite length closed subrepresentations with all irreducible constituents being constituents of unitary continuous principal series. Then \cite[Conj.\ 4.2.2]{fondamental} can be stated {\it verbatim}, and implies in particular that $\wh S(U^p,E)[\p^\Sigma]\ord$ should be of finite length. Note that, due to the closure process, it is not clear {\it a priori} that all irreducible constituents of $\wh S(U^p,E)[\p^\Sigma]\ord$ are still constituents of unitary continuous principal series.

The proof of \cite[Thm.\ 4.4.8]{fondamental} is too sketchy and moreover the representation $\wh S(U^p,E)_{\p^\Sigma}$ at the end of the proof is not a Banach space (it is not necessarily complete, as it is just some localisation), hence one cannot apply \cite[Cor.\ 4.3.11]{fondamental} to it. One can fix our proof of \cite[Thm.\ 4.4.8]{fondamental} (by working instead with the localisation at a maximal ideal of the complete integral Hecke algebra as in \cite[\S5.2]{emerton-local-global}), but in any case this result is now a special case of Theorem \ref{thm:main} (together with Theorem \ref{BHS3}).

Finally, due to the above comment on the (corrected) definition of $\wh S(U^p,E)[\p^\Sigma]\ord$, \cite[Rk.~4.4.9(a)]{fondamental} should be ignored.

\bibliography{BH-conj}

\newcommand{\etalchar}[1]{$^{#1}$}
\begin{thebibliography}{CEG{\etalchar{+}}18}

\bibitem[BC14]{BC-ord}
John Bergdall and Przemyslaw Chojecki.
\newblock Ordinary representations and companion points for {$U(3)$} in the
  indecomposable case.
\newblock Preprint, 2014.

\bibitem[BC18]{BC-adjunction}
John Bergdall and Przemyslaw Chojecki.
\newblock {An adjunction formula for the {E}merton-{J}acquet functor}.
\newblock {\em Israel J. Math.}, 223:1--52, 2018.

\bibitem[BE10]{BE}
Christophe Breuil and Matthew Emerton.
\newblock Repr\'esentations {$p$}-adiques ordinaires de {${\rm GL}_2(\mathbb
  {Q}_p)$} et compatibilit\'e local-global.
\newblock {\em Ast\'erisque}, 331:255--315, 2010.

\bibitem[BGR84]{MR746961}
S.~Bosch, U.~G\"untzer, and R.~Remmert.
\newblock {\em Non-{A}rchimedean analysis}, volume 261 of {\em Grundlehren der
  Mathematischen Wissenschaften [Fundamental Principles of Mathematical
  Sciences]}.
\newblock Springer-Verlag, Berlin, 1984.
\newblock A systematic approach to rigid analytic geometry.

\bibitem[BH15]{fondamental}
Christophe Breuil and Florian Herzig.
\newblock Ordinary representations of {$G(\mathbb{Q}_p)$} and fundamental
  algebraic representations.
\newblock {\em Duke Math. J.}, 164:1271--1352, 2015.

\bibitem[BHS17a]{BHS3}
Christophe Breuil, Eugen Hellmann, and Benjamin Schraen.
\newblock A local model for the trianguline variety and applications.
\newblock Preprint, 2017.

\bibitem[BHS17b]{BHS2}
Christophe Breuil, Eugen Hellmann, and Benjamin Schraen.
\newblock Smoothness and classicality on eigenvarieties.
\newblock {\em Invent. Math.}, 209:197--274, 2017.

\bibitem[Brea]{breuil-ext1}
Christophe Breuil.
\newblock {$\Ext^1$} localement analytique et compatibilit\'e local-global.
\newblock To appear in Amer. J. Math.

\bibitem[Breb]{Breuil-peking}
Christophe Breuil.
\newblock Towards the locally analytic socle for {${\rm GL}_n$}.
\newblock Notes of course given at BICMR in Beijing, 2013.

\bibitem[Bre15]{socle2}
Christophe Breuil.
\newblock Vers le socle localement analytique pour {${\rm GL}_n$} {II}.
\newblock {\em Math. Ann.}, 361:741--785, 2015.

\bibitem[Bre16]{socle1}
Christophe Breuil.
\newblock Socle localement analytique {I}.
\newblock {\em Ann. Inst. Fourier (Grenoble)}, 66:633--685, 2016.

\bibitem[CDP14]{CDP}
Pierre Colmez, Gabriel Dospinescu, and Vytautas Pa{\v s}k{\= u}nas.
\newblock The $p$-adic local {L}anglands correspondence for {$\GL_2(\qp)$}.
\newblock {\em Camb. J. Math.}, 2:1--47, 2014.

\bibitem[CEG{\etalchar{+}}16]{MR3529394}
Ana Caraiani, Matthew Emerton, Toby Gee, David Geraghty, Vytautas Pa{\v s}k{\=
  u}nas, and Sug~Woo Shin.
\newblock Patching and the {$p$}-adic local {L}anglands correspondence.
\newblock {\em Camb. J. Math.}, 4:197--287, 2016.

\bibitem[CEG{\etalchar{+}}18]{CEGGPS2}
Ana Caraiani, Matthew Emerton, Toby Gee, David Geraghty, Vytautas Pa{\v s}k{\=
  u}nas, and Sug~Woo Shin.
\newblock Patching and the $p$-adic {Langlands} program for
  $\mathrm{GL}_2(\mathbb{Q}_p)$.
\newblock {\em Compos. Math.}, 154(3):503--548, 2018.

\bibitem[CHT08]{CHT}
Laurent Clozel, Michael Harris, and Richard Taylor.
\newblock Automorphy for some {$\ell$}-adic lifts of automorphic mod {$\ell$}
  {Galois} representations.
\newblock {\em Publ. Math. Inst. Hautes \'Etudes Sci.}, 108:1--181, 2008.

\bibitem[Col10]{colmez-une-variable}
Pierre Colmez.
\newblock Fonctions d'une variable {$p$}-adique.
\newblock {\em Ast\'erisque}, 330:13--59, 2010.

\bibitem[Con99]{conrad-rigid}
Brian Conrad.
\newblock Irreducible components of rigid spaces.
\newblock {\em Ann. Inst. Fourier (Grenoble)}, 49:473--541, 1999.

\bibitem[Din17]{Ding}
Yiwen Ding.
\newblock Formes modulaires {$p$}-adiques sur les courbes de {S}himura
  unitaires et compatibilit\'{e} local-global.
\newblock {\em M\'{e}m. Soc. Math. Fr. (N.S.)}, (155):viii+245, 2017.

\bibitem[dS16]{deshalit-mahler}
Ehud de~Shalit.
\newblock Mahler bases and elementary {$p$}-adic analysis.
\newblock {\em J. Th\'eor. Nombres Bordeaux}, 28:597--620, 2016.

\bibitem[Eme]{jacquet2}
Matthew Emerton.
\newblock Jacquet modules of locally analytic representations of {$p$}-adic
  reductive groups. {II}. {T}he relation to parabolic induction.
\newblock To appear in J. Inst. Math. Jussieu.

\bibitem[Eme05]{MR2181093}
Matthew Emerton.
\newblock {$p$}-adic {$L$}-functions and unitary completions of representations
  of {$p$}-adic reductive groups.
\newblock {\em Duke Math. J.}, 130:353--392, 2005.

\bibitem[Eme06a]{jacquet1}
Matthew Emerton.
\newblock Jacquet modules of locally analytic representations of {$p$}-adic
  reductive groups. {I}. {C}onstruction and first properties.
\newblock {\em Ann. Sci. \'Ecole Norm. Sup. (4)}, 39:775--839, 2006.

\bibitem[Eme06b]{Em1}
Matthew Emerton.
\newblock On the interpolation of systems of eigenvalues attached to
  automorphic {Hecke} eigenforms.
\newblock {\em Invent. Math.}, 164, 2006.

\bibitem[Eme07]{MR2392361}
Matthew Emerton.
\newblock Locally analytic representation theory of {$p$}-adic reductive
  groups: a summary of some recent developments.
\newblock In {\em {$L$}-functions and {G}alois representations}, volume 320 of
  {\em London Math. Soc. Lecture Note Ser.}, pages 407--437. Cambridge Univ.
  Press, Cambridge, 2007.

\bibitem[Eme11]{emerton-local-global}
M.~Emerton.
\newblock Local-global compatibility in the $p$-adic {L}anglands programm for
  $\mathrm{GL}_{2/\Q}$.
\newblock Preprint, 2011.

\bibitem[Eme17]{locallyanalytic-memoir}
Matthew Emerton.
\newblock Locally analytic vectors in representations of locally {$p$}-adic
  analytic groups.
\newblock {\em Mem. Amer. Math. Soc.}, 248, 2017.

\bibitem[FdL99]{feaux-diss}
Christian~Tobias F\'eaux~de Lacroix.
\newblock Einige {R}esultate \"uber die topologischen {D}arstellungen
  {$p$}-adischer {L}iegruppen auf unendlich dimensionalen {V}ektorr\"aumen
  \"uber einem {$p$}-adischen {K}\"orper.
\newblock In {\em Schriftenreihe des {M}athematischen {I}nstituts der
  {U}niversit\"at {M}\"unster. 3. {S}erie, {H}eft 23}, volume~23 of {\em
  Schriftenreihe Math. Inst. Univ. M\"unster 3. Ser.}, pages x+111. Univ.
  M\"unster, Math. Inst., M\"unster, 1999.

\bibitem[Hau]{Hauseux-BH}
Julien Hauseux.
\newblock Sur une conjecture de {B}reuil-{H}erzig.
\newblock To appear in J. Reine Angew. Math.

\bibitem[Hum08]{humphreys-bgg}
James~E. Humphreys.
\newblock {\em Representations of semisimple {L}ie algebras in the {BGG}
  category {$\mathcal{O}$}}, volume~94 of {\em Graduate Studies in
  Mathematics}.
\newblock American Mathematical Society, Providence, RI, 2008.

\bibitem[OS15]{OS}
Sascha Orlik and Matthias Strauch.
\newblock On {J}ordan-{H}\"older series of some locally analytic
  representations.
\newblock {\em J. Amer. Math. Soc.}, 28:99--157, 2015.

\bibitem[Sch02]{nfa}
Peter Schneider.
\newblock {\em Nonarchimedean functional analysis}.
\newblock Springer Monographs in Mathematics. Springer-Verlag, Berlin, 2002.

\bibitem[Sch11]{schraen-gl3}
Benjamin Schraen.
\newblock Repr\'esentations localement analytiques de {${\rm GL}_3(\mathbb
  {Q}_p)$}.
\newblock {\em Ann. Sci. \'Ec. Norm. Sup\'er. (4)}, 44:43--145, 2011.

\bibitem[ST02a]{MR1900706}
P.~Schneider and J.~Teitelbaum.
\newblock Banach space representations and {I}wasawa theory.
\newblock {\em Israel J. Math.}, 127:359--380, 2002.

\bibitem[ST02b]{MR1887640}
Peter Schneider and Jeremy Teitelbaum.
\newblock Locally analytic distributions and {$p$}-adic representation theory,
  with applications to {${\rm GL}_2$}.
\newblock {\em J. Amer. Math. Soc.}, 15:443--468, 2002.

\bibitem[ST03]{MR1990669}
Peter Schneider and Jeremy Teitelbaum.
\newblock Algebras of {$p$}-adic distributions and admissible representations.
\newblock {\em Invent. Math.}, 153:145--196, 2003.

\bibitem[ST05]{ST-duality}
Peter Schneider and Jeremy Teitelbaum.
\newblock Duality for admissible locally analytic representations.
\newblock {\em Represent. Theory}, 9:297--326, 2005.

\bibitem[Tho12]{Thor}
Jack Thorne.
\newblock On the automorphy of {$\ell$}-adic {Galois} representations with
  small residual image.
\newblock {\em J. Inst. Math. Jussieu}, 11:855--920, 2012.

\end{thebibliography}
\bibliographystyle{alpha}

\end{document}